\documentclass[11pt,letterpaper]{amsart}

\usepackage{graphicx}
\usepackage{tabularx}
\usepackage{amssymb}
\usepackage{array}
\usepackage{amsthm}
\usepackage{mathtools}
\usepackage[hidelinks]{hyperref}
\usepackage{verbatim}
\usepackage{hyperref}
\usepackage{fullpage}
\usepackage{enumitem}
\usepackage{dsfont}
\usepackage{setspace}
\usepackage{verbatim}
\usepackage{stackengine}
\usepackage{tikz-cd}
\usepackage{cleveref}
\usepackage{amsmath}
\usepackage[backend=biber,style=alphabetic]{biblatex}
\usetikzlibrary{positioning}
\newtheorem{theorem}{Theorem}[section]

\newtheorem{lemma}[theorem]{Lemma}
\newtheorem{corollary}[theorem]{Corollary}

\theoremstyle{definition}
\newtheorem{definition}[theorem]{Definition}
\newtheorem{question}[theorem]{Question}
\newtheorem{observation}[theorem]{Observation}

\newtheorem{definition*}[theorem]{Definition*}
\newtheorem{lemma*}[theorem]{Lemma*}
\newtheorem{corollary*}[theorem]{Corollary*}

\newtheorem{question5.7}[theorem]{Question 5.7}
\newtheorem{remark}[theorem]{Remark}

\newtheoremstyle{quest}
  {}
  {}
  {}
  {0pt}
  {\bfseries}
  {}
  { }
  {\thmname{#1}\thmnumber{ #2}\textnormal{\thmnote{ (#3)}}}
\theoremstyle{quest}
\newtheorem*{question-non}{Question}

\numberwithin{equation}{subsection}
\title{On $\mathrm{BPI}$ in Symmetric Extensions Part 1}
\author{Brian Ransom}
\begin{document}
\maketitle
\begin{abstract}
Historically, proofs of $\mathrm{BPI}$ in models without choice have relied on a contradiction framework that was introduced by Halpern in \cite{halpern1962contributions}. We introduce the filter extension property for permutation models and symmetric extensions, which formalizes the na\"ive approach to extend arbitrary filters to ultrafilters by repeatedly extending filters by minimal increments. We use this framework to give the first direct proof of $\mathrm{BPI}$ in the generalized Cohen model $N(I,Q)$ -- a model that adds a Dedekind-finite set of mutually $Q$-generic filters over a ground model $M\vDash\mathrm{ZFC}$. In the case that the index set $I$ is large, we adapt Harrington's proof of the Halpern-L\"auchli theorem to prove the result. We then extend the results from \cite{karagila2020have} to show that $I$ can be assumed to be large without loss of generality. The approach given by Harrington's proof is essentially dynamical, and we show that this technique can be used in permutation models to reprove a direction of Blass' theorem from \cite{blass1986prime}: that a dynamical condition called the Ramsey property is sufficient for $\mathrm{BPI}$ to hold in a permutation model. We then introduce a dynamical generalization of the Ramsey property called the virtual Ramsey property, which abstracts core features of our adaptation of Harrington's proof, and we prove that the virtual Ramsey property is sufficient for $\mathrm{BPI}$ to hold in a symmetric extension.
\end{abstract}

\section{Introduction}
A central tension in the study of infinite objects is our limited ability to describe them. Over $\mathrm{ZF}$, non-effective axioms have shown broad success in addressing this tension. This raises the question: to what extent can the different tools that augment our ability to describe infinite objects, or processes, describe one another? In this paper, we revisit the relationship between the axiom of choice ($\mathrm{AC}$) and the ultrafilter lemma -- that every filter, on every set, can be extended to an ultrafilter. The latter theorem has many known equivalents over $\mathrm{ZF}$, two of which are the compactness theorem for first-order logic, and the Boolean prime ideal theorem ($\mathrm{BPI}$). By historical convention, we use the label $\mathrm{BPI}$ throughout the paper.

It was first shown in \cite{halpern1971boolean} that $\mathrm{BPI}$ is strictly weaker than $\mathrm{AC}$; here it was shown in particular that $\mathrm{BPI}$ holds in the Cohen model, and more generally that $\mathrm{BPI}$ holds in any model of a certain theory that abstracts properties of the Cohen model. Informally, $\mathrm{BPI}$ expresses the ability to make arbitrarily many \textit{finitely bounded} choices that are altogether \textit{locally compatible} with one another. Put this way, the result from \cite{halpern1971boolean} may seem surprising: locally compatible choices do not need to be made ``one at a time" via some wellordering by $\mathrm{AC}$. Put another way, the result from \cite{halpern1971boolean} may seem expected: one can get the sense that first-order compactness should not be able to wellorder an arbitrary set, since we can only exert local control over the construction. This tension captures much of the typical intrigue and difficulty of questions in the study of fragments of $\mathrm{AC}$.

Historically, it has taken significant effort to prove $\mathrm{BPI}$ in models without $\mathrm{AC}$. In his thesis \cite{halpern1962contributions}, Halpern gave the first proof that $\mathrm{BPI}$ was independent of $\mathrm{AC}$ over $\mathrm{ZFA}$, by showing that $\mathrm{BPI}$ holds in the ordered Mostowski model of $\mathrm{ZFA}$.\footnote{$\mathrm{ZFA}$ denotes the theory $\mathrm{ZF}$ with the existence of ``atoms;"  see subsection 2.1 of this paper or \cite{jech2008axiom} for more details.}$^,$\footnote{This portion of his thesis was later rewritten as the paper \cite{halpern1964independence}.} Halpern's argument employs a contradiction framework that appears, to the author's knowledge, in every other proof of $\mathrm{BPI}$ in a choiceless model to date, including both permutation models of $\mathrm{ZFA}$ and symmetric extensions of $\mathrm{ZF}$. In this framework, one extends an arbitrary hereditarily symmetric filter $F$, on an arbitrary hereditarily symmetric set $x$, to a maximal hereditarily symmetric filter $U\supseteq F$. By assuming that $U$ is not an ultrafilter, a combinatorial argument is used to produce the contradiction that $F$ must not have had the finite intersection property. The nature of this combinatorial argument depends on the context of the model.\footnote{Pincus made these exact observations regarding the known proofs of $\mathrm{BPI}$ in models without choice, in \cite{pincus1997dense}.}

In \cite{halpern1971boolean}, this contradiction is derived as an application of the (original) version of the Halpern-L\"auchli theorem from \cite{halpern1966partition}. In \cite{repicky2015proof}, a new proof of $\mathrm{BPI}$ in the Cohen model is given, where the application of the Halpern-L\"auchli theorem in \cite{halpern1971boolean} is replaced with a simpler finite induction argument to produce the desired contradiction; this version of the proof remains the most efficient. Recently, in \cite{stefanovic2023alternatives}, the proof from \cite{halpern1971boolean} is generalized to prove $\mathrm{BPI}$ in the model $N(I,Q)$, called the generalized Cohen model, which adds a Dedekind-finite set of mutually $Q$-generic filters. A new proof of the result that $N(I,Q)\vDash\mathrm{BPI}$ is one of the main goals of this paper. Notably, in the case of permutation models, Blass used Halpern's contradiction framework in \cite{blass1986prime} to prove $\mathrm{BPI}$ from the assumption of a dynamical condition called the Ramsey property (moreover, he proved the two to be equivalent).

While Halpern's contradiction framework has been widely influential, the underlying forcing used to define symmetric extensions has often made its application quite technical.\footnote{Comments along these lines are made, for instance, in \cite{pincus1977addingdep} and in \cite{jech2008axiom}.} To simplify these arguments would be of inherent interest, but it would also help address several contemporary questions relating to $\mathrm{BPI}$ in choiceless models. The question of preserving $\mathrm{BPI}$ in iterated symmetric extensions, for instance, is raised in \cite{karagila2019iterated}. This question is related to longstanding open conjectures from \cite{pincus1977adding}: that $\mathrm{BPI}$ holds in a candidate model for the independence of $\mathrm{AC}$ from $\mathrm{ZF}+\mathrm{BPI}+\mathrm{DC}$ (an independence which Pincus later proved in \cite{pincus1977addingdep}), and that a certain related variant of the Halpern-L\"auchli theorem holds. Questions like these seem to have remained open because generalizing the existing techniques would be exceedingly complicated.

With this in mind, we can state the main theses of this paper:
\begin{enumerate}
    \item[$\mathrm{T}1$] that it is natural and clarifying to prove $\mathrm{BPI}$ in choiceless models by direct arguments;
    \item[$\mathrm{T}2$] that the methodology from Harrington's proof of the Halpern-L\"auchli theorem can be adapted to give these direct proofs that $N(I,Q)\vDash\mathrm{BPI}$, when $I$ is large;
    \item[$\mathrm{T3}$] that the theorems from \cite{karagila2020have} can be extended to show that $I$ can be assumed to be large in $N(I,Q)$ without loss of generality;
    \item[$\mathrm{T}4$] that our adaptation of Harrington's proof offers a dynamical perspective on $\mathrm{BPI}$ in $N(I,Q)$, which can be expressed in terms of a generalization of the Ramsey property that implies $\mathrm{BPI}$ in symmetric extensions.
\end{enumerate}
Roughly, the theses $\mathrm{T}1,\mathrm{T}2,\mathrm{T}3,$ and $\mathrm{T}4$ are respectively addressed in sections 3, 4, 5, and 6.

\subsection{Outline}
In section 2, we state the basic definitions for permutation models of $\mathrm{ZFA}$, symmetric extensions of $\mathrm{ZF}$, the notion of ``infinite" that we will use without choice, and the generalized Cohen model $N(I,Q)$ from \cite{stefanovic2023alternatives}. We aim to provide a complete enough background that only knowledge of forcing is assumed to read the paper. Knowledge of the generalized Cohen model is not needed until section 4.

In section 3, we consider a na\"ive approach to directly prove $\mathrm{BPI}$ in models without choice, by asking how one would coarsen the standard proof of $\mathrm{BPI}$ in a model of $\mathrm{ZFC(A)}$. Our answer to this question is to introduce the filter extension property (for permutation models or symmetric systems). The filter extension property formalizes the notion that one can extend symmetric filters on symmetric sets by minimal increments, in a way that is stable under the iteration of extensions that result in ultrafilters. Using a theorem from \cite{blass1986prime}, we show that $\mathrm{BPI}$ in a permutation model or symmetric extension also implies the respective form of the filter extension property. Accordingly, the known proofs of $\mathrm{BPI}$ that use Halpern's contradiction framework can be rephrased as contradiction proofs of the filter extension property. The approach given in this paper diverges by considering ways to directly prove the filter extension property.

In section 4, we directly prove the filter extension property for the generalized Cohen model $N(I,Q)$, in the case that $I$ is assumed to be sufficiently large. We do so by adapting Harrington's proof of the Halpern-L\"auchli theorem. In particular, we show that the ability to extend a hereditarily symmetric name $\dot{F}$ for a filter on a hereditarily symmetric set $\dot{x}$ can be coordinated by appealing to a (not necessarily symmetric) name $\dot{U}$ for an ultrafilter on $\dot{x}$ that extends $\dot{F}$. This establishes a mechanism by which the ability to extend filters to ultrafilters in the symmetric extension $N(I,Q)$ can be ``inherited" from the ability to extend filters to ultrafilters in the outer model $M[G]\vDash\mathrm{ZFC}$.

Harrington's proof is celebrated for its conceptual simplicity. For the reader's benefit, we state the relevant version of the Halpern-L\"auchli theorem in Appendix $\mathrm{A}$, and broadly sketch Harrington's proof in Appendix $\mathrm{B}$. Our sketch in Appendix $\mathrm{B}$ is written to emphasize certain features that will influence our adaptation of the proof, and is based on the exposition of Harrington's proof given in \cite{todorchevich1995some}. If the reader is not already familiar with Harrington's proof, we recommend that they skim Appendix $\mathrm{B}$ before reading section 4, then return to fully read Appendices $\mathrm{A}$ and $\mathrm{B}$ afterward. While Appendices $\mathrm{A}$ and $\mathrm{B}$ can be skipped entirely if the reader wishes, they provide helpful context and motivation for our main arguments in section 4.

In \cite{karagila2020have}, Karagila and Schlicht essentially show that the index $I$ in the Cohen model $N(I,2^{<\omega})$ can be taken to be arbitrarily large, without loss of generality. In question 5.5 from \cite{karagila2020have}, they ask whether (up to a change in notation) the results in their paper also hold for generalized Cohen model $N(I,Q)$, and the speculate that this should be the case. In section 5, we answer a case of this question positively, and essentially show that the index $I$ in $N(I,Q)$ can be taken to be arbitrarily large, without loss of generality. Our proofs differ from the ones given in \cite{karagila2020have}. We prove that, given a model $M\vDash\mathrm{ZFC}$ and sets $I,J\in M$, if there is a bijection $b:I\to J$ in an outer model $V$ in which $M$ is a class, then $b$ can be used to produce an explicit isomorphism in $V$ between the relativized forcing relations $\Vdash^M_{N(I,Q)}$ and $\Vdash^M_{N(J,Q)}$. We establish the relationship between such an isomorphism in $V$ and the existence of an elementary embedding in $M$ from $\Vdash_{N(I,Q)}$ to $\Vdash_{N(J,Q)}$. These maps, in $M$ and in $V$, can then be used to answer the relevant case of question 5.5 from \cite{karagila2020have}, and to prove a sharp version of the full theorem that $N(I,Q)\vDash\mathrm{BPI}$. We conclude section 5 by remarking on the fact that a case of our main results answers the $8^\mathrm{th}$ question on Karagila's website, and thus fills a gap to complete an alternate new proof of $\mathrm{BPI}$ in $N(I,Q)$, based on theorems from \cite{blass1979injectivity}.

In section 6, we formalize the notion that our adaptation of Harrington's proof gives a dynamical account of $\mathrm{BPI}$ in $N(I,Q)$. In the case of permutation models, we apply a version of our proof to show that the filter extension property (and hence $\mathrm{BPI}$) follows from the Ramsey property. This provides a new direct proof of one direction of Blass' theorem from \cite{blass1986prime}: that the Ramsey property implies $\mathrm{BPI}$ in permutation models. We then introduce the virtual Ramsey property for symmetric systems -- a dynamical condition that both generalizes the Ramsey property and abstracts core features of our adaptation of Harrington's proof. We prove that the virtual Ramsey property implies the filter extension property for symmetric systems, and hence is sufficient for $\mathrm{BPI}$ to hold in a symmetric extension. In this sense, we show that $\mathrm{BPI}$ can be proven in (all) permutation models and in (a certain class of) symmetric extensions by the same dynamical argument: the (virtual) Ramsey property implies the filter extension property, which in turn implies $\mathrm{BPI}$. We conclude the section by using this more general framework to revisit our proofs from sections 4 and 5. When $I$ is large, we confirm that the symmetric system for $N(I,Q)$ has the virtual Ramsey property. We then show that the virtual Ramsey property fails for the standard symmetric system that presents the Cohen model $N(\omega,2^{<\omega})$, despite the fact that $\mathrm{BPI}$ still holds in $N(\omega,2^{<\omega})$. From this, we infer that a symmetric system can have the filter extension property but not the virtual Ramsey property. We conclude by posing a question about the relationship between the virtual Ramsey property and $\mathrm{BPI}$ in symmetric extensions.

In section 7, we sketch two future directions for this work that are in preparation. For the reader's convenience, we also recall three questions that are posed throughout the paper.

If the reader wishes to quickly read the parts of this paper that deal with permutation models, note that subsections 2.1, 3.1, and 6.1 isolate our treatment of permutation models. The proofs in these subsections can be read independently of the rest of the paper. These subsections can also be read first as a gentle introduction to our treatment of symmetric extensions.

I would like to thank my advisors Isaac Goldbring and Martin Zeman for their time, support, and guidance while this work was being carried out.

\section{Preliminaries}
\subsection{Permutation Models}
Permutation models of $\mathrm{ZFA}$ are the precursors to symmetric extensions. They were introduced to give methods by which models without choice could be produced before the introduction of forcing. A comprehensive treatment of this topic can be found in \cite{jech2008axiom}.

The language of $\mathrm{ZF(C)A}$ includes the binary predicates $=,\in$, the constant symbol $0$ used to denote the empty set, and the constant symbol $A$ used to denote the set of atoms. There are several standard modifications that $\mathrm{ZFCA}$ makes to $\mathrm{ZFC}$. The primary adjustment is to include the axiom $\mathrm{A}$, which states that the ``atoms" in $A$ have no elements, but are distinct from one another and from the empty set. The rest of the $\mathrm{ZFC}$ axioms are adjusted accordingly to form the theory $\mathrm{ZFCA}$; $\mathrm{ZFA}$ is this theory without choice.

\begin{definition}
    Given a model $M\vDash\mathrm{ZFCA}$, with $A$ its set of atoms and $\mathrm{ON}$ the class of ordinals, for $\alpha\in \mathrm{ON}$, one can recursively define
    \begin{enumerate}[label=(\roman*)]
        \item $V_0(A)=A$,
        \item $V_{\alpha+1}(A)=2^{V_\alpha(A)}$,\footnote{Throughout the paper, we use the convention that $2^x$ denotes the powerset of a set $x$.}
        \item $V_\alpha(A)=\bigcup_{\beta<\alpha}V_\beta(A)$ if $\alpha$ is a limit ordinal.
    \end{enumerate}
\end{definition}

\begin{theorem}
    Given a model $M\vDash\mathrm{ZFCA}$, with $A$ its set of atoms,
    $$
    M=\bigcup_{\alpha\in\mathrm{ON}}V_\alpha(A).
    $$
\end{theorem}

\begin{definition}
    For $M\vDash\mathrm{ZFCA}$ and $x\in M$, define $\mathrm{rnk}(x)=\mathrm{inf}\{\alpha\ |\ x\in V_\alpha\}$.
\end{definition}

Throughout the paper, we will write $H\preceq K$ to denote that $H$ is a subgroup of $K$.

\begin{definition}
    For a model $M\vDash\mathrm{ZFCA}$, let $\mathrm{Aut}(A)$ be the group of permutations of the set $A$ of atoms in $M$. Given $\mathcal{G}\preceq\mathrm{Aut}(A)$, recursively define the action of $\pi\in\mathcal{G}$ on $M$ by 
    $$
    \pi x=\{\pi y\ |\ y\in x\}.
    $$
\end{definition}

\begin{lemma}\label{automorphisms of permutation model lemma}
    Every $\pi\in\mathrm{Aut}(A)$, acting on $M$ as above, is an $\in$-respecting automorphism of $M$.
\end{lemma}

\begin{definition}
    Let $M\vDash\mathrm{ZFCA}$, $A\in M$ be the set of atoms, and $\mathcal{G}\preceq\mathrm{Aut}(A)$. $\mathcal{F}$ is called a \textit{normal filter of subgroups} of $\mathcal{G}$ if
    \begin{enumerate}[label=(\roman*)]
        \item Every $K\in \mathcal{F}$ is a subgroup of $\mathcal{G}$,
        \item $\mathcal{G}\in\mathcal{F}$,
        \item if $K\in \mathcal{F}$ and $K\preceq H\preceq\mathcal{G}$, then $H\in\mathcal{F}$,
        \item if $K,H\in\mathcal{F}$ then $K\cap H\in\mathcal{F}$,
        \item if $\pi\in\mathcal{G}$ and $K\in\mathcal{F}$ then $\pi K\pi^{-1}\in\mathcal{F}$ (this is called the ``normality" condition).
    \end{enumerate}
\end{definition}

\begin{definition}
    For $x\in M\vDash\mathrm{ZFCA}$, let $\mathrm{sym}_\mathcal{G}(x)=\{\pi\in\mathcal{G}\ |\ \pi x=x\}$.
\end{definition}

\begin{definition}
    For any $y\in M\vDash\mathrm{ZFCA}$ and $K\preceq\mathrm{Aut}(A)$, let
    $$
    \mathrm{Orb}(y,K)=\{\pi y\ |\ \pi\in K\}.
    $$
\end{definition}

\begin{definition}
    Given a model $M\vDash\mathrm{ZFCA},$ a group $\mathcal{G}\preceq \mathrm{Aut}(A)$, and $\mathcal{F}$ a normal filter of subgroups of $\mathcal{G}$, $x\in M$ is called \textit{$\mathcal{F}$-symmetric} if $\mathrm{sym}_\mathcal{G}(x)\in\mathcal{F}$. Recursively, $x\in M$ is called and \textit{hereditarily $\mathcal{F}$-symmetric} if $x$ is $\mathcal{F}$-symmetric and every $y\in x$ is hereditarily $\mathcal{F}$-symmetric.
\end{definition}

Intuitively, a set $x\in M$ is $\mathcal{F}$-symmetric if it is fixed by an $\mathcal{F}$-large subgroup of $\mathcal{G}$. The normality condition ensures that if $x$ is $\mathcal{F}$-symmetric, then so is $\pi x$ for every $\pi\in \mathcal{G}$. The intent is that $\mathcal{G}$ and $\mathcal{F}$ can be selected so that hereditarily $\mathcal{F}$-symmetric sets cannot witness particular consequences of $\mathrm{AC}$, because they cannot sufficiently tell atoms apart.

\begin{definition}
    Let $M\vDash\mathrm{ZFCA}$, $A\in M$ be the set of atoms, $\mathcal{G}$ be a group of permutations of $A$, and $\mathcal{F}$ be a normal filter on $\mathcal{G}$. Define $M(A,\mathcal{G},\mathcal{F})=\{x\in M\ |\ x\mathrm{\ is\ hereditarily\ }\mathcal{F}-\mathrm{symmetric}\}$.
\end{definition}

\begin{theorem}
$M(A,\mathcal{G},\mathcal{F})\vDash\mathrm{ZFA}$.
\end{theorem}

\subsection{Symmetric Extensions}
A symmetric extension is an inner model of a $\mathbb{P}$-extension $M[G]$, defined analogously to a permutation model by using a group $\mathcal{G}\preceq\mathrm{Aut}(\mathbb{P})$ and a normal filter $\mathcal{F}$ of subgroups of $\mathcal{G}$. By considering automorphisms of the underlying forcing $\mathbb{P}$, symmetric extensions avoid the need to use atoms. This results in a model of $\mathrm{ZF}$, rather than $\mathrm{ZFA}$. The reader can again refer to \cite{jech2008axiom} for a complete exposition of this topic.

\begin{definition}\label{symmetric system definition}
    A \textit{symmetric system} is a tuple $\langle\mathbb{P},\mathcal{G},\mathcal{F}\rangle$, where $\mathbb{P}$ is a poset, $\mathcal{G}\preceq\mathrm{Aut}(\mathbb{P})$ is a subgroup of automorphisms of $\mathbb{P}$, and $\mathcal{F}$ is a normal filter on $\mathcal{G}$, as in definition 2.6.
\end{definition}

\begin{definition}
    For a poset $\mathbb{P}$ in a model $M\vDash\mathrm{ZFC}$, let $M^\mathbb{P}$ denote the class of $\mathbb{P}-$names in $M$.
\end{definition}

\begin{definition}
    Given $\pi\in\mathcal{G}\preceq\mathrm{Aut}(\mathbb{P})$, recursively define the action of $\pi$ on $M^\mathbb{P}$ by
    $$
    \pi\dot{x}=\{\langle\pi p,\pi\dot{z}\rangle\ |\ \langle p,\dot{z}\rangle\in\dot{x}\}.
    $$
\end{definition}

\begin{lemma}[Symmetry Lemma]\label{symmetry lemma}
    For $\pi\in\mathcal{G}\preceq\mathrm{Aut}(\mathbb{P})$,
    $$
    p\Vdash\phi(\dot{x}_0,...,\dot{x}_n)\iff\pi p\Vdash\phi(\pi\dot{x}_0,...,\pi\dot{x}_n).
    $$
\end{lemma}

\begin{definition}
    For a $\mathbb{P}$-name $\dot{x}$ and a group $\mathcal{G}\preceq\mathrm{Aut}(\mathbb{P})$, let $\mathrm{sym}_\mathcal{G}(\dot{x})=\{\pi\in\mathcal{G}\ |\ \pi\dot{x}=\dot{x}\}$.
\end{definition}

\begin{definition}
    For any $\mathbb{P}$-name $\dot{y}$, any condition $p\in\mathbb{P}$, and any $K\preceq\mathrm{Aut}(\mathbb{P})$, let
    $$
    \mathrm{Orb}(\langle p,\dot{y}\rangle,K)=\{\langle \pi p,\pi \dot{y}\rangle\ |\ \pi\in K\}.
    $$
\end{definition}

\begin{definition}
    For a symmetric system $\langle\mathbb{P},\mathcal{G},\mathcal{F}\rangle$, a $\mathbb{P}$-name
    $\dot{x}$ is \textit{$\mathcal{F}$-symmetric} if $\mathrm{sym}_\mathcal{G}(\dot{x})\in\mathcal{F}$. Recursively, we say that $\dot{x}$ is \textit{hereditarily $\mathcal{F}$-symmetric} if $\mathrm{sym}_\mathcal{G}(\dot{x})\in\mathcal{F}$ and every $\dot{y}\in\mathrm{rng}(\dot{x})$ is hereditarily $\mathcal{F}$-symmetric.
\end{definition}

As in the case of permutation models, the intent is that a symmetric system can be defined so that hereditarily $\mathcal{F}$-symmetric names cannot be forced to witness particular consequences of $\mathrm{AC}$.

\begin{definition}
     Given a symmetric system $\langle\mathbb{P},\mathcal{G},\mathcal{F}\rangle$, let $\mathrm{HS}_\mathcal{F}$ be the class of hereditarily $\mathcal{F}$-symmetric $\mathbb{P}$-names (with $\mathbb{P}$, $\mathcal{G}$, and the ground model $M$ understood from context).
\end{definition}

\begin{definition}
    Given a ground model $M\vDash\mathrm{ZFC}$, a symmetric system $\langle\mathbb{P},\mathcal{G},\mathcal{F}\rangle$, and an $(M,\mathbb{P})$-generic filter $G$, let $\mathrm{HS}_\mathcal{F}^G$ denote the class of interpretations of the names in $\mathrm{HS}_\mathcal{F}$ by $G$. We call $\mathrm{HS}_\mathcal{F}^G$ a \textit{symmetric extension} of $M$, presented by $\langle\mathbb{P},\mathcal{G},\mathcal{F}\rangle$.
    \end{definition}

\begin{theorem}
    If $M\vDash\mathrm{ZF(C)}$, $\langle\mathbb{P},\mathcal{G},\mathcal{F}\rangle$ is a symmetric system, and $G$ is a $\mathbb{P}$-generic filter, then $M\subseteq\mathrm{HS}_\mathcal{F}^G$ and
    $$
    \mathrm{HS}_\mathcal{F}^G\vDash\mathrm{ZF}.
    $$
\end{theorem}

\begin{remark}
    For a symmetric extension $N$ based on the symmetric system $\langle\mathbb{P},\mathcal{G},\mathcal{F}\rangle$, one can define the forcing relation $\Vdash_{\mathrm{HS}_\mathcal{F}}$, or $\Vdash_N$, by relativizing the forcing relation for $\mathbb{P}$ over $M$ to the class $\mathrm{HS}_\mathcal{F}$. In this paper, we will often use $\Vdash$ to mean either, and in any such context the correct formalization will be clear. 
\end{remark}

\subsection{The Finite Intersection Property}
Throughout the paper, we refer to the finite intersection property ($\mathrm{FIP}$): $F$ has the finite intersection property if the intersection of finitely many elements of $F$ is infinite. It is well-known that ``is infinite" is not a first-order definable concept, but in a model $M\vDash \mathrm{ZFC}$ there are many known and equivalent definitions of infinite that reference the natural numbers $\omega\in M$.

In a model $N\vDash\mathrm{ZF}$, many of these definitions of infinite are known to no longer be equivalent. Because the models we study in this paper do not have choice, we need to specify what is meant by the $\mathrm{FIP}$. In both permutation models and in symmetric extensions, note that $\omega$ is still definable in the model. We can give the following definition for ``infinite."

\begin{definition}\label{Infinite definition}
    In a model $M\vDash\mathrm{ZF(A)}$, we say that a set $x$ is \textit{infinite} if there are no bijections between $x$ and any $n\in\omega$.
\end{definition}

Note that this definition of ``infinite" is absolute between a permutation model or symmetric extension and their relevant outer model of $\mathrm{ZFC(A)}$.

\subsection{The Generalized Cohen Model}
This presentation follows \cite{stefanovic2023alternatives}, where this generalization of the Cohen model was introduced.

\begin{definition}
    A set $x\in M\vDash\mathrm{ZF(A)}$ is called \textit{Dedekind-finite} if there are no injections from $\omega$ into $x$ in $M$. In this paper, we will further require that Dedekind-finite sets are infinite, in the sense of definition \ref{Infinite definition}.
\end{definition}

The Cohen model primarily introduces a Dedekind-finite set of mutually generic Cohen reals; the generalized Cohen model, denoted $N(I,Q)$ for some index set $I$ and poset $Q$, primarily introduces a Dedekind-finite set of mutually $Q$-generic filters, which is defined with respect to an index set $I$ in the ground model. We will first define a poset that adds $I$-many mutually $Q$-generic filters.

\begin{definition}\label{P(I,Q)}
    For a poset $Q$ and an index set $I$, let $\mathbb{P}(I,Q)$ denote the partial order of finite functions from $I$ to $Q$ ordered by
    $$
    (\forall p,q\in\mathbb{P}(I,Q))(p\leq q)\iff(\mathrm{dom}(q)\subseteq\mathrm{dom}(p)\wedge(\forall i\in\mathrm{dom}(q))p(i)\leq q(i))).
    $$
\end{definition}

\begin{definition}
    Given a set $X$, let $[X]^\kappa=\{A\subseteq X\ |\ |A|=\kappa\}$ and $[X]^{<\kappa}=\bigcup_{\alpha<\kappa}[X]^\alpha.$
\end{definition}

\begin{definition}
    For a group $\mathcal{G}\preceq\mathrm{Aut}(I)$ and $E\in[I]^{<\omega}$, let $\mathrm{fix}_\mathcal{G}(E)=\{\pi\in\mathcal{G}\ |\ \pi\mathrm{\ fixes\ }E\mathrm{\ pointwise}\}.$
\end{definition}

\begin{definition}\label{Generalized Cohen System}
    Let $M\vDash\mathrm{ZFC}$, and define the symmetric system $\langle\mathbb{P}(I,Q),\mathcal{G},\mathcal{F}\rangle$ by setting
    \begin{enumerate}[label=(\roman*)]
        \item $\mathbb{P}(I,Q)$ as in definition \ref{P(I,Q)}, with $I$ an arbitrary index and $Q$ an arbitrary poset,
        \item $\mathcal{G}$ the group of finite permutations of $I$, acting on $p\in\mathbb{P}(I,Q)$ by $\mathrm{dom}(\pi p)=\pi\mathrm{dom}(p)$ and
    $$
    (\forall i\in\mathrm{dom}(p))\pi p(\pi i)=p(i),\footnote{By this group action, $\mathcal{G}$ can be formally identified with a subgroup of $\mathrm{Aut}(\mathbb{P}(I,Q))$, as is required in definition \ref{symmetric system definition}.}
    $$
        \item the normal filter $\mathcal{F}$ on $\mathcal{G}$ generated by subgroups of the form $\mathrm{fix}_\mathcal{G}(E)$, for $E\in[I]^{<\omega}$.\footnote{For every $K\in\mathcal{F}$, this guarantees that we can find some $E\in[I]^{<\omega}$ for which $\mathrm{fix}_\mathcal{G}(E)\preceq K\preceq\mathcal{G}$. In general, however, $K\in\mathcal{F}$ may not equal $\mathrm{fix}_\mathcal{G}(E)$ for any $E\in[I]^{<\omega}$.}
    \end{enumerate}
    If we would like to emphasize the index $I$, we write $\langle\mathbb{P}(I,Q),\mathcal{G}_I,\mathcal{F}_I\rangle$. We call the corresponding symmetric extension of $M$ a generalized Cohen model, and denote it by $N(I,Q)$.
\end{definition}

Throughout the paper, whenever we refer to $\langle\mathbb{P}(I,Q),\mathcal{G},\mathcal{F}\rangle$, we will mean the symmetric system from definition \ref{Generalized Cohen System}. We typically omit the $\mathbb{P}(I,Q)$-generic filter $G$ used to define $N(I,Q)$, but if we choose to emphasize this filter, we will write $N^G(I,Q)$.

\begin{definition}\label{Dedekind Finite A}
    For $n\in I$, let $\dot{a}_n=\{\langle p,\check{q}\rangle\ |\ p(n)=q\}$, and let  $\dot{A}=\{\langle1,\dot{a}_n\rangle\ |\ n\in I\}$.
\end{definition}

\begin{theorem}
    For a generalized Cohen model $N^G(I,Q)$ over a ground model $M\vDash\mathrm{ZFC}$, $\dot{a}_n^G\in N^G(I,Q)$ is $Q$-generic over $M$ and  $\dot{A}^G\in N^G(I,Q)$ and is Dedekind-finite in $N^G(I,Q)$.
\end{theorem}

\begin{definition}\label{pi_vec{alpha}}
    For $D\in[I]^{<\omega}$, disjoint sets $H_i\subseteq I$ (for $i\in D$), $\vec{\alpha}=\{\langle i,\alpha_i\rangle\ |\ i\in D\}\in\prod_{i\in D}H_i$, and $\mathcal{G}$ as in definition \ref{Generalized Cohen System}, define the permutation $\pi_{\vec{\alpha}}\in\mathcal{G}$ by
    $$
    \pi_{\vec{\alpha}}=\prod_{i\in D}(i\ \alpha_i).
    $$
\end{definition}

\begin{lemma}\label{Minimal Supports Lemma}
    For the symmetric system $\langle\mathbb{P}(I,Q),\mathcal{G},\mathcal{F}\rangle$, if $\dot{x}\in\mathrm{HS}_\mathcal{F}$, there is a set $X\in[I]^{<\omega}$ so that $\mathrm{fix}_\mathcal{G}(X)\preceq\mathrm{sym}_\mathcal{G}(\dot{x})$ and $X\subseteq X'$ for every other $X'\in[I]^{<\omega}$ so that $\mathrm{fix}_\mathcal{G}(X')\preceq\mathrm{sym}_\mathcal{G}(\dot{x})$.
\end{lemma}
\begin{proof}
    See \cite{jech2008axiom} for the proof of ``least supports" in ``support models."
\end{proof}

\begin{definition}\label{Support Definition}
    $X\in[I]^{<\omega}$ as in lemma 2.26 is called the support of $\dot{x}$; we write $\mathrm{supp}(\dot{x})=X$.
\end{definition}

\begin{definition}
    For $p\in \mathbb{P}(I,Q)$, let $\mathrm{supp}(p)=\mathrm{dom}(p)$.
\end{definition}

\begin{definition}
    Given $p\in\mathbb{P}(I,Q)$ and $X\subseteq I$, define the \textit{restriction of $p$ to $X$} to be the (possibly empty) function $p\upharpoonright_X=\{\langle i,p(i)\rangle\ |\ i\in X\}$.
\end{definition}

The following lemma states an essential property that is used throughout the paper.

\begin{lemma}\label{Support Restriction Lemma}
    For the symmetric system $\langle\mathbb{P}(I,Q),\mathcal{G},\mathcal{F}\rangle$, if $\dot{x}_0,...,\dot{x}_n\in HS_\mathcal{F}$ have supports $X_0,...,X_n$, and $p\Vdash\phi(\dot{x}_0,...,\dot{x}_n),$ then for $X\supseteq\bigcup_{i\leq n}X_i$,
    $$
    p\upharpoonright_X\Vdash\phi(\dot{x}_0,...,\dot{x}_n).
    $$
\end{lemma}
\begin{proof}
    By the symmetry lemma \ref{symmetry lemma}, if $\pi \in\mathrm{fix}_\mathcal{G}(X)$, then $\pi p\Vdash\phi(\dot{x}_0,...,\dot{x}_n).$ Because the set of conditions $D=\{\pi p\ |\ \pi\in\mathrm{fix}_\mathcal{G}(X)\}$ is dense below $p\upharpoonright_X$ in $\mathbb{P}(I,Q)$, we get that $p\upharpoonright_X\Vdash\phi(\dot{x}_0,...,\dot{x}_n).$
\end{proof}

\begin{lemma}\label{image support lemma}
    Let $\langle\mathbb{P}(I,Q),\mathcal{G},\mathcal{F}\rangle$ be a symmetric system for a generalized Cohen model, $p\in\mathbb{P}(I,Q)$, and $\dot{x}\in\mathrm{HS}_\mathcal{F}$. If $\mathrm{supp}(p),\mathrm{supp}(\dot{x})\subseteq D\in[I]^{<\omega}$, then for every $\pi,\tau\in\mathcal{G}$ with $\pi\upharpoonright_D=\tau\upharpoonright_D$,
    $$
    \pi p=\tau p\wedge\pi\dot{x}=\tau\dot{x}.
    $$
\end{lemma}
\begin{proof}
    Because $\pi\upharpoonright_D=\tau\upharpoonright_D$, we have that $\pi^{-1}\tau\in\mathrm{fix}_\mathcal{G}(D)$. Then,
    $$
    \pi p=\pi(\pi^{-1}\tau p)=\tau p\wedge\pi \dot{x}=\pi(\pi^{-1}\tau \dot{x})=\tau \dot{x}.
    $$
\end{proof}

\begin{theorem}[Halpern, L\'evy,  \cite{halpern1971boolean}, Stefanovi\'c, \cite{stefanovic2023alternatives}\footnote{This is proven for the Cohen model in \cite{halpern1971boolean}; in \cite{stefanovic2023alternatives}, the proof from \cite{halpern1971boolean} is shown to generalize $N(I,Q)$.}]
    $N(I,Q)\vDash \mathrm{BPI}$.
\end{theorem}

\section{The Filter Extension Property}

In a model of $\mathrm{ZFC(A)}$, the standard proof of $\mathrm{BPI}$ involves showing that
\begin{enumerate}
    \item an arbitrary filter $F$, on an arbitrary set $x$, can be extended to include an arbitrary subset $y\subseteq x$ or its complement $x\setminus y$, and
    \item Zorn's lemma can be applied to produce a maximal filter $U\supseteq F$ on $x$; by $(1)$, the maximality of $U$ ensures that it ``completes" the extension of $F$ and must be an ultrafilter on $x$.
\end{enumerate}
While the full axiom of choice fails in permutation models and symmetric extensions, it can still be found to hold up to a coarse degree -- certain choice functions are also hereditarily symmetric. In this section, we will introduce a coarsening of $(1)$ called the filter extension property (for permutation models and for symmetric systems), which can be used to prove $\mathrm{BPI}$ by a coarse version of $(2)$.\footnote{The filter extension property for permutation models and for symmetric extensions are two distinct conditions. We will typically just write ``the filter extension property," and from context the version will be clear.} We will also show that the filter extension property, respectively in permutation models and symmetric extensions, is necessary for $\mathrm{BPI}$ to hold in these models. Subsection 3.1 contains these proofs for permutation models; subsection 3.2 contains their counterparts for symmetric extensions.

\subsection{The Filter Extension Property for Permutation Models}
Definition \ref{filter extension property for permutation models} can be compared to $(1)$ from the standard proof of $\mathrm{BPI}$ in a model of $\mathrm{ZFC(A)}$; theorem \ref{FEP proves BPI in permutation models} can be compared to $(2)$. We can motivate definition \ref{filter extension property for permutation models} as follows. Given $F,x\in M(A,\mathcal{G},\mathcal{F})$, with $F$ a filter on $x$, we would like to mimic the step (1) by considering how to extend $F$ by some minimal increment while retaining the $\mathrm{FIP}$. If this can be done in a manner that is stable under iteration, we can then apply a coarse version of $(2)$ to produce an ultrafilter extending $F$.

Note that if $w\in M(A,\mathcal{G},\mathcal{F})$, then
$$
w=\bigcup_{y\in w}\mathrm{Orb}(y,\mathrm{sym}_\mathcal{G}(w)).
$$
This is essentially a partition of $w$, up to the fact that some orbits may be repeated in the union above. In this sense, one can view these orbits as the building blocks of (hereditarily) symmetric sets. Moreover, this gives the sense in which choice holds coarsely in a permutation model: if $f\in M\vDash\mathrm{ZFCA}$ is a wellordering of the set $\{\mathrm{Orb}(y,\mathrm{sym}_\mathcal{G}(w))\ |\ y\in w\}$, then $\mathrm{sym}_\mathcal{G}(f)=\mathrm{sym}_\mathcal{G}(w)$, and so $f\in M(A,\mathcal{G},\mathcal{F})$. Thus, while $w$ itself may not be wellorderable in $M(A,\mathcal{G},\mathcal{F})$, it is still possible to wellorder the partition of $w$ given by the $\mathrm{sym}_\mathcal{G}(w)$-orbits.

To extend the filter $F$ by a minimal increment thus involves selecting a subset $y\in M(A,\mathcal{G},\mathcal{F})$ of $x$, a group $K\preceq\mathrm{sym}_\mathcal{G}(F)\cap\mathrm{sym}_\mathcal{G}(x)$, and extending to
$$
F\cup \mathrm{Orb}(y,K)
$$
or
$$
F\cup\mathrm{Orb}(x\setminus y,K),
$$
while ensuring that the finite intersection property $(\mathrm{FIP}$) holds in the chosen extension.\footnote{Here we mean ``minimal" more in a figurative than a literal sense, as it may be possible to choose a smaller group $K$ with this same property.} If this process can be iterated, using the same group $K\preceq\mathrm{sym}_\mathcal{G}(F)\cap\mathrm{sym}_\mathcal{G}(x)$ at every stage of the iteration, the coarse form of choice in $M(A,\mathcal{G},\mathcal{F})$ will guarantee that the process terminates in an ultrafilter $U\supseteq F$ on $x$ in $M(A,\mathcal{G},\mathcal{F})$, with $K\preceq\mathrm{sym}(U)$. Definition \ref{filter extension property for permutation models} expresses that the hypotheses in this sketch hold.

\begin{definition}\label{filter extension property for permutation models}
    A permutation model $M(A,\mathcal{G},\mathcal{F})$ has the \textit{filter extension property} (for permutation models) if $\mathcal{F}$ has a base of subgroups $K$, for which the following property holds: for every $F,x,y\in M(A,\mathcal{G},\mathcal{F})$ such that $K\preceq\mathrm{sym}_\mathcal{G}(F)\cap\mathrm{sym}_\mathcal{G}(x)$ and
    $$
    M(A,\mathcal{G},\mathcal{F})\vDash ``F\subseteq2^x\mathrm{\ has\ the\ FIP}"\wedge y\subseteq x,
    $$
    there is a $y'\in\{y,x\setminus y\}$ such that
    $$
    M(A,\mathcal{G},\mathcal{F})\vDash ``F\cup\mathrm{Orb}(y',K)\subseteq2^x\mathrm{\ has\ the\ FIP}."
    $$
\end{definition}

Theorem \ref{FEP proves BPI in permutation models} uses the fact that the filter extension property guarantees groups $K\in\mathcal{F}$ that are stable under iterations of filter extensions, in the sense described above.

\begin{theorem}\label{FEP proves BPI in permutation models}
    If $M\vDash\mathrm{ZFCA}$ and $M(A,\mathcal{G},\mathcal{F})$ has the filter extension property, then
    $$
    M(A,\mathcal{G},\mathcal{F})\vDash\mathrm{BPI}.
    $$
\end{theorem}
\begin{proof}
    Let $F,x\in M(A,\mathcal{G},\mathcal{F})$ so that
    $$
    M(A,\mathcal{G},\mathcal{F})\vDash ``F\subseteq2^x\mathrm{\ has\ the\ FIP},"
    $$
    and let $K\preceq\mathrm{sym}_\mathcal{G}(F)\cap\mathrm{sym}_\mathcal{G}(x)$ witness the filter extension property. At the cost of lengthening a standard application of Zorn's lemma, we would like our proof to emphasize the notion that $K$ can be used to extend $F$, by minimal increments, in a manner that is stable under iteration.
    
    In the outer model $M\vDash\mathrm{ZFCA}$, use choice to wellorder the subsets $y\subseteq x$ that appear in $M(A,\mathcal{G},\mathcal{F})$; denote this wellordering as $\{y_{\lambda}|\lambda<\kappa\}$. Using choice again in $M$, define the sequence $F_{\lambda}$, for $\lambda\leq\kappa$, in the following manner.
    
    Let $F_0$ be the filter $F$ on $x$. At successor stages $\lambda+1$, choose between
    $$
    F_{\lambda+1}=F_{\lambda}\cup\mathrm{Orb}(y_\lambda,K)
    $$
    and
    $$
    F_{\lambda+1}=F_\lambda\cup\mathrm{Orb}(x\setminus y_\lambda,K),
    $$
    so that $F_{\lambda+1}$ has the $\mathrm{FIP}$ and $K\preceq \mathrm{sym}_\mathcal{G}(F_{\lambda+1})\cap\mathrm{sym}_\mathcal{G}(x)$. If we suppose inductively that $F_{\lambda}$ meets these constraints, then the fact that $K$ witnesses the filter extension property ensures that $F_{\lambda+1}$ can be chosen to meet these constraints as well.

    At limit stages $\lambda$, let
    $$
    F_{\lambda}=\bigcup_{\xi<\lambda}F_{\xi}.
    $$
    If we again suppose inductively that for every $\xi<\xi'<\lambda$, $K\preceq \mathrm{sym}_\mathcal{G}(F_{\xi})\cap\mathrm{sym}_\mathcal{G}(x)$, $F_\xi\subseteq F_{\xi'}$, and $F_{\xi}$ has the $\mathrm{FIP}$, then $K\preceq\mathrm{sym}_\mathcal{G}(F_\lambda)\cap\mathrm{sym}_\mathcal{G}(x)$ and $F_\lambda$ has the $\mathrm{FIP}$ as well.

    Then, $U=F_\kappa\supseteq F$ is an ultrafilter on $x$ in $M(A,\mathcal{G},\mathcal{F})$, with $K\preceq\mathrm{sym}_\mathcal{G}(U)$.
\end{proof}

Note that theorem 3.2 is essentially the same as the step in Halpern's contradiction framework to consider a filter $U\in M(A,\mathcal{G},\mathcal{F})$ on $x$ that extends $F$ and is maximal with respect to having some specified $K\preceq\mathrm{sym}_\mathcal{G}(U)$. The next step in Halpern's contradiction framework, to produce a contradiction from the assumption that there is a subset $y\in M(A,\mathcal{G},\mathcal{F})$ of $x$ for which $y\not\in U$ and $x\setminus y\not\in U$, can be rephrased as a contradiction proof of the filter extension property.

We will now prove that the filter extension property is necessary for $\mathrm{BPI}$ to hold in a permutation model, whereby this relationship to Halpern's contradiction framework can be seen to be inevitable. Note that the filter extension property ensures that filters can be extended in a manner based \textit{uniformly} on witness groups $K\preceq\mathrm{sym}_\mathcal{G}(F)\cap\mathrm{sym}_\mathcal{G}(x)$, irrespective of the filter $F$ or the set $x$. This uniformity can be explained in terms of the uniform method given in \cite{blass1979injectivity} to define ultrafilters in models of $\mathrm{SVC}+\mathrm{BPI}$.

\begin{definition}[Blass, \cite{blass1979injectivity}]
    Let $\mathrm{SVC}$ (\textit{small violations of choice}) be the axiom that there is a set $S$ so that, for every set $a$, there is an ordinal $\alpha$ and a function from a subset of $S\times \alpha$ onto $a$.\footnote{In \cite{blass1979injectivity}, the original statement is that there is a function from $S\times \alpha$ onto $a$. Blass notes that this is equivalent to the form stated here, since all points outside a subset of $S\times \alpha$ can be mapped to a single element of $a$.\label{Blass footnote}}
\end{definition}

\begin{theorem}[Blass, \cite{blass1979injectivity}]\label{SVC iff AC in extension}
    A model of $\mathrm{ZF}$ satisfies $\mathrm{SVC}$ if and only if some generic extension is a model of $\mathrm{ZFC}$.
\end{theorem}

We will sketch Blass' proof of theorem \ref{Blass SVC BPI Theorem}, and emphasize definability considerations that we will reference in our proofs of \ref{fep equiv to bpi in permutation models}$(i)\Rightarrow(ii)$ and \ref{BPI equivalent to FEP for sym ext}$(i)\Rightarrow(ii)$, that $\mathrm{BPI}$ implies the respective versions of the filter extension properties.

\begin{theorem}[Blass, \cite{blass1979injectivity}]\label{Blass SVC BPI Theorem}
    If $S$ witnesses $\mathrm{SVC}$, then $\mathrm{BPI}$ holds if and only if there is an ultrafilter $U'$ on the set $S^{<\omega}$ so that $\{p\in S^{<\omega}|s\in\mathrm{rng}(p)\}\in U'$ for every $s\in S$.
\end{theorem}
\begin{proof}[Proof sketch]
    The ``only if" direction is clear, so we will prove the ``if" direction. Let $B$ be a boolean algebra, $S$ witness $\mathrm{SVC}$, $f:S\times\alpha\to B$ be a surjection (as is guaranteed in footnote \ref{Blass footnote}), and $U'$ be an ultrafilter as in the theorem statement. We claim that there is a prime ideal $I$ of $B$ that is definable with respect to the parameters $f,S$ and $U'$. After proving this claim, we will show how to modify the argument to get the same result if the domain of $f$ is an arbitrary subset of $S\times\alpha$.
    
    Every $p\in S^{<\omega}$ is itself a wellordering of $\mathrm{rng}(p)$. Based on this, wellorder $\mathrm{rng}(p)\times\alpha$ by the lexicographical order; this is definable with respect to the parameter $p$. Then, wellorder $f(\mathrm{rng}(p)\times\alpha)$ by sending points to the least element of their preimage in $\mathrm{rng}(p)\times\alpha$; this is definable with respect to the parameters $p$ and $f$. Let $B_p$ be the subalgebra of $B$ generated by $f(\mathrm{rng}(p)\times\alpha)$. Define a prime ideal $I_p$ of $B_p$ based on the wellorder of $f(\mathrm{rng}(p)\times\alpha)$, with respect to the parameters $p$ and $f$. One can then define, with respect to the parameters $f,S$, and $U'$, the prime ideal of $B$:
    $$
    I=\{b\in B\ |\ \{p\in S^{<\omega}\ |\ b\in I_p\}\in U'\}.
    $$
    
    If we had taken $\mathrm{dom}(f)\subsetneq S\times\alpha$, the same proof works to define a prime ideal $I$ of $B$, if one takes the convention that $f(\mathrm{rng}(p)\times\alpha)=\{f(y)\ |\ y\in(\mathrm{rng}(p)\times\alpha)\cap\mathrm{dom}(f)\}$.
\end{proof}

For the reader's convenience, we will also sketch Blass' proof of theorem \ref{SVC in permutation models}.

\begin{theorem}[Blass, \cite{blass1979injectivity}]\label{SVC in permutation models}
    Permutation models satisfy $\mathrm{SVC}$.
\end{theorem}

\begin{proof}[Proof sketch]
    Consider the permutation model $M(A,\mathcal{G},\mathcal{F})$. For $H,K\in \mathcal{F}$, with $H\preceq K$, let $y_{H,K}$ be a set in $M(A,\mathcal{G},\mathcal{F})$ for which $\mathrm{sym}_\mathcal{G}(y_{H,K})=H$. Let 
    $$
    s=\bigcup_{H,K\in\mathcal{F},H\preceq K}\mathrm{Orb}(y_{H,K},K),\mathrm{\ and\ }S=\bigcup_{y\in s}\mathrm{Orb}(y,\mathcal{G}),
    $$
    so that $\mathrm{sym}_\mathcal{G}(S)=\mathcal{G}$, and thus $S\in M(A,\mathcal{G},\mathcal{F})$. For every $x\in M(A,\mathcal{G},\mathcal{F})$, we have that
    $$
    x=\bigcup_{y\in x}\mathrm{Orb}(y,\mathrm{sym}_\mathcal{G}(x)),
    $$
    which is again a partition of $x$ up to the repetition of orbits. Let $\alpha=|\{\mathrm{Orb}(y,\mathrm{sym}_\mathcal{G}(x))\ |\ y\in x\}|$. For every $y\in x$, there is a natural bijection between $\mathrm{Orb}(y,\mathrm{sym}_\mathcal{G}(x))$ and $\mathrm{Orb}(y_{H,K},K)\subseteq S$, with $H=\mathrm{sym}_\mathcal{G}(y)$ and $K=\mathrm{sym}_\mathcal{G}(x)$. Accordingly, there is a natural bijection $f\in M(A,\mathcal{G},\mathcal{F})$ between a subset of $S\times\alpha$ and $x$, so that $\mathrm{sym}_\mathcal{G}(f)=\mathrm{sym}_\mathcal{G}(x)$. Thus, $S$ witnesses $\mathrm{SVC}$ for $M(A,\mathcal{G},\mathcal{F})$.
\end{proof}

We can now extend theorem \ref{FEP proves BPI in permutation models} to an equivalence.

\begin{theorem}\label{fep equiv to bpi in permutation models}
    If $M\vDash\mathrm{ZFCA}$, the following are equivalent:
    \begin{enumerate}[label=(\roman*)]
        \item $M(A,\mathcal{G},\mathcal{F})\vDash\mathrm{BPI}$,
        \item $M(A,\mathcal{G},\mathcal{F})$ has the filter extension property.
    \end{enumerate}
\end{theorem}
\begin{proof}
The direction $(ii)\Rightarrow(i)$ is given in theorem \ref{FEP proves BPI in permutation models}, so we will prove the remaining direction $(i)\Rightarrow(ii)$. Fix a permutation model $M(A,\mathcal{G},\mathcal{F})\vDash\mathrm{BPI}$. By theorem \ref{Blass SVC BPI Theorem}, there is some $S\in M(A,\mathcal{G},\mathcal{F})$ that witnesses $\mathrm{SVC}$. In $M(A,\mathcal{G},\mathcal{F})$, let $B$ be a Boolean algebra, let $S$ witness $\mathrm{SVC}$, let $f$ be a surjection from a subset of $S\times\alpha$ to $B$, and let $U'$ be the ultrafilter on $S^{<\omega}$ from the statement of theorem \ref{Blass SVC BPI Theorem}. By the proof sketch of theorem \ref{SVC in permutation models}, we can assume that $\mathrm{sym}_\mathcal{G}(S)=\mathcal{G}$ and $\mathrm{sym}_\mathcal{G}(f)=\mathrm{sym}_\mathcal{G}(B)$. By the proof sketch of theorem \ref{Blass SVC BPI Theorem}, we have a formula $\phi(w_0,w_1,w_2,w_3)$ so that the following is a prime ideal in $B$:
    $$
    I=\{b\ |\ \phi(b,f,S,U')\}.
    $$
    
Given a filter $F$ on a set $x$ in $M(A,\mathcal{G},\mathcal{F})$, one can let $B=2^x/F$ to get a first-order formula $\varphi(w_0,w_1,w_2,w_3)$ for which
    \begin{equation}\label{Ultrafilter definition}
    U=\{y\ |\ \varphi(y,f,S,U')\}
    \end{equation}
is an ultrafilter on $x$ extending $F$. By this definition and by lemma \ref{automorphisms of permutation model lemma}, we have that

$$
\mathrm{sym}_\mathcal{G}(f)\cap\mathrm{sym}_\mathcal{G}(S)\cap\mathrm{sym}_\mathcal{G}(U')\preceq \mathrm{sym}_\mathcal{G}(U).
$$
Because $\mathrm{sym}_\mathcal{G}(S)=\mathcal{G}$ and $\mathrm{sym}_\mathcal{G}(2^x/F)=\mathrm{sym}_\mathcal{G}(f)$, we then have
\begin{equation}\label{sym(U) bound}
    \mathrm{sym}_\mathcal{G}(x)\cap\mathrm{sym}_\mathcal{G}(F)\cap\mathrm{sym}_\mathcal{G}(U')\preceq \mathrm{sym}_\mathcal{G}(U).
\end{equation}

Now let $M(A,\mathcal{G},\mathcal{F})\vDash\mathrm{BPI}$, $S\in M(A,\mathcal{G},\mathcal{F})$ witness $\mathrm{SVC}$, $U'$ be the ultrafilter on $S^{<\omega}$ witnessing theorem \ref{Blass SVC BPI Theorem}, and $K=\mathrm{sym}_\mathcal{G}(U')$. The normal filter $\mathcal{F}$ has a base of subgroups $\{H\cap K|H\in\mathcal{F}\}$; we will show that this base witnesses the filter extension property. Let $H\in\mathcal{F}$, and let $F\in M(A,\mathcal{G},\mathcal{F})$ be a filter on a set $x\in M(A,\mathcal{G},\mathcal{F})$, so that $H\cap K\preceq\mathrm{sym}_\mathcal{G}(F)\cap\mathrm{sym}_\mathcal{G}(x)$. Let $U\supseteq F$ be the ultrafilter on $x$ in $M(A,\mathcal{G},\mathcal{F})$, defined by the formula $\varphi$ as in line (\ref{Ultrafilter definition}). For every $y\subseteq x$ in $M(A,\mathcal{G},\mathcal{F})$, there is a $y'\in\{y,x\setminus y\}$ so that
    $$
    F\cup\mathrm{Orb}(y',\mathrm{sym}_\mathcal{G}(U'))\subseteq U.
    $$
From line (\ref{sym(U) bound}), we can conclude that
    $$
    H\cap K\preceq \mathrm{sym}_\mathcal{G}(F)\cap\mathrm{sym}_\mathcal{G}(x)\cap\mathrm{sym}_\mathcal{G}(U')\preceq\mathrm{sym}_\mathcal{G}(U),
    $$
and so, in particular,
    $$
    M(A,\mathcal{G},\mathcal{F})\vDash``F\cup\mathrm{Orb}(y',H\cap K)\subseteq 2^x\mathrm{\ has\ the\ FIP}."
    $$
    Thus, $H\cap K$ is a witness of the filter extension property. Because $H\cap K$ was an arbitrary group in the specified base for $\mathcal{F}$, $M(A,\mathcal{G},\mathcal{F})$ has the filter extension property.
\end{proof}

\subsection{The Filter Extension Property for Symmetric Systems}
We can now generalize these definitions and proofs to the context of symmetric extensions. Definition \ref{FEP for sym systems definition} can be compared to definition \ref{filter extension property for permutation models} and to $(1)$ from the standard proof of $\mathrm{BPI}$ in a model of $\mathrm{ZFC}$; theorem \ref{fep implies bpi in sym ext} can be compared to theorem \ref{FEP proves BPI in permutation models} and to $(2)$; theorem \ref{BPI equivalent to FEP for sym ext} can be compared to theorem \ref{fep equiv to bpi in permutation models}. We will misuse notation slightly by writing $2^{\dot{x}}$ to denote a name for the powerset of $\dot{x}$ and $\dot{x}\setminus\dot{y}$ to denote a name for the complement of $\dot{y}$ in $\dot{x}$.

To motivate definition \ref{FEP for sym systems definition}, let $\langle\mathbb{P},\mathcal{G},\mathcal{F}\rangle$ be a symmetric system, and let $\dot{F},\dot{x},\dot{y}\in\mathrm{HS}_\mathcal{F}$ be names for which
$$
1\Vdash``\dot{F}\subseteq2^{\dot{x}}\mathrm{\ has\ the\ FIP}"\wedge\dot{y}\subseteq\dot{x}.
$$
We would again like to consider how to extend $\dot{F}$ by some minimal increment to match $(1)$, in a manner that is stable under iteration to produce a name $\dot{U}$ that is forced to be an ultrafilter on $\dot{x}$.

Note that if $\dot{w}\in\mathrm{HS}_\mathcal{F}$, then
$$
\dot{w}=\bigcup_{\langle q,\dot{y}\rangle\in\dot{w}}\mathrm{Orb}(\langle q,\dot{y}\rangle,\mathrm{sym}_\mathcal{G}(\dot{w})),
$$
which can be expressed as a partition up to some repetition of orbits. We can draw the same conclusions as in the case for permutation models: these orbits are the building blocks for (hereditarily) symmetric names, and one can identify a name $\dot{f}$ for a wellordering of the $\mathrm{sym}_\mathcal{G}(\dot{w})$-orbits that partition $\dot{w}$, so that $\mathrm{sym}_\mathcal{G}(\dot{f})=\mathrm{sym}_\mathcal{G}(\dot{w})$.

In this case, extending $\dot{F}$ by some minimal increment involves selecting some condition $q$, some name $\dot{y}$ for a subset of $\dot{x}$, and some group $K\preceq\mathrm{sym}_\mathcal{G}(\dot{F})\cap\mathrm{sym}_\mathcal{G}(\dot{x})$, then extending to either
$$
\dot{F}\cup\mathrm{Orb}(\langle q,\dot{y}\rangle,K)
$$
or
$$
\dot{F}\cup\mathrm{Orb}(\langle q,\dot{x}\setminus\dot{y}\rangle,K),
$$
so that $1$ forces that the $\mathrm{FIP}$ is still met.\footnote{Again, we mean ``minimal" more figuratively than literally, as it may be possible to choose a smaller group $K$ with this same property.} We would again like to use the same group $K$ at every stage of an iteration that terminates in a name $\dot{U}$, with $K\preceq\mathrm{sym}_\mathcal{G}(\dot{U})$, and so that
$$
1\Vdash_{\mathrm{HS}_{\mathcal{F}}}``\dot{U}\supseteq\dot{F}\mathrm{\ is\ an\ ultrafilter\ on\ }\dot{x}."
$$
In this context, doing so will involve extending to include a dense set of conditions $q$ that force whether a given name $\dot{y}$ or its complement $\dot{x}\setminus\dot{y}$ appears in $\dot{U}$. Definition \ref{FEP for sym systems definition} expresses that the hypotheses in this sketch hold.

\begin{definition}\label{FEP for sym systems definition}
    A symmetric system $\langle\mathbb{P},\mathcal{G},\mathcal{F}\rangle$ has the \textit{filter extension property} (for symmetric systems, over a ground model $M$) if $\mathcal{F}$ has a base of subgroups $K$ for which the following holds: for every $\dot{F},\dot{x},\dot{y}\in \mathrm{HS}_\mathcal{F}$ such that $K\preceq \mathrm{sym}(\dot{F})\cap\mathrm{sym}(\dot{x})$ and
    $$
    1\Vdash``\dot{F}\subseteq2^{\dot{x}}\mathrm{\ has\ the\ FIP}"\wedge\dot{y}\subseteq\dot{x},
    $$
    and for every $p\in\mathbb{P}(I,Q)$, there is a $q\leq p$ and a $y'\in\{\dot{y},\dot{x}\setminus\dot{y}\}$ so that
    $$
    1\Vdash``\dot{F}\cup\mathrm{Orb}(\langle q,y'\rangle,K)\mathrm{\ has\ the\ FIP}."
    $$
\end{definition}

When proving that the filter extension property for a symmetric system implies $\mathrm{BPI}$ in the corresponding symmetric extension, we would like to emphasize the notion that definition \ref{FEP for sym systems definition} allows for the stable iteration of the process of extending hereditarily symmetric names for filters by minimal increments. This again comes at the cost of lengthening a standard application of Zorn's lemma. Lemma \ref{extend to include new subset} shows that the filter extension property implies the ability to extend a name for a filter $\dot{F}\in\mathrm{HS}_\mathcal{F}$ on a set $\dot{x}\in\mathrm{HS}_\mathcal{F}$ to choose between an arbitrary subset $\dot{y}\in\mathrm{HS}_\mathcal{F}$ and its complement in $\dot{x}$.

\begin{lemma}\label{extend to include new subset}
    If $M\vDash\mathrm{ZFC}$ and if $K\in\mathcal{F}$ is a witness of the filter extension property for $\langle\mathbb{P},\mathcal{G},\mathcal{F}\rangle$, then for every $\dot{F},\dot{x},\dot{y}\in\mathrm{HS}_\mathcal{F}$ such that $K\preceq\mathrm{sym}_\mathcal{G}(\dot{F})\cap\mathrm{sym}_\mathcal{G}(\dot{x})$ and
    $$
    1\Vdash``\dot{F}\subseteq2^{\dot{x}}\mathrm{\ has\ the\ FIP}"\wedge\dot{y}\subseteq\dot{x},
    $$
    there is a name $\dot{F}_{\dot{y}}\supseteq\dot{F}$, with $K\preceq\mathrm{sym}_\mathcal{G}(\dot{F}_{\dot{y}})\cap\mathrm{sym}_\mathcal{G}(\dot{x})$, so that
    $$
    1\Vdash``\dot{F}_{\dot{y}}\subseteq2^{\dot{x}}\mathrm{\ has\ the\ FIP}"\wedge(\dot{y}\in\dot{F}_{\dot{y}}\vee\dot{x}\setminus\dot{y}\in\dot{F}_{\dot{y}}).
    $$
\end{lemma}
\begin{proof}
    Fix $\dot{F},\dot{x},\dot{y}\in\mathrm{HS}_\mathcal{F}$ and $K\preceq\mathrm{sym}_\mathcal{G}(\dot{F})\cap\mathrm{sym}_\mathcal{G}(\dot{x})$ as in the statement of the lemma. In the ground model $M\vDash\mathrm{ZFC}$, let $\{p_\alpha\ |\ \alpha<\beta\}$ denote a wellordering of the poset $\mathbb{P}$. Inductively define the sequence $F_\alpha,$ for $\alpha\leq \beta$, as follows.

    Let $F_0=\dot{F}$. For $\alpha\leq\beta$, assume inductively that for every $\gamma<\alpha$, $K\preceq\mathrm{sym}_\mathcal{G}(\dot{F}_\gamma)\cap\mathrm{sym}_\mathcal{G}(\dot{x})$ and
    $$
    1\Vdash``\dot{F}_\gamma\subseteq2^{\dot{x}}\mathrm{\ has\ the\ FIP}."
    $$
    If $\alpha=\gamma+1$, use filter extension property and $\mathrm{AC}$ in $M$ to choose $q\leq p_\gamma$ and $y'\in\{\dot{y},\dot{x}\setminus\dot{y}\}$ so that
    $$
    1\Vdash``F_\gamma\cup\mathrm{Orb}(\langle q,y'\rangle,K)\subseteq2^{\dot{x}}\mathrm{\ has\ the\ FIP},"
    $$
    and set
    $$
    \dot{F}_\alpha=F_\gamma\cup\mathrm{Orb}(\langle q,y'\rangle,K).
    $$
    By the inductive hypothesis on $\dot{F}_\gamma$, the filter extension property can be applied to fix $\dot{F}_{\alpha}$ in the manner given above, in which case $\dot{F}_\alpha$ can been seen to also meet the inductive hypothesis.
    
    If $\alpha$ is a limit ordinal, set
    $$
    \dot{F}_\alpha=\bigcup_{\gamma<\alpha}\dot{F}_\gamma.
    $$
    By the inductive hypothesis for $\dot{F}_\gamma$, with $\gamma<\alpha$, and the fact that $\dot{F}_\gamma\subseteq\dot{F}_{\gamma'}$ when $\gamma<\gamma'$, we can see that $\dot{F}_\alpha$ also meets the inductive hypothesis.
    
    In this case, because we have added a dense set of conditions that decide between $\dot{y}$ and $\dot{x}\setminus\dot{y}$, $\dot{F}_{\dot{y}}=\dot{F}_\beta$ has the desired properties to satisfy the lemma.
\end{proof}

\begin{theorem}\label{fep implies bpi in sym ext}
    If $M\vDash\mathrm{ZFC}$ and $\langle\mathbb{P},\mathcal{G},\mathcal{F}\rangle$ has the filter extension property over $M$, then $\mathrm{BPI}$ holds in the corresponding symmetric extension of $M$.
\end{theorem}
\begin{proof}
    Note first that if $N$ is a symmetric extension based on the symmetric system $\langle\mathbb{P},\mathcal{G},\mathcal{F}\rangle$ and the $\mathbb{P}-$generic filter $G$, and if
    $$
    N\vDash``F\subseteq2^x\mathrm{\ has\ the\ FIP}"\wedge y\subseteq x,
    $$
    then there are names $\dot{F},\dot{y},\dot{x}\in\mathrm{HS}_\mathcal{F}$ for which
    $$
    1\Vdash``\dot{F}\subseteq2^{\dot{x}}\mathrm{\ has\ the\ FIP}"\wedge\dot{y}\subseteq\dot{x}
    $$
    and $\dot{F}^G=F,\dot{x}^G=x,$ and $\dot{y}^G=y$. Thus, the constraints given in the filter extension property are comprehensive for extending filters in symmetric extensions.

    Now, in the ground model $M\vDash\mathrm{ZFC}$, fix a wellordering $\{\dot{y}_\lambda\ |\ \lambda<\kappa\}$ of the names in $\mathrm{HS}_\mathcal{F}$ that are forced by $1\in\mathbb{P}$ to be subsets of $\dot{x}$. In a manner similar to the proof of lemma \ref{extend to include new subset}, one can inductively define a sequence $\dot{F}_\lambda$, for $\lambda\leq\kappa$, with $\dot{F}_0=\dot{F}$, and so that for every $\lambda\leq\kappa$, $K\preceq\mathrm{sym}_\mathcal{G}(\dot{F}_\lambda)\cap\mathrm{sym}_\mathcal{G}(\dot{x}),$
    $$
    1\Vdash``\dot{F}_\lambda\subseteq2^{\dot{x}}\mathrm{\ has\ the\ FIP},"
    $$
    and so that
    $$
    1\Vdash\dot{y}_\lambda\in\dot{F}_{\lambda+1}\vee\dot{x}\setminus\dot{y}_\lambda\in\dot{F}_{\lambda+1}.
    $$
    In this case, lemma \ref{extend to include new subset} justifies the ability to extend from $\dot{F}_\lambda$ to $\dot{F}_{\lambda+1}$ in such a manner; at limit stages, unions are taken. We can then see that for $\dot{U}=\dot{F}_\kappa$, $K\preceq\mathrm{sym}_\mathcal{G}(\dot{U})$ and so $\dot{U}\in\mathrm{HS}_\mathcal{F}$, and
    $$
    1\Vdash_{\mathrm{HS}_\mathcal{F}}\dot{U}\supseteq\dot{F}\mathrm{\ is\ an\ ultrafilter\ on\ }\dot{x}."
    $$    
\end{proof}

\begin{remark}
    One might wonder whether it is possible to prove theorem \ref{fep implies bpi in sym ext} while only assuming a ground model $M\vDash\mathrm{ZF}+\mathrm{BPI}$. This would not work; it is consistent that one can force with a poset $\mathbb{P}$ over a model $M\vDash\mathrm{ZF}+\mathrm{BPI}$ to add a filter that does not extend to any ultrafilter in the extension $M[G]$. Forcing with $\mathbb{P}$ can be expressed in terms of a symmetric system $\langle\mathbb{P},\mathcal{G},\mathcal{F}\rangle$ in which $\mathcal{G}$ and $\mathcal{F}$ are trivial, and which trivially has the filter extension property. In \cite{hall2013existence}, for instance, the authors force over the Cohen model (a known model of $\mathrm{ZF}+\mathrm{BPI}$), and identify a filter on $\omega$ in the resulting extension that does not extend to an ultrafilter. One can check that the main obstacle preventing a proof of theorem \ref{fep implies bpi in sym ext} using only $\mathrm{BPI}$, is that there are infinitely many options to take $q\leq p$ that witnesses an instance of the filter extension property.
\end{remark}

For the same reason as the case for permutation models, the use of Halpern's contradiction framework to prove $\mathrm{BPI}$ in symmetric extensions can be seen to subsume theorem \ref{fep implies bpi in sym ext}, and can be rephrased as a contradiction proof that filter extension property holds. Using the same approach as in theorem \ref{fep equiv to bpi in permutation models}, we can prove that the filter extension property is necessary for $\mathrm{BPI}$ in symmetric extensions, and thus that this connection is inevitable.

\begin{theorem}[Blass, \cite{blass1979injectivity}]\label{SVC in sym ext}
    Symmetric extensions satisfy $\mathrm{SVC}$.\footnote{Recall that we take symmetric extensions to be based on set-sized symmetric systems and ground models of $\mathrm{ZFC}$.}
\end{theorem}

\begin{proof}[Proof sketch of theorem 3.12]
    Let $\langle\mathbb{P},\mathcal{G},\mathcal{F}\rangle$ be a symmetric system. For $p\in\mathbb{P}$ and $H,K\in \mathcal{F}$, let $\dot{y}_{p,H,K}$ be a $\mathbb{P}-$name so that $\mathrm{sym}_\mathcal{G}(\dot{y}_{p,H,K})=H$ and for every distinct pair of cosets $\pi H,\tau H\in K/H$,
    \begin{equation}\label{disjoint forcing line}
    1\Vdash\pi \dot{y}_{p,H,K}\neq\tau \dot{y}_{p,H,K}.
    \end{equation}
    Let 
    $$
    \dot{s}=\bigcup_{p\in \mathbb{P},H,K\in\mathcal{F},H\preceq K}\mathrm{Orb}(\langle p,\dot{y}_{p,H,K}\rangle,K),\mathrm{\ and\ }\dot{S}=\bigcup_{\langle q,\dot{z}\rangle\in\dot{s}}\mathrm{Orb}(\langle q,\dot{z}\rangle,\mathcal{G}),
    $$
    so that $\mathrm{sym}_\mathcal{G}(S)=\mathcal{G}$, and thus $\dot{S}\in\mathrm{HS}_\mathcal{F}$. For every $\dot{x}\in\mathrm{HS}_\mathcal{F}$, we have the identity
    $$
    \dot{x}=\bigcup_{\langle p,\dot{z}\rangle\in\dot{x}}\mathrm{Orb}(\langle p,\dot{z}\rangle,\mathrm{sym}_\mathcal{G}(\dot{x})),
    $$
    which is a partition of $\dot{x}$ up to the repetition of orbits. Let $\alpha=|\{\mathrm{Orb}(\langle p,\dot{z}\rangle,\mathrm{sym}_\mathcal{G}(\dot{x}))\ |\ \langle p,\dot{z}\rangle\in\dot{x}\}|.$ For every $\langle p,\dot{z}\rangle\in\dot{x}$, there is a natural bijection from the name $\mathrm{Orb}(\langle p,\dot{y}_{p,H,K}\rangle,K)\subseteq \dot{S}$, with $H=\mathrm{sym}_\mathcal{G}(\dot{z})$ and $K=\mathrm{sym}_\mathcal{G}(\dot{x})$, to the name $\mathrm{Orb}(\langle p,\dot{z}\rangle,\mathrm{sym}_\mathcal{G}(\dot{x}))$. By line \ref{disjoint forcing line}, this bijection can be used to define a name $\dot{g}$ that is forced by 1 to be a surjection from $\mathrm{Orb}(\langle p,\dot{y}_{p,H,K}\rangle,K)$ onto $\mathrm{Orb}(\langle p,\dot{z}\rangle,\mathrm{sym}_\mathcal{G}(\dot{x}))$, with $\mathrm{sym}_\mathcal{G}(\dot{g})=\mathrm{sym}_\mathcal{G}(\dot{x})$. We can thus identify a name $\dot{f}\in\mathrm{HS}_\mathcal{F}$ for a surjection from a subset of $\dot{S}\times\check{\alpha}$ onto $\dot{x}$, with $\mathrm{sym}_\mathcal{G}(\dot{f})=\mathrm{sym}_\mathcal{G}(\dot{x})$, and so $1\Vdash_{\mathrm{HS}}\dot{S}\mathrm{\ satisfies\ }\mathrm{SVC}$.
\end{proof}

\begin{theorem}\label{BPI equivalent to FEP for sym ext}
    If $M\vDash\mathrm{ZFC}$, then the following are equivalent:
    \begin{enumerate}[label=(\roman*)]
        \item $\mathrm{BPI}$ holds in the symmetric extension of $M$ based on $\langle\mathbb{P},\mathcal{G},\mathcal{F}\rangle$,
        \item $\langle\mathbb{P},\mathcal{G},\mathcal{F}\rangle$ has the filter extension property.
    \end{enumerate}
\end{theorem}
\begin{proof}
    Part $(ii)\Rightarrow(i)$ is given in theorem \ref{fep implies bpi in sym ext}. Given theorem \ref{SVC in sym ext}, that symmetric extensions also model $\mathrm{SVC}$, the proof of $(i)\Rightarrow(ii)$ is the same argument as our proof of theorem \ref{fep equiv to bpi in permutation models}$(i)\Rightarrow(ii)$. Let $N$ be the symmetric extension of a ground model $M\vDash\mathrm{ZFC}$ based on the symmetric system $\langle\mathbb{P},\mathcal{G},\mathcal{F}\rangle$, and let $N\vDash\mathrm{BPI}$. By theorems \ref{Blass SVC BPI Theorem} and \ref{SVC in sym ext}, the first-order formula $\varphi(w_0,w_1,w_2,w_3)$ from our proof of theorem \ref{fep equiv to bpi in permutation models} part $(i)\Rightarrow(ii)$ can be used to define ultrafilters extending filters in $N$. Let $\dot{F}\in\mathrm{HS}_\mathcal{F}$ be a name for a filter on $\dot{x}\in\mathrm{HS}_\mathcal{F}$, let $\dot{S}\in\mathrm{HS}_\mathcal{F}$ be a name for a witness of $\mathrm{SVC}$ in $N$, let $\dot{U}'\in\mathrm{HS}_\mathcal{F}$ be a name for an ultrafilter witnessing \ref{Blass SVC BPI Theorem} in $N$, and let $\dot{f}\in\mathrm{HS}_\mathcal{F}$ be a name for a surjection from a subset of $\dot{S}\times\check{\alpha}$ onto the Boolean algebra $2^{\dot{x}}/\dot{F}$.\footnote{The correct formalizations of these names should be clear.} As in the proof skeetch of theorem \ref{SVC in sym ext}, take $\mathrm{sym}_\mathcal{G}(\dot{S})=\mathcal{G}$ and $\mathrm{sym}_\mathcal{G}(\dot{F})\cap\mathrm{sym}_\mathcal{G}(\dot{x})\preceq\mathrm{sym}_\mathcal{G}(2^{\dot{x}}/\dot{F})=\mathrm{sym}_\mathcal{G}(\dot{f})$. Then,
    $$
    \dot{U}=\{\langle p,\dot{y}\rangle|p\Vdash\varphi(\dot{y},\dot{f},\dot{S},\dot{U}')\}
    $$
    is a name so that
    $$
    1\Vdash_N``\dot{U}\supseteq\dot{F}\mathrm{\ is\ an\ ultrafilter\ on\ }\dot{x}."
    $$
    By the symmetry lemma \ref{symmetry lemma}, we have that
    $$
    \mathrm{sym}_\mathcal{G}(\dot{f})\cap\mathrm{sym}_\mathcal{G}(\dot{S})\cap\mathrm{sym}_\mathcal{G}(\dot{U}')\preceq\mathrm{sym}_\mathcal{G}(\dot{U}),
    $$
    from which it follows that
    $$
    \mathrm{sym}_\mathcal{G}(\dot{F})\cap\mathrm{sym}_\mathcal{G}(\dot{x})\cap\mathrm{sym}_\mathcal{G}(\dot{U}')\preceq\mathrm{sym}_\mathcal{G}(\dot{U}).
    $$

    We claim again that, for $K=\mathrm{sym}_\mathcal{G}(\dot{U}')$, the base of subgroups $\{H\cap K\ |\ H\in \mathcal{F}\}$ of $\mathcal{F}$ witnesses the filter extension property. Let $H\in\mathcal{F}$, and consider $\dot{F},\dot{x},\dot{y}$ meeting the hypotheses of the filter extension property, with $H\cap K\preceq\mathrm{sym}_\mathcal{G}(\dot{F})\cap\mathrm{sym}_\mathcal{G}(\dot{x})$. Let $\dot{U}$ be a name for an ultrafilter on $\dot{x}$ extending $\dot{F}$, defined by $\varphi$ as above. For every $p\in\mathbb{P}$, there is a $q\leq p$ and a $y'\in\{\dot{y},\dot{x}\setminus\dot{y}\}$ so that
    $$
    \dot{F}\cup\mathrm{Orb}(\langle q,y'\rangle,\mathrm{sym}_\mathcal{G}(\dot{U}))\subseteq\dot{U}.
    $$
    Because $H\cap K\preceq\mathrm{sym}_\mathcal{G}(\dot{F})\cap\mathrm{sym}_\mathcal{G}(\dot{x})\cap\mathrm{sym}_\mathcal{G}(\dot{U}')\preceq\mathrm{sym}_\mathcal{G}(\dot{U})$, in particular
    $$
    1\Vdash``\dot{F}\cup\mathrm{Orb}(\langle q,y'\rangle,H\cap K)\mathrm{\ has\ the\ FIP}."
    $$
    Thus, $H\cap K$ witnesses the filter extension property, which completes the proof.
\end{proof}

We can now turn our attention to developing direct methods to prove the filter extension property for permutation models and symmetric systems. In the next section, this is done for the symmetric system $\langle\mathbb{P}(I,Q),\mathcal{G},\mathcal{F}\rangle$ for the generalized Cohen model, when $I$ is sufficiently large, by adapting Harrington's proof of the Halpern-L\"auchli theorem.

\section{$\mathrm{BPI}$ in the Generalized Cohen Model $N(I,Q)$, When $I$ is Large}
In the previous section, we noted that the known contradiction proofs of $\mathrm{BPI}$ in models without choice can be rephrased as contradiction proofs of the filter extension property. In \cite{halpern1971boolean}, a contradiction proof was given for that $\mathrm{BPI}$ holds in the Cohen model; in \cite{stefanovic2023alternatives}, it was shown that this contradiction proof can be modified to prove $\mathrm{BPI}$ in the generalized Cohen model $N(I,Q)$. In terms of the filter extension property, the proofs from \cite{halpern1971boolean} and \cite{stefanovic2023alternatives} can be restated as proofs of the following theorem and corollary.

\begin{theorem}[Hapern, L\'evy, \cite{halpern1971boolean}, Stefanovi\'c, \cite{stefanovic2023alternatives}]\label{Halpern-Levy FEP Theorem}

    Let $\langle\mathrm{P}(I,Q),\mathcal{G},\mathcal{F}\rangle$ be a symmetric system for a Generalized Cohen model and let $\dot{F},\dot{x},\dot{y}\in \mathrm{HS}_\mathcal{F}$, with $\mathrm{supp}(\dot{F})=F$, $\mathrm{supp}(\dot{x})=X$, and $K=\mathrm{fix}_\mathcal{G}(F\cup X)$, such that
    $$
    1\Vdash``\dot{F}\subseteq2^{\dot{x}}\mathrm{\ has\ the\ FIP}"\wedge\dot{y}\subseteq\dot{x}.
    $$
    Then, for every $p\in\mathbb{P}(I,Q)$, there is a $q\leq p$ and a $y'\in\{\dot{y},\dot{x}\setminus\dot{y}\}$ so that
    $$
    1\Vdash``\dot{F}\cup\mathrm{Orb}(\langle q,y'\rangle,K)\mathrm{\ has\ the\ FIP}."
    $$
\end{theorem}

\begin{corollary}[Halpern, L\'evy, \cite{halpern1971boolean}, Stefanovi\'c, \cite{stefanovic2023alternatives}]\label{N(I,Q) models BPI from FEP}
    $N(I,Q)\vDash\mathrm{BPI}$.
\end{corollary}
\begin{proof}
    Theorem \ref{Halpern-Levy FEP Theorem} shows that the symmetric system $\langle\mathbb{P}(I,Q),\mathcal{G},\mathcal{F}\rangle$ has the filter extension property, and specifically that it is witnessed by the base of subgroups $\{\mathrm{fix}_\mathcal{G}(E)\ |\ E\in[I]^{<\omega}\}$ of $\mathcal{F}$.

    To see why, let $K=\mathrm{fix}_\mathcal{G}(E)$, let $\dot{F},\dot{x},\dot{y}\in\mathrm{HS}_\mathcal{F}$ with $K\preceq\mathrm{sym}_\mathcal{G}(\dot{F})\cap\mathrm{sym}_\mathcal{G}(\dot{x})$ and
    $$
    1\Vdash``\dot{F}\subseteq2^{\dot{x}}\mathrm{\ has\ the\ FIP}"\wedge\dot{y}\subseteq\dot{x},
    $$
    and let $p\in\mathbb{P}(I,Q)$. By theorem \ref{Halpern-Levy FEP Theorem}, for $\mathrm{supp}(\dot{F})=F$ and $\mathrm{supp}(\dot{x})=X$, there is a $q\leq p$ and $y'\in\{\dot{y},\dot{x}\setminus\dot{y}\}$ for which
    $$
    1\Vdash``\dot{F}\cup\mathrm{Orb}(\langle q,y'\rangle,\mathrm{fix}_\mathcal{G}(F\cup X))\mathrm{\ has\ the\ FIP}."
    $$
    By lemma \ref{Minimal Supports Lemma}, because $K=\mathrm{fix}_\mathcal{G}(E)\preceq\mathrm{sym}_\mathcal{G}(\dot{F})\cap\mathrm{sym}_\mathcal{G}(\dot{x})$, we can conclude that $X\subseteq E$, $F\cap X\subseteq E$, and $K=\mathrm{fix}_\mathcal{G}(E)\preceq\mathrm{fix}_\mathcal{G}(F\cup X)$. In this case, the line above allows us to conclude that
    $$
    1\Vdash``\dot{F}\cup\mathrm{Orb}(\langle q,y'\rangle,K)\mathrm{\ has\ the\ FIP}."
    $$
\end{proof}

The main objective of this section is to directly prove theorem \ref{Halpern-Levy FEP Theorem}, in the case that $I$ is taken to be sufficiently large. We do so by adapting Harrington's proof of the Halpern-L\"auchli theorem.

\subsection{Adapting Harrington's Proof}
In Harrington's proof, given a level product $\otimes_{i<d}T$ and a coloring $c:\otimes_{i<d}T\to\{\mathrm{green,\ red}\}$, one searches through a $\mathbb{P}(\theta,T)$-name $\dot{U}$ for an ultrafilter on $\omega$, and takes $\theta$ large enough to identify a node $\vec{t}\in\otimes_{i<d}T$ and a color $c'\in\{\mathrm{green,\ red}\}$ that are regulated by $\dot{U}$ to witness the Halpern-L\"auchli theorem for $\otimes_{i<d}T$ and $f$.

Suppose that $\dot{F},\dot{x},$ and $\dot{y}$ meet the hypotheses of theorem \ref{Halpern-Levy FEP Theorem}. The main idea of this subsection is that one can search through a $\mathbb{P}(I,Q)$-name $\dot{U}$ for an ultrafilter on $\dot{x}$ that extends $\dot{F}$, and take $I$ large enough to identify a condition $q\in\mathbb{P}(\theta,Q)$ and a $y'\in\{\dot{y},\dot{x}\setminus\dot{y}\}$, so that $\dot{U}$ can regulate that
$$
1\Vdash``\dot{F}\cup\mathrm{Orb}(\langle q,y'\rangle,\mathrm{fix}_\mathcal{G}(F\cup X))\subseteq2^{\dot{x}}\mathrm{\ has\ the\ FIP}."
$$
This is made precise in our proof of theorem \ref{N(theta,Q) Models BPI}. For the rest of the subsection, we will write $\theta$ instead of $I$ to match the notation from Harrington's proof, and to indicate that we wish to take $I=\theta$ to be sufficiently large.

The lemmas and definitions preceding theorem \ref{N(theta,Q) Models BPI} establish several basic combinatorial features that are needed in our adaptation of Harrington's proof. Lemma \ref{Product Ramsey Theorem}, definition \ref{Type Definition}, and lemma \ref{Compatible Type Ramsey Lemma} can be seen to be simple modifications made to the exposition preceding Harrington's proof in \cite{todorchevich1995some}. Note that our definition of a condition ``type" differs slightly from theirs, and that our lemma \ref{Compatible Type Ramsey Lemma} is essentially their lemma 6.8. See Appendix $\mathrm{C}$ for the proof of lemma \ref{Compatible Type Ramsey Lemma} as a corollary of lemma 6.8 from \cite{todorchevich1995some}.

\begin{lemma}\label{Product Ramsey Theorem}
    For every pair of cardinals $\kappa,\lambda$, there is a cardinal $\theta$ so that for every $d\in\omega$, and every function of the form
    $$
    c:\theta^d\to\kappa,
    $$
    there are disjoint sets $\mathcal{H}_0,...,\mathcal{H}_{d-1}\subseteq\theta$ of cardinality $\lambda$ with $f$ is constant on $\prod_{i<d}\mathcal{H}_i$.
\end{lemma}
\begin{proof}
This is essentially due to the Erd\H{o}s-Rado theorem, as stated in line 17.32 of \cite{erdos2011combinatorial}. On page 280 of \cite{todorcevic2007walks}, the cardinals $\Theta_d$ are defined, which satisfy the lemma for fixed $d$, in the case that $\kappa=\lambda=\omega$. It is noted in the remark on page 52 of \cite{todorchevich1995some} that the Erd\H{o}s-Rado theorem implies that $\Theta_d\leq \beth_{d-1}^+$. This idea is then applied in theorem 11 of \cite{stefanovic2023alternatives} to prove that there is a $\theta_d$ satisfying the lemma for arbitrary $\kappa,\lambda$ and fixed $d$; this is a pidgeonhole proof that relies on the fact that powersets grow in cardinality. Taking $\theta=\bigcup_{d<\omega}\theta_d$ then satisfies the lemma.
\end{proof}

\begin{definition}\label{Type Definition}
    Two conditions $p,q\in\mathbb{P}(\theta,Q)$ are said to be of the same \textit{type} if there is an order-preserving map $\pi$ from $\mathrm{dom}(p)$ to $\mathrm{dom}(q)$, so that $p(\xi)=q(e(\xi))$. We use the term \textit{type} to refer to an equivalence class in $\mathbb{P}(\theta,Q)$ under this relation.
\end{definition}

\begin{lemma}\label{Compatible Type Ramsey Lemma}
    For every type $t$ of a condition in $\mathbb{P}(\theta,Q)$ and integers $m$ and $d$, there is an integer $M=M(t,m,d)$ such that for every set $\{p_{\vec{\alpha}}\ |\ \vec{\alpha}\in M^d\}$ of elements of $\mathbb{P}(\theta,Q)$ of type $t$, there exist $H_i\subseteq M$ (for $i<d$) of size $m$, so that the conditions in $\{p_{\vec{\alpha}}\ |\ \vec{\alpha}\in \prod_{i<d}H_i\}$ are pairwise compatible.
\end{lemma}
    \begin{proof}
    See Appendix $\mathrm{C}$, where this is derived from the analogous lemma 6.8 of \cite{todorchevich1995some}.
\end{proof}

When choosing a $q\leq p$ to witness theorem \ref{Halpern-Levy FEP Theorem}, one can strengthen $q$ to refine which images of $y'\in\{\dot{y},\dot{x}\setminus\dot{y}\}$ can mutually appear in an interpretation of the name $\mathrm{Orb}(\langle q,y'\rangle, K)$. If $q(i)\perp q(j)$ for some $i,j\in\mathrm{supp}(q)$, then for $\pi=(i\ j)$, $q$ and $\pi q$ are incompatible and $y',\pi y'$ will never simultaneously appear in an interpretation of $\mathrm{Orb}(\langle q,y'\rangle, K)$. We introduce the definition below to capture this important property.

\begin{definition}\label{Distinguishing Supports Definition}
    A condition $p\in\mathbb{P}(I,Q)$ \textit{distinguishes its support} if for every distinct $i,j\in\mathrm{supp}(p)$, $p(i)\perp p(j)$.
\end{definition}

Taking conditions to distinguish their supports will also help to manage certain restrictions and images of conditions that are relevant to our proof; lemmas \ref{Distinguished Restriction Lemma} and \ref{Disjoint J_i Lemma} express these properties.

\begin{lemma}\label{Distinguished Restriction Lemma}
    Let $r\in\mathbb{P}(\theta,Q)$, $K\in\mathcal{F}$, and $\pi_0,...,\pi_n\in K$. Suppose that the conditions $r_0,...,r_n$ are mutually compatible, distinguish their supports, and that for every $i\leq n$, and $r_i\upharpoonright_{\mathrm{supp}(\pi_ir)}=\pi_ir$. Then $\bigcup_{i\leq n}r_i$ restricts to $\bigcup_{i\leq n}\pi_ir$ over the latter condition's support.
\end{lemma}

\begin{proof}
    For each $i\leq n$, we have that $r_i\leq \pi_ir$, and so
    $$
    \bigcup_{i\leq n}r_i\leq \bigcup_{i\leq n}\pi_ir.
    $$
    Suppose for the sake of contradiction that the lemma were false. By the above line, this would imply that there are $i,j\leq n$ and $k\in \theta$ for which $k\in\mathrm{supp}(r_i)\cap\mathrm{supp}(\pi_jr)$ and $r_i(k)\lneq\pi_jr(k)$. Because $r_i$ and $r_j$ distinguish their supports, so do their restrictions $\pi_ir$ and $\pi_jr$; in this case, $\mathrm{rng}(r_i),\mathrm{rng}(\pi_ir),$ and $\mathrm{rng}(\pi_jr)$ are antichains in $Q$. Moreover, since $\pi_j\pi_i^{-1}(\pi_ir)=\pi_jr$ and since $r_i$ restricts to $\pi_ir$, we have that $\mathrm{rng}(\pi_jr)=\mathrm{rng}(\pi_ir)\subseteq\mathrm{rng}(r_i)$. If $r_i(k)\leq \pi_jr(k)$, these two conditions are compatible in $Q$. This can only happen if $r_i(k)$ and $\pi_jr(k)$ each take the same value in the antichain $\mathrm{rng}(r_i)$, in which case $r_i(k)=\pi_jr(k)$. It is thus not possible that $r_i(k)\lneq\pi_jr(k)$.
    
\end{proof}

\begin{lemma}\label{Disjoint J_i Lemma}
    Suppose that $q\in\mathbb{P}(\theta,Q)$ distinguishes its support, and that for $\delta\in[\mathcal{G}]^{<\omega}$, the conditions in the set $\{\pi q\ |\ \pi\in\delta\}$ are mutually compatible. Then for $i\in\mathrm{supp}(q)$, the sets
    $$
    J_i=\{\pi(i)\ |\ \pi\in\delta\}\subseteq\theta
    $$
    are mutually disjoint.
\end{lemma}
\begin{proof}
    For each $\pi(i)\in J_i$, $\pi q(\pi(i))=q(i)$ -- in this sense, $J_i$ collects the indices over which the conditions $\pi q$, for $\pi\in \delta$, copy the output of $q(i)$. Thus, if $k\in J_i\cap J_{j}$, with $i\neq j$, then for some $\pi,\tau\in\delta$ we have that $\pi q(k)=q(i)$ and $\tau q(k)=q(j)$. Because $\pi q$ and $\tau q$ are compatible, it must be the case that $\pi q(k)=q(i)$ and $\tau q(k)=q(j)$ are compatible in $Q$. This contradicts the fact that $q$ distinguishes its support, and so $J_i\cap J_{j}=\emptyset$.
\end{proof}

We can now prove the main theorem of this section. For the sake of comparison, the language and formatting of our proof are modeled after the exposition of Harrington's proof from \cite{todorchevich1995some}.

\begin{theorem}[Halpern, L\'evy, \cite{halpern1971boolean}, Stefanovi\'c, \cite{stefanovic2023alternatives}]\label{N(theta,Q) Models BPI}
    Theorem 4.1 holds for $\langle\mathbb{P}(\theta,Q),\mathcal{G},\mathcal{F}\rangle$ when $\theta$ is large enough to satisfy lemma \ref{Product Ramsey Theorem}, with $\kappa=|Q|$.
\end{theorem}
\begin{proof}
Let $\dot{F},\dot{x},\dot{y}\in\mathrm{HS}_\mathcal{F}$, $K=\mathrm{fix}_\mathcal{G}(F\cup X),$ and $p\in\mathbb{P}(\theta,Q)$ be as in the statement of theorem 4.1, with $F,X,Y\in [\theta]^{<\omega}$ the supports of $\dot{F},\dot{x},$ and $\dot{y}$. Without loss of generality, we may assume that $\mathrm{supp}(p)\subseteq F\cup X\cup Y$. To see why, note that for $p'=p\upharpoonright_{F\cup X\cup Y}$, our proof will produce witnesses $q'\leq p'$ and $y'\in\{\dot{y},\dot{x}\setminus\dot{y}\}$, with $\mathrm{supp}(q')=F\cup X\cup Y$, so that
$$
1\Vdash``\dot{F}\cup\mathrm{Orb}(\langle q',y'\rangle,K)\subseteq2^{\dot{x}}\mathrm{\ has\ the\ FIP}."
$$
In this case, $q=q'\cup p$ is well-defined, and
$$
1\Vdash \mathrm{Orb}(\langle q,y'\rangle,K)\subseteq\mathrm{Orb}(\langle q',y'\rangle,K).
$$
As such, $q\leq p$ and
$$
1\Vdash``\dot{F}\cup\mathrm{Orb}(\langle q',y'\rangle,K)\subseteq2^{\dot{x}}\mathrm{\ has\ the\ FIP}."
$$
We can thus assume without loss of generality that $\mathrm{supp}(p)\subseteq F\cup X\cup Y$, and produce witnesses $q\leq p$ and $y'\in\{\dot{y},\dot{x}\setminus\dot{y}\}$ to theorem 4.1, with $\mathrm{supp}(q)=F\cup X\cup Y$. This will simplify the calculation of supports throughout the proof.

For $i\in Y\setminus(F\cup X)$, let $\theta_i\subseteq\theta\setminus(F\cup X)$ be pairwise disjoint sets of cardinality $\theta$. Recall from definition \ref{pi_vec{alpha}} that for
$$
\vec{\alpha}=\{\langle i,\alpha_i\rangle|i\in Y\setminus(F\cup X)\}\in\prod_{i\in Y\setminus(F\cup X)}\theta_i
$$
we define
$$
\pi_{\vec{\alpha}}=\prod_{i\in Y\setminus(F\cup X)}(i\ \alpha_i)\in K=\mathrm{fix}_\mathcal{G}(F\cup X).
$$
The membership $\pi_{\vec{\alpha}}\in K$ follows from the fact that each transposition $(i\ \alpha_i)$ cycles $i\in Y\setminus(F\cup X)$ with $\alpha_i\in \theta_i\subseteq\theta\setminus(F\cup X)$, and thus fixes $F\cup X$.

Let $\dot{U}$ be a (not necessarily symmetric) $\mathbb{P}(\theta,Q)$-name for an ultrafilter on $\dot{x}$ extending $\dot{F}$. For every $\vec{\alpha}\in \prod_{i\in Y\setminus(F\cup X)}\theta_i$, either
\begin{enumerate}
    \item $\pi_{\vec{\alpha}} p\Vdash\pi_{\vec{\alpha}}({\dot{x}}\setminus{\dot{y}})\in\dot{U}$, or
    \item there is some $q_{\vec{\alpha}}\leq p$ so that $\pi_{\vec{\alpha}}q_{\vec{\alpha}}\Vdash\pi_{\vec{\alpha}}{\dot{y}}\in\dot{U}$.
\end{enumerate}

In the second case, extend $q_{\vec{\alpha}}$ if necessary so that $F\cup X\cup Y\subseteq \mathrm{supp}(q_{\vec{\alpha}})$ and so that $q_{\vec{\alpha}}$ distinguishes its support. Recall that we assume $\theta$ to be large enough to satisfy lemma \ref{Product Ramsey Theorem} for $\kappa=|Q|$. Then, for $i\in Y\setminus(F\cup X)$, there are disjoint infinite sets $\mathcal{H}_i\subseteq\theta_i$ so that either
\begin{enumerate}
    \item[(a)] the first alternative above happens for all $\vec{\alpha}\in \prod_{i\in Y\setminus(F\cup X)}\mathcal{H}_i$, or
    \item[(b)] there is a type $t$ and a condition $q$ with $\mathrm{supp}(q)=F\cup X\cup Y$ that distinguishes its support, so that for all $\vec{\alpha}\in \prod_{i\in Y\setminus(F\cup X)}\mathcal{H}_i$, $\pi_{\vec{\alpha}}q_{\vec{\alpha}}$ satisfies the second alternative above, is of type $t$, and we have the restriction $q_{\vec{\alpha}}|_{F\cup X\cup Y}=q$.\footnote{It is straightforward to verify in $M$ that there are $\kappa=|Q|$ many possibilities specified in the outcome $(b)$.}
\end{enumerate}
Note that $q\leq p$, as is required by theorem \ref{Halpern-Levy FEP Theorem}. Without loss of generality, we assume the second case and argue that
\begin{equation}\label{FEP conclusion, Main Proof}
1\Vdash``\dot{F}\cup \mathrm{Orb}(\langle q,\dot{y}\rangle,K)\subseteq2^{\dot{x}}\mathrm{\mathrm{\ has\ the\ FIP}.}"
\end{equation}

By lemma \ref{Compatible Type Ramsey Lemma}, for every $m\in\omega$ and $i\in Y\setminus(F\cup X)$, we can find subsets $H^m_i\subseteq\mathcal{H}_i$ of cardinality $m$, so that for $\vec{\alpha}\in\prod_{i\in Y\setminus(F\cup X)}H^m_i$, the conditions $\pi_{\vec{\alpha}}q_{\vec{\alpha}}$ are mutually compatible. We claim that for every $S\subseteq\prod_{i\in Y\setminus(F\cup X)}H^m_i$,
\begin{equation}\label{Search Conclusion, Main Proof}
    1\Vdash``\dot{F}\cup\{\langle\pi_{\vec{\alpha}}q,\pi_{\vec{\alpha}}\dot{y}\rangle|\vec{\alpha}\in S\}\subseteq2^{\dot{x}}\mathrm{\ has\ the\ FIP}."
\end{equation}
Because the conditions $q_{\vec{\alpha}}$, for $\vec{\alpha}\in S$ satisfy case $(b)$ above, we have that
$$
\bigcup_{\vec{\alpha}\in S}\pi_{\vec{\alpha}}q_{\vec{\alpha}}\Vdash(\forall\vec{\alpha}\in S)\pi_{\vec{\alpha}}\dot{y}\in\dot{U}.
$$
Because $1\Vdash\dot{F}\subseteq\dot{U}$, we can conclude in particular that
\begin{equation}\label{Forced Search Conclusion, Main Proof}
\bigcup_{\vec{\alpha}\in S}\pi_{\vec{\alpha}}q_{\vec{\alpha}}\Vdash``\dot{F}\cup\{\langle1,\pi_{\vec{\alpha}}\dot{y}\rangle|\vec{\alpha}\in S\}\subseteq2^{\dot{x}}\mathrm{\ has\ the\ FIP}."
\end{equation}
We would like to use lemma \ref{Support Restriction Lemma} to restrict the condition in line (\ref{Forced Search Conclusion, Main Proof}) to the support of the names that appear in the statement being forced. We can first calculate this union of supports to be
$$
    \mathrm{supp}(\dot{F})\cup\mathrm{supp}(\dot{x})\cup\bigcup_{\vec{\alpha}\in S}\mathrm{supp}(\pi_{\vec{\alpha}}\dot{y})=
$$
$$
F\cup X\cup \bigcup_{\vec{\alpha}\in S}\pi_{\vec{\alpha}}Y=
$$
$$
\bigcup_{\vec{\alpha}\in S}\mathrm{supp}(\pi_{\vec{\alpha}}q)=
$$
$$
    \mathrm{supp}(\bigcup_{\vec{\alpha}\in S}\pi_{\vec{\alpha}}q).
$$
One can verify that the hypothesis of lemma \ref{Distinguished Restriction Lemma} are met when the conditions $r_0,...,r_n$ are replaced with $\pi_{\vec{\alpha}}q_{\vec{\alpha}}$, for $\vec{\alpha}\in S$, and the conditions $\pi_0r,...,\pi_nr$ are replaced with $\pi_{\vec{\alpha}}q$, for $\vec{\alpha}\in S$. In this case, one can apply the lemma to conclude that $\bigcup_{\vec{\alpha}\in S}\pi_{\vec{\alpha}}q_{\vec{\alpha}}$ restricts to $\bigcup_{\vec{\alpha}\in S}\pi_{\vec{\alpha}}q$ over the latter condition's support. Given that $\mathrm{supp}(\bigcup_{\vec{\alpha}\in S}\pi_{\vec{\alpha}}q)$ is the support of the names in line (\ref{Forced Search Conclusion, Main Proof}), we can use lemma \ref{Support Restriction Lemma} to conclude that
\begin{equation}\label{Refined Forcing Search Conclusion, Main Proof}
\bigcup_{\vec{\alpha}\in S}\pi_{\vec{\alpha}}q\Vdash``\dot{F}\cup\{\langle1,\pi_{\vec{\alpha}}\dot{y}\rangle|\vec{\alpha}\in S\}\subseteq2^{\dot{x}}\mathrm{\ has\ the\ FIP}."   
\end{equation}
Because $S\subseteq\prod_{i\in Y\setminus(F\cup X)}H^m_i$ was arbitrary, line (\ref{Refined Forcing Search Conclusion, Main Proof}) holds for every $S'\subseteq S$, which allows us to conclude line (\ref{Search Conclusion, Main Proof}).

To complete our proof of line (\ref{FEP conclusion, Main Proof}), let $\Delta\in[\mathrm{Orb}(\langle q,\dot{y}\rangle,K)]^{<\omega}$ be arbitrary, and take $\delta\in[K]^{<\omega}$ so that $\Delta=\{\langle\pi q,\pi\dot{y}\rangle\ |\ \pi\in\delta\}$. We may assume without loss of generality that the conditions in the set $\{\pi q\ |\ \pi\in\delta\}$ are mutually compatible, so that the names in $\mathrm{rng}(\Delta)$ can mutually appear in an interpretation of $\mathrm{Orb}(\langle q,\dot{y}\rangle,K)$. If we can derive that
\begin{equation}\label{Shift Conclusion, Main Proof}
1\Vdash``\dot{F}\cup\Delta\subseteq2^{\dot{x}}\mathrm{\ has\ the\ FIP},"
\end{equation}
then \ref{FEP conclusion, Main Proof} follows from the fact that $\Delta\in[\mathrm{Orb}(\langle q,\dot{y}\rangle,K)]^{<\omega}$ was taken to be arbitrary.

We will show that, essentially, the lines (\ref{Search Conclusion, Main Proof}) and (\ref{Shift Conclusion, Main Proof}) are the same up to an application of the symmetry lemma by some $\sigma\in K=\mathrm{fix}_\mathcal{G}(F\cup X)$. We claim that $m\in\omega$ can be taken large enough to define a $\sigma\in K$, for which
\begin{equation}\label{shift by sigma}
\sigma\delta=\{\sigma\pi\ |\ \pi\in\delta\}\subseteq\{\pi_{\vec{\alpha}}|\vec{\alpha}\in\prod_{i\in Y\setminus(F\cup X)}H^m_i\}.
\end{equation}
If this can be done, let $S=\sigma\delta$ and apply \ref{Search Conclusion, Main Proof} to get
\begin{equation}\label{pre-shifted conclusion, main proof}
1\Vdash``\dot{F}\cup\sigma\Delta\subseteq2^{\dot{x}}\mathrm{\ has\ the\ FIP}."
\end{equation}
Because $\sigma\in K=\mathrm{fix}_\mathcal{G}(F\cup X)$, we can apply the symmetry lemma \ref{symmetry lemma} with $\sigma^{-1}$ to line (\ref{pre-shifted conclusion, main proof}), while fixing $\dot{F}$ and $\dot{x}$, to get the desired statement in line (\ref{Shift Conclusion, Main Proof}), which completes the proof.

We will now justify the existence of such a $\sigma\in K$. Fix an enumeration $\delta=\{\pi_j\ |\ j\leq n\}$, and for $j\leq n$, define
$$
\vec{\beta}_j=\{\langle i,\pi_j(i)\rangle|i\in Y\setminus(F\cup X)\},
$$
so that $\pi_j\upharpoonright_{F\cup X\cup Y}=\pi_{\vec{\beta}_j}$. Because $\mathrm{supp}(\dot{y})\subseteq\mathrm{supp}(q)= F\cup X\cup Y$, lemma \ref{image support lemma} shows that $\pi_jq=\pi_{\vec{\beta}_j}q$ and $\pi_j\dot{y}=\pi_{\vec{\beta}_j}\dot{y}$. Without loss of generality, we can thus assume that $\delta=\{\pi_{\vec{\beta}_j}\ |\ j\leq n\}$. For $i\in Y\setminus(F\cup X)$, define the sets
$$
J_i=\{\pi(i)|\pi\in\delta\}.
$$
Because we assume that the conditions in $\{\pi q\ |\ \pi\in\delta\}$ are mutually compatible, by lemma \ref{Disjoint J_i Lemma} the sets $J_i$, for $i\in Y\setminus(F\cup X)$, are mutually disjoint. For every $i\in Y\setminus(F\cup X)$, because every $\pi\in\delta\subseteq K=\mathrm{fix}_\mathcal{G}(F\cup X)$, $J_i$ is also disjoint from $F\cup X$. For every $j\leq n$, $\vec{\beta}_j\in\prod_{i\in Y\setminus(F\cup X)}J_i$, which gives the inclusion
$$
\delta=\{\pi_{\vec{\beta}_j}|j\leq n\}\subseteq\{\pi_{\vec{\gamma}}|\vec{\gamma}\in\prod_{i\in Y(\setminus(F\cup X)}J_i\}.
$$

Now let $m=\mathrm{max}_{i\in Y\setminus(F\cup X)}|J_i|$ and for $i\in Y\setminus(F\cup X)$, consider the sets $H^m_i\subseteq\mathcal{H}_i\subseteq\theta_i\subseteq\theta\setminus(F\cup X)$ defined above. For $i\in Y\setminus(F\cup X)$, let $f_i:J_i\to H^m_i$ be an injection. Define $\sigma\in K=\mathrm{fix}_\mathcal{G}(F\cup X)$ by
$$
\sigma=\prod_{i\in Y\setminus(F\cup X)}\prod_{a\in J_i}(a\ f_i(a)).
$$
Because the sets $J_i$ are mutually disjoint, and the functions $f_i$ are injections into mutually disjoint codomains, $\sigma$ is well-defined as an automorphism of the index set $\theta$, an in particular it is an injection from $\prod_{i\in Y\setminus(F\cup X)}J_i$ into $\prod_{i\in Y\setminus(F\cup X)}H^m_i$. Moreover, for $i\in Y\setminus(F\cup X)$, the sets $J_i$ and $H^m_i$ are disjoint from $F\cup X$, so $\sigma\in K=\mathrm{fix}_\mathcal{G}(F\cup X)$. We can then conclude that
$$
\sigma\delta=\{\sigma\pi_j|j\leq n\}=\{\sigma\pi_{\vec{\beta}_j}|j\leq n\}\subseteq\{\sigma\pi_{\vec{\gamma}}|\vec{\gamma}\in\prod_{i\in Y(\setminus(F\cup X)}J_i\}\subseteq\{\pi_{\vec{\alpha}}|\vec{\alpha}\in\prod_{i\in Y\setminus(F\cup X)}H^m_i\},
$$
which satisfies line (\ref{shift by sigma}) and completes the proof.
\end{proof}

\begin{remark}\label{search and shift remark}
    We call this style of proof a ``search and shift" argument, since we search through \textit{any} name $\dot{U}$ for an ultrafilter on $\dot{x}$ that extends $\dot{F}$, then shift the outcomes of this search (stated in line (\ref{Search Conclusion, Main Proof})) to verify that for every relevant $\Delta\in[\mathrm{Orb}(\langle q,y'\rangle,K)]^{<\omega}$,
    $$
    1\Vdash``\dot{F}\cup\Delta\subseteq2^{\dot{x}}\mathrm{\ has\ the\ FIP},"
    $$
    In this sense, $\mathrm{BPI}$ in $N(\theta,Q)$ can be seen as inherited from the outer model $M[G]\vDash\mathrm{ZFC}$, with arbitrary names for ultrafilters coordinating the extension of names for symmetric filters.
\end{remark}

\section{$\mathrm{BPI}$ in the Generalized Cohen Model $N(I,Q)$, When $I$ is Small}
In this section, we use theorem \ref{N(theta,Q) Models BPI} to prove the general case of theorem \ref{Halpern-Levy FEP Theorem}. To do so, we will expand upon corollary \ref{Karagila-Schlicht Cohen Index Theorem}, originally from \cite{karagila2020have}. To state the corollary in full context, we cite the following definition and theorem.

\begin{definition}[Karagila, Schlicht, \cite{karagila2020have}]
    For sets $X$ and $Y$, let $\mathrm{Col}_{\mathrm{inj}}(X,Y)$ be the poset of well-orderable partial injections from $X$ into $Y$, ordered by reverse-inclusion.
\end{definition}

\begin{theorem}[Karagila, Schlicht, \cite{karagila2020have}]\label{KS G,F Theorem}
    Let $M\vDash\mathrm{ZFC}$, and $I,J\in M$ be infinite. Let $G$ be $\mathbb{P}(I,2^{<\omega})$-generic, and let $A\in N(I,2^{<\omega})$ be the standard Dedekind-finite set of Cohen-generic reals. If $F$ is $\mathrm{Coll}_{\mathrm{inj}}(J,A)$-generic over $N^G(I,Q)$, then $F$ is $\mathbb{P}(J,2^{<\omega})$-generic over $M$.
\end{theorem}

\begin{corollary}[Karagila, Schlicht, \cite{karagila2020have}]\label{Karagila-Schlicht Cohen Index Theorem}
Let $M\vDash\mathrm{ZFC}$ and let $I,J\in M$ be infinite. For every $\mathbb{P}(I,2^{<\omega})$-generic $G$, there is a $\mathbb{P}(J,2^{<\omega})$-generic $F$ (as in theorem \ref{KS G,F Theorem}), so that
$$
N^G(I,2^{<\omega})=N^F(J,2^{<\omega}).
$$
\end{corollary}
\noindent One can apply corollary \ref{Karagila-Schlicht Cohen Index Theorem} to see that the index $I$ in $N(I,2^{<\omega})$ can be taken to be arbitrarily large, without loss of generality. From theorem \ref{N(theta,Q) Models BPI} and corollary \ref{Karagila-Schlicht Cohen Index Theorem}, it immediately follows that $N(I,2^{<\omega})\vDash\mathrm{BPI}$, for every $I$.

We can sketch a quick, informal account of theorem \ref{KS G,F Theorem} and corollary \ref{Karagila-Schlicht Cohen Index Theorem}. One can view the generics $F$ and $G$ from theorem \ref{KS G,F Theorem} as essentially containing two pieces of information: a new set $A$ of mutually generic Cohen reals, and a bijection witnessing the cardinality of $A$. In theorem \ref{KS G,F Theorem}, one takes the generic filters $F$ and $G$ to have the same set $A$ of Cohen reals, but to disagree on the bijection from $A$ to a cardinal. Because the inner models $N^G(I,Q)\subseteq M[G]$ and $N^F(J,Q)\subseteq M[F]$ retain $A$ but lose its bijection to any cardinal, we arrive at the same model $N^G(I,2^{<\omega})=N^F(J,2^{<\omega})$ by removing the point of disagreement between $M[G]$ and $M[F]$.

With this sketch in mind, one would expect that the proofs of theorem \ref{KS G,F Theorem} and corollary \ref{Karagila-Schlicht Cohen Index Theorem} do not depend on the poset $2^{<\omega}$, and that the same results hold for $N(I,Q)$. Question 5.5 from \cite{karagila2020have} essentially asks whether this intuition holds. Here, the authors note that their ``proofs involved in a meaningful way the Cohen forcing [$\mathbb{P}(\omega,2^{<\omega})$] itself," but that ``on its face it seems that
the proof uses more of the fact that we take a finite-support product of infinitely
many copies of the same forcing, rather than the specific properties of the Cohen
forcing." They then ask (up to a change in notation) whether the results in \cite{karagila2020have} hold for $N(I,Q)$. We will answer this question positively, in the case of corollary \ref{Karagila-Schlicht Cohen Index Theorem}:

\begin{theorem}\label{Karagila-Schlicht Generalized Cohen Index Theorem}
Let $M\vDash\mathrm{ZFC}$ and let $I,J\in M$ be infinite. For every $\mathbb{P}(I,Q)$-generic $G$, there is a $\mathbb{P}(J,Q)$-generic $F$ so that
$$
N^G(I,Q)=N^F(J,Q).
$$
\end{theorem}
\noindent This formalizes the notion that the index set $I$ in $N(I,Q)$ can be taken to be arbitrarily large, without loss of generality. Two relevant corollaries are immediate.

\begin{corollary}[Halpern, L\'evy, \cite{halpern1971boolean}, Stefanovi\'c, \cite{stefanovic2023alternatives}]\label{N(I,Q) models BPI by applying KS}
    $N(I,Q)\vDash\mathrm{BPI}$.\footnote{This is the same statement as corollary \ref{N(I,Q) models BPI from FEP}, but we number this corollary distinctly since we use a different proof, and because we cite an established case of the proof of corollary \ref{N(I,Q) models BPI from FEP}.}
\end{corollary}
\begin{proof}
    Given $N^G(I,Q)$, take $\theta$ large enough for theorem \ref{N(theta,Q) Models BPI} to hold, and let $F$ be $\mathbb{P}(\theta,Q)$-generic so that $N^G(I,Q)=N^F(\theta,Q)$. By theorem \ref{N(theta,Q) Models BPI} and the proof of corollary \ref{N(I,Q) models BPI from FEP}, in the particular case of $N(\theta,Q)$, we have that $N^G(I,Q)=N^F(\theta,Q)\vDash\mathrm{BPI}$.
\end{proof}

\begin{corollary}[Halpern, L\'evy, \cite{halpern1971boolean}, Stefanovi\'c, \cite{stefanovic2023alternatives}]\label{N(I,Q) FEP (weak)}
    The symmetric system $\langle\mathbb{P}(I,Q),\mathcal{G},\mathcal{F}\rangle$ has the filter extension property for every $I$.
\end{corollary}
\begin{proof}
    Apply theorem \ref{BPI equivalent to FEP for sym ext} to corollary \ref{N(I,Q) models BPI by applying KS}
\end{proof}

From here, there are two clear loose ends in our pursuit of $\mathrm{BPI}$ in the generalized Cohen model: the need to formally prove theorem \ref{Karagila-Schlicht Generalized Cohen Index Theorem}, and the fact that corollary \ref{N(I,Q) FEP (weak)} is not the sharp version of the filter extension property given by theorem \ref{Halpern-Levy FEP Theorem}.

The proofs of theorem \ref{KS G,F Theorem} and corollary \ref{Karagila-Schlicht Cohen Index Theorem} in \cite{karagila2020have} are clear and insightful. We expect that they generalize in a straightforward way to prove theorem \ref{Karagila-Schlicht Generalized Cohen Index Theorem}, but for the purposes of this section, it will be advantageous to offer a new perspective on theorem \ref{Karagila-Schlicht Generalized Cohen Index Theorem}. The approach we take was inspired by a comment made in \cite{karagila2020have}, in the paragraph immediately following their remark 4.4. There, the authors point out that their result does not imply a bijection between $\omega$ and $\kappa$ (or, in our case, $I$ and $J$) in the ground model $M$, only that $\kappa$ is countable in is some ambient outer model of $M$. This suggests that one can study the results in \cite{karagila2020have} from the perspective of an outer model with a bijection between $I$ and $J$. We will use $M$ to denote our intended ground model, and $V$ to denote such an outer model of $M$ with a bijection $b:I\to J$, and in which $M$ is a class. In $V$, based on the bijection $b$, there are natural notions of isomorphism between $\mathbb{P}(I,Q)$ and $\mathbb{P}(J,Q)$, between $V^{\mathbb{P}(I,Q)}$ and $V^{\mathbb{P}(J,Q)}$, and thus between $\Vdash_{N(I,Q)}^V$ and $\Vdash_{N(J,Q)}^V$. We prove that, in $V$, the natural isomorphism between $V^{\mathbb{P}(I,Q)}$ and $V^{\mathbb{P}(J,Q)}$ restricts to an isomorphism between the relativized classes $(\mathrm{HS}_{\mathcal{F}_I})^M$ and $(\mathrm{HS}_{\mathcal{F}_J})^M$.\footnote{As noted in definition \ref{Generalized Cohen System}, we can include the subscripts $\langle\mathbb{P}(I,Q),\mathcal{G}_I,\mathcal{F}_I\rangle$ and $\langle\mathbb{P}(J,Q),\mathcal{G}_J,\mathcal{F}_J\rangle$ to emphasize the difference in index sets $I$ and $J$.} From this, we show that the isomorphism in $V$ between $\Vdash_{N(I,Q)}^V$ and $\Vdash_{N(J,Q)}^V$ restricts to an isomorphism between the relativized forcing relations $\Vdash_{N(I,Q)}^M$ and $\Vdash^M_{N(J,Q)}$. Even though this isomorphism is not in $M$ itself, the fact that the relations $\Vdash_{N(I,Q)}^M$ and $\Vdash^M_{N(J,Q)}$ are relativized to $M$ allows us to draw conclusions about the relationship between $\Vdash_{N(I,Q)}$ and $\Vdash_{N(J,Q)}$ in $M$. In particular, this allows us to transfer witnesses from theorem \ref{N(theta,Q) Models BPI} to theorem \ref{Halpern-Levy FEP Theorem}, and to prove theorem \ref{Karagila-Schlicht Generalized Cohen Index Theorem}.

To carry out this approach in the outer model $V$, it will be useful to have a thorough understanding of the relationship between $\mathrm{HS}_{\mathcal{F}_I}$ and $\mathrm{HS}_{\mathcal{F}_J}$ that can be described in the ground model $M$. In subsection 2.1, we introduce the maps in $M$
$$
\Sigma:\mathrm{HS}_{\mathcal{F}_I}\to\mathrm{HS}_{\mathcal{F}_J},\ \Gamma:\mathrm{HS}_{\mathcal{F}_J}\to\mathrm{HS}_{\mathcal{F}_I}.
$$
These maps capture, in a strong sense, the ability of $M$ to compare the classes $\mathrm{HS}_{\mathcal{F}_I}$ and $\mathrm{HS}_{\mathcal{F}_J}$, and hence to compare the symmetric systems $\langle\mathbb{P}(I,Q),\mathcal{G}_I,\mathcal{F}_I\rangle$ and $\langle\mathbb{P}(J,Q),\mathcal{G}_J,\mathcal{F}_J\rangle$.

In subsection 2.2, we use the maps $\Sigma$ and $\Gamma$ to analyze the isomorphisms that extend from a bijection $b:I\to J$ in an outer model $V$ of $M$, in which $M$ is a class. We then carry out the details of the approach sketched above, which completes our full proof that $N(I,Q)\vDash\mathrm{BPI}$, and of theorem \ref{Halpern-Levy FEP Theorem}. This approach incidentally allows us to answer a related question: to what extent can $M$ internally express a notion of isomorphism between $N(I,Q)$ and $N(J,Q)$? We prove that the isomorphism between $\Vdash^M_{N(I,Q)}$ and $\Vdash^M_{N(J,Q)}$ in $V$ is closely related to the map $\Sigma\in M$; from this, it follows that $\Sigma$ is (the nontrivial part of) an elementary embedding from $\Vdash_{N(I,Q)}$ to $\Vdash_{N(J,Q)}$ in $M$.

In subsection 2.3, with our full proof that $N(I,Q)\vDash\mathrm{BPI}$ in hand, we fill a gap to complete a separate new proof that $N(I,Q)\vDash\mathrm{BPI}$, based on the characterization given by Blass in theorem \ref{Blass SVC BPI Theorem}. This alternate approach is outlined as question 8 on Karagila's website, and is answered/completed by a special case of our main theorems from this section. We conclude by proving a related maximal generalization of theorem \ref{Halpern-Levy FEP Theorem}.

\subsection{The Maps $\Sigma$ and $\Gamma$ in $M$}
As stated above, we introduce the maps in $M$
$$
\Sigma:\mathrm{HS}_{\mathcal{F}_I}\to\mathrm{HS}_{\mathcal{F}_J},\ \Gamma:\mathrm{HS}_{\mathcal{F}_J}\to\mathrm{HS}_{\mathcal{F}_I}
$$
to capture the ability of $M$ to compare $\mathrm{HS}_{\mathcal{F}_I}$ with $\mathrm{HS}_{\mathcal{F}_J}$, and hence to compare $\langle\mathbb{P}(I,Q),\mathcal{G}_I,\mathcal{F}_I\rangle$ with $\langle\mathbb{P}(J,Q),\mathcal{G}_J,\mathcal{F}_J\rangle$. In this manner, we can begin to address the following more general questions, in the case of the symmetric systems $\langle\mathbb{P}(I,Q),\mathcal{G}_I,\mathcal{F}_I\rangle$ and $\langle\mathbb{P}(J,Q),\mathcal{G}_J,\mathcal{F}_J\rangle$.

\begin{question}\label{Iso btw sym systems question}
    What is a good notion of isomorphism between symmetric systems?
\end{question}

\begin{question}\label{Iso btw sym systems in M question}
    To what extent can question \ref{Iso btw sym systems question} be answered using maps in the given model $M\vDash\mathrm{ZFC}$, as opposed to maps in some outer model $V$ of $M$?
\end{question}

\noindent By corollary \ref{Karagila-Schlicht Cohen Index Theorem}, question \ref{Iso btw sym systems question} and \ref{Iso btw sym systems in M question} are nontrivial: the symmetric systems $\langle\mathbb{P}(I,2^{<\omega}),\mathcal{G}_I,\mathcal{F}_I\rangle$ and $\langle\mathbb{P}(J,2^{<\omega}),\mathcal{G}_J,\mathcal{F}_J\rangle$ present the same symmetric extension $N^G(I,Q)=N^F(J,Q)$, despite there being no isomorphism in $M$ between $\mathbb{P}(I,Q)$ and $\mathbb{P}(J,Q)$ when $|I|\neq |J|$ in $M$.

We will assume without loss of generality that $I\subseteq J$ and $|I|<|J|$ in $M$. Given this assumption, we have that $\mathbb{P}(I,Q)\subseteq\mathbb{P}(J,Q)$. Our definition of the map $\Sigma$ aims to provide a canonical sense in which a name $\dot{x}\in\mathrm{HS}_{\mathcal{F}_I}$ can be modified to produce a corresponding name in $\mathrm{HS}_{\mathcal{F}_J}$.

\begin{definition}\label{Sigma Definition}
    Let $I,J$ be infinite, with $I\subseteq J$ and $|I|<|J|$. Define the ``spreading" function
    $$
    \Sigma:\mathrm{HS}_{\mathcal{F}_I}\to\mathrm{HS}_{\mathcal{F}_J}
    $$
    by recursively setting
    $$
    \Sigma(\dot{x})=\bigcup_{\langle p,\dot{z}\rangle\in\dot{x}}\mathrm{Orb}(\langle p,\Sigma(\dot{z})\rangle,\mathrm{fix}_{\mathcal{G}_J}(\mathrm{supp}(\dot{x}))).
    $$
\end{definition}

We can think of $\Sigma$ as hereditarily ``spreading" the name $\dot{x}\in\mathrm{HS}_{\mathcal{F}_I}$ across the index set $J$, while fixing the part of $\dot{x}$ above its support. Ultimately, in theorem \ref{Sigma Elementary Embedding Theorem}, we will address the relevant case of question \ref{Iso btw sym systems in M question} by showing that $\Sigma$ can be used to provide an elementary embedding in $M$ between the relativized forcing relations $\Vdash_{N(I,Q)}$ and $\Vdash_{N(J,Q)}$. Lemma \ref{Sigma into HS_FJ^I} justifies that $\mathrm{rng}(\Sigma)\subseteq\mathrm{HS}_{\mathcal{F}_J}$.

\begin{lemma}\label{Sigma into HS_FJ^I}
    If $\dot{x}\in\mathrm{HS}_{\mathcal{F}_I}$, then $\Sigma(\dot{x})\in\mathrm{HS}_{\mathcal{F}_J}$ and $\mathrm{supp}(\Sigma(\dot{x}))\subseteq \mathrm{supp}(\dot{x})\subseteq I$.\footnote{In corollary \ref{Supports under Sigma, Gamma corollary}, we strengthen this to show that $\mathrm{supp}(\dot{x})=\mathrm{supp}(\Sigma(\dot{x}))$.}
\end{lemma}
\begin{proof}
    Assume inductively that for every $\dot{z}\in\mathrm{rng}(\dot{x})$, $\Sigma(\dot{z})\in\mathrm{HS}_{\mathcal{F}_J}$. From the definition of $\Sigma$, it is apparent that $\mathrm{fix}_{\mathcal{G}_J}(\mathrm{supp}(\dot{x}))\preceq\mathrm{sym}_{\mathcal{G}_J}(\Sigma(\dot{x}))$. Because $I\subseteq J$, note that $\mathrm{supp}(\dot{x})\in [I]^{<\omega}\subseteq [J]^{<\omega}$, and so $\mathrm{fix}_{\mathcal{G}_J}(\mathrm{supp}(\dot{x}))\in\mathcal{F}_J$. Thus, $\Sigma(\dot{x})$ is $\mathcal{F}_J$-symmetric, and by the inductive assumption, $\Sigma(\dot{x})\in\mathrm{HS}_{\mathcal{F}_J}$. By lemma \ref{Minimal Supports Lemma}, $\mathrm{supp}(\Sigma(\dot{x}))$ is minimal, and so $\mathrm{supp}(\Sigma(\dot{x}))\subseteq\mathrm{supp}(\dot{x})\subseteq I$.
\end{proof}

An important ingredient in the proofs of our main theorems is that every name in $\mathrm{HS}_{\mathcal{F}_J}$ can be described, in $M$, in terms of $\mathrm{HS}_{\mathcal{F}_I}$, $\Sigma$, and $\mathcal{G}_J$ (see corollary \ref{HS_FJ in terms of HS_FI, Sigma, G_J Corollary}). To show this, it will be helpful to define a ``gathering" function $\Gamma:\mathrm{HS}_{\mathcal{F}_J}\to\mathrm{HS}_{\mathcal{F}_I}$, and prove that the restriction of $\Gamma$ to the class $\{\dot{w}\in\mathrm{HS}_{\mathcal{F}_J}\ |\ \mathrm{supp}(\dot{w})\subseteq I\}$ is the inverse of $\Sigma$ (see theorem \ref{Sigma, Gamma Inverses on support I}).

\begin{definition}\label{Gamma definition}
    Let $I,J$ be infinite, with $I\subseteq J$ and $|I|<|J|$. Recursively define the ``gathering" function
    $$
    \Gamma:\mathrm{HS}_{\mathcal{F}_J}\to\mathrm{HS}_{\mathcal{F}_I}
    $$
    by setting
    $$
    \Gamma(\dot{w})=\{\langle p,\Gamma(\dot{z})\rangle\ |\ \langle p,\dot{z}\rangle\in\dot{w}\wedge\mathrm{supp}(p),\mathrm{supp}(\dot{z})\subseteq I\}.
    $$
\end{definition}

We can think of $\Gamma$ as hereditarily ``gathering" the subset of $\dot{w}$ that can be stated with respect to the index set $I\subseteq J$, and thus can be expressed as a $\mathbb{P}(I,Q)$-name. 

\begin{lemma}\label{Gamma commutes with G_J}
    For $\dot{w}\in\mathrm{HS}_{\mathcal{F}_J}$ and $\pi\in\mathcal{G}_I$, $\pi\Gamma(\dot{w})=\Gamma(\pi\dot{w})$.
\end{lemma}
\begin{proof}
    Let $\dot{w}\in\mathrm{HS}_{\mathcal{F}_J}$ and inductively assume that for every $\dot{z}\in\mathrm{rng}(\dot{w})$ and every $\pi\in\mathcal{G}_I$, $\pi\Gamma(\dot{z})=\Gamma(\pi\dot{z})$. We will first show that $\pi\Gamma(\dot{w})\subseteq\Gamma(\pi\dot{w})$. Suppose that $\langle\pi p,\pi\Gamma(\dot{z})\rangle\in\pi\Gamma(\dot{w})$, with $\langle p,\Gamma(\dot{z})\rangle\in\Gamma(\dot{w})$, and $\langle p,\dot{z}\rangle\in\dot{w}$ with $\mathrm{supp}(p),\mathrm{supp}(\dot{w})\subseteq I$. Then $\langle \pi p,\pi\dot{w}\rangle\in\pi\dot{w}$, and because $\pi\in\mathcal{G}_I$, $\mathrm{supp}(\pi p),\mathrm{supp}(\pi\dot{w})\subseteq I$. From this, we see that $\langle \pi p,\Gamma(\pi\dot{z})\rangle\in\Gamma(\pi\dot{w})$, and by the inductive hypothesis $\langle \pi p,\pi\Gamma(\dot{z})\rangle\in\Gamma(\pi\dot{w})$.

    To prove the containment $\Gamma(\pi\dot{w})\subseteq\pi\Gamma(\dot{w})$, let $\langle \pi p,\Gamma(\pi\dot{z})\rangle\in\Gamma(\pi\dot{w})$, with $\langle \pi p,\pi\dot{z}\rangle\in\pi\dot{w}$ and $\mathrm{supp}(\pi p),\mathrm{supp}(\pi\dot{z})\subseteq I$. Then $\langle p,\dot{z}\rangle\in \dot{w}$, and because $\pi,\pi^{-1}\in\mathcal{G}_I$, $\mathrm{supp}( p),\mathrm{supp}(\dot{z})\subseteq I$. By the definition of $\Gamma(\dot{w})$, $\langle p,\Gamma(\dot{z})\rangle\in\Gamma(\dot{w})$, and so $\langle \pi p,\pi\Gamma(\dot{z})\rangle\in\pi\Gamma(\dot{w})$. By the inductive hypothesis, $\langle\pi p,\Gamma(\pi\dot{z})\rangle\in\pi\Gamma(\dot{w})$, which completes the proof.
\end{proof}

\begin{lemma}\label{Gamma into HS_FI}
    If $\dot{w}\in\mathrm{HS}_{\mathcal{F}_J}$ and $\mathrm{supp}(\dot{w})\subseteq I$, then $\Gamma(\dot{w})\in\mathrm{HS}_{\mathcal{F}_I}$ and $\mathrm{supp}(\Gamma(\dot{w}))\subseteq\mathrm{supp}(\dot{w})$.\footnote{In corollary \ref{Supports under Sigma, Gamma corollary}, we strengthen this to show that $\mathrm{supp}(\dot{w})=\mathrm{supp}(\Gamma(\dot{w}))$.}
\end{lemma}
\begin{proof}
    Suppose inductively that for every $\dot{z}\in\mathrm{rng}(\dot{w})$ with $\mathrm{supp}(\dot{z})\subseteq I$, $\Gamma(\dot{z})\in\mathrm{HS}_{\mathcal{F}_I}$. From this, we can see that $\Gamma(\dot{w})\in M^{\mathbb{P}(I,Q)}$, and that $\Gamma(\dot{z})\in\mathrm{HS}_{\mathcal{F}_I}$ if $\mathrm{sym}_{\mathcal{G}_I}(\Gamma(\dot{w}))\in\mathcal{F}_I$. Let $\mathrm{supp}(\dot{w})=W\in[I]^{<\omega}$; we will show that $\mathrm{fix}_{\mathcal{G}_I}(W)\preceq\mathrm{sym}_{\mathcal{G}_I}(\Gamma(\dot{w}))$, and thus that $\mathrm{sym}_{\mathcal{G}_I}(\Gamma(\dot{w}))\in\mathcal{F}_I$.
    
    We claim that for every $\mathrm{Orb}(\langle p,\dot{z}\rangle,\mathrm{fix}_{\mathcal{G}_J}(W))\subseteq\dot{w}$, it is possible to assume without loss of generality that $\mathrm{supp}(p),\mathrm{supp}(\dot{z})\subseteq I$. To see why, let $\mathrm{supp}(p)\cup\mathrm{supp}(\dot{z})=S$; since $W\subseteq I$, there is a $\pi\in \mathrm{fix}_{\mathcal{G}_J}(W)$ so that $\pi S\subseteq I$. Then $\mathrm{supp}(\pi p),\mathrm{supp}(\pi\dot{z})\subseteq I$, and
    $$
    \mathrm{Orb}(\langle p,\dot{z}\rangle,\mathrm{fix}_{\mathcal{G}_J}(W))=\mathrm{Orb}(\langle \pi p,\pi \dot{z}\rangle,\mathrm{fix}_{\mathcal{G}_J}(W)).
    $$
     As a result, we can write
    \begin{equation}\label{w in terms of I supports}
    \dot{w}=\bigcup_{\langle p,\dot{z}\rangle\in \dot{w}\wedge\mathrm{supp}(p),\mathrm{supp}(\dot{z})\subseteq I}\mathrm{Orb}(\langle q,\dot{z}\rangle,\mathrm{fix}_{\mathcal{G}_J}(W)).
    \end{equation}
    Given that $\mathrm{supp}(p),\mathrm{supp}(\dot{z})\subseteq I$, note that
    $$
    \{\langle q,\dot{y}\rangle\in\mathrm{Orb}(\langle p,\dot{z}\rangle,\mathrm{fix}_{\mathcal{G}_J}(W))\ |\ \mathrm{supp}( q),\mathrm{supp}(\dot{y})\subseteq I\}=\mathrm{Orb}(\langle p,\dot{z}\rangle,\mathrm{fix}_{\mathcal{G}_I}(W)).
    $$
    By this, and from the definition of $\Gamma$, we have that
    $$
    \Gamma(\dot{w})=\bigcup_{\langle p,\dot{z}\rangle\in\dot{w}\wedge\mathrm{supp}(p),\mathrm{supp}(\dot{z})\subseteq I}\{\langle\pi p,\Gamma(\pi\dot{z})\rangle\ |\ \pi\in\mathrm{fix}_{\mathcal{G}_I}(W)\}.
    $$
    By lemma \ref{Gamma commutes with G_J}, since $\pi\in\mathrm{fix}_{\mathcal{G}_I}(W)\preceq\mathcal{G}_I$, we have $\Gamma(\pi\dot{z})=\pi\Gamma(\dot{z})$, so this can be rewritten as
    \begin{equation}\label{Gamma(w) form}
    \Gamma(\dot{w})=\bigcup_{\langle p,\dot{z}\rangle\in\dot{w}\wedge\mathrm{supp}(p),\mathrm{supp}(\dot{z})\subseteq I}\mathrm{Orb}(\langle p,\Gamma(\dot{z})\rangle,\mathrm{fix}_{\mathcal{G}_I}(W)).
    \end{equation}
    From this, it is clear that $\mathrm{fix}_{\mathcal{G}_I}(W)\preceq\mathrm{sym}_{\mathcal{G}_I}(\Gamma(\dot{w}))$, which completes the proof.
\end{proof}

\begin{remark}
    Generally, it is true that if $\dot{w}\in\mathrm{HS}_{\mathcal{F}_J}$, then $\Gamma(\dot{w})\in\mathrm{HS}_{\mathcal{F}_I}$, as the notation from definition \ref{Gamma definition} suggests. Because we will not use this in any proof, we leave this detail to the reader.
\end{remark}

\begin{theorem}\label{Sigma, Gamma Inverses on support I}
    If $\dot{x}\in \mathrm{HS}_{\mathcal{F}_I}$, then $\Gamma(\Sigma(\dot{x}))=\dot{x}$; if   $\dot{w}\in\mathrm{HS}_{\mathcal{F}_J}$ and $\mathrm{supp}(\dot{w})\subseteq I$, then $\Sigma(\Gamma(\dot{w}))=\dot{w}$.
\end{theorem}
\begin{proof}
To prove the first part of the theorem, let $\dot{x}\in\mathrm{HS}_{\mathcal{F}_I}$ with $\mathrm{supp}(\dot{x})=X$, and suppose inductively that for every $\dot{z}\in\mathrm{rng}(\dot{x})$, $\Gamma(\Sigma(\dot{z}))=\dot{z}$. From the definition of $\Sigma$ and $\Gamma$, we can write
$$
\Sigma(\dot{x})=\bigcup_{\langle p,\dot{z}\rangle\in\dot{x}}\mathrm{Orb}(\langle p,\Sigma(\dot{z})\rangle,\mathrm{fix}_{\mathcal{G}_J}(X)),
$$
and
$$
\Gamma(\Sigma(\dot{x}))=\{\langle \pi p,\Gamma(\pi\Sigma(\dot{z}))\rangle\ |\ \langle p,\dot{z}\rangle\in\dot{x}\wedge \pi\in\mathrm{fix}_{\mathcal{G}_J}(X)\wedge\mathrm{supp}(\pi p),\mathrm{supp}(\pi\Sigma(\dot{z}))\subseteq I\}.
$$
To break this expression down further, we claim, for$\langle p,\dot{z}\rangle\in\dot{x}$ and $\pi\in\mathcal{G}_J$, that $\mathrm{supp}(\pi p),\mathrm{supp}(\pi\Sigma(\dot{z}))\subseteq I$ if and only if there is a $\tau\in\mathrm{fix}_{\mathcal{G}_I}(X)$ so that $\pi p=\tau p$ and $\pi\Sigma(\dot{z})=\tau\Sigma(\dot{z})$. To prove this, let $\langle \pi p,\pi\Sigma(\dot{z})\rangle\in\Sigma(\dot{x})$, with $\langle p,\dot{z}\rangle\in\dot{x}$ and $\pi\in\mathrm{fix}_{\mathcal{G}_J}(X)$. Then $\mathrm{supp}(p)\subseteq I$, and by lemma \ref{Sigma into HS_FJ^I}, $\mathrm{supp}(\Sigma(\dot{z}))\subseteq I$. The claim then holds by lemma \ref{image support lemma}. As a result, we can write
$$
\Gamma(\Sigma(\dot{x}))=\{\langle \tau p,\Gamma(\tau\Sigma(\dot{z}))\rangle\ |\ \langle p,\dot{z}\rangle\in\dot{x}\wedge\tau\in\mathrm{fix}_{\mathcal{G}_I}(X)\}.
$$
By lemma \ref{Gamma commutes with G_J}, since $\tau\in\mathcal{G}_I$, we can replace this with 
$$
\Gamma(\Sigma(\dot{x}))=\{\langle \tau p,\tau\Gamma(\Sigma(\dot{z}))\rangle\ |\ \langle p,\dot{z}\rangle\in\dot{x}\wedge\tau\in\mathrm{fix}_{\mathcal{G}_I}(X)\}.
$$
By the inductive hypothesis, this can again be replaced with
$$
    \Gamma(\Sigma(\dot{x}))=\{\langle \tau p,\tau\dot{z}\rangle\ |\ \langle p,\dot{z}\rangle\in\dot{x}\wedge\tau\in\mathrm{fix}_{\mathcal{G}_I}(X)\}=\bigcup_{\langle p,\dot{z}\rangle\in\dot{x}}\mathrm{Orb}(\langle p,\dot{z}\rangle,\mathrm{fix}_{\mathcal{G}_I}(X))=\dot{x},
$$
which proves the first part of the theorem.

To prove the second part of the theorem, let $\mathrm{supp}(\dot{w})=W$, and suppose inductively for every $\dot{z}\in\mathrm{rng}(\dot{w})$ with $\mathrm{supp}(\dot{z})\subseteq I$, that $\Sigma(\Gamma(\dot{z}))=\dot{z}$. Because $W\subseteq I$, we meet the hypothesis of lemma \ref{Gamma into HS_FI}; then, as in lines (\ref{w in terms of I supports}) and (\ref{Gamma(w) form}), we can express $\dot{w}$ and $\Gamma(\dot{w})$ as
$$
    \dot{w}=\bigcup_{\langle p,\dot{z}\rangle\in \dot{w}\wedge\mathrm{supp}(p),\mathrm{supp}(\dot{z})\subseteq I}\mathrm{Orb}(\langle q,\dot{z}\rangle,\mathrm{fix}_{\mathcal{G}_J}(W)).
$$
and
$$
    \Gamma(\dot{w})=\bigcup_{\langle p,\dot{z}\rangle\in\dot{w}\wedge\mathrm{supp}(p),\mathrm{supp}(\dot{z})\subseteq I}\mathrm{Orb}(\langle p,\Gamma(\dot{z})\rangle,\mathrm{fix}_{\mathcal{G}_I}(W)).
$$
From this, and the definition of $\Sigma$, we can see that
$$
\Sigma(\Gamma(\dot{w}))=\bigcup_{\langle p,\dot{z}\rangle\in\dot{w}\wedge\mathrm{supp}(p),\mathrm{supp}(\dot{z})\subseteq I}\mathrm{Orb}(\langle p,\Sigma(\Gamma(\dot{z}))\rangle,\mathrm{fix}_{\mathcal{G}_J}(W)).
$$
To complete the proof, apply the inductive hypothesis to the names $\dot{z}$ in the subscript of the union, for which $\mathrm{supp}(\dot{z})\subseteq I$, to get
$$
\Sigma(\Gamma(\dot{w}))=\bigcup_{\langle p,\dot{z}\rangle\in\dot{w}\wedge\mathrm{supp}(p),\mathrm{supp}(\dot{z})\subseteq I}\mathrm{Orb}(\langle p,\dot{z}\rangle,\mathrm{fix}_{\mathcal{G}_J}(W))=\dot{w}.
$$

\end{proof}

\begin{corollary}\label{HS_FJ in terms of HS_FI, Sigma, G_J Corollary}
    If $\dot{w}\in\mathrm{HS}_{\mathcal{F}_J}$, then $\pi^{-1}\Sigma(\Gamma(\pi\dot{w}))=\dot{w}$ for every $\pi\in\mathcal{G}_J$ so that $\pi\mathrm{supp}(\dot{w})\subseteq I$.
\end{corollary}
\begin{proof}
    Since $\mathrm{supp}(\pi\dot{w})=\pi\mathrm{supp}(\dot{w})\subseteq I$, by theorem \ref{Sigma, Gamma Inverses on support I}, $\pi^{-1}(\Sigma(\Gamma(\pi\dot{w})))=\pi^{-1}(\pi\dot{w})=\dot{w}.$
\end{proof}

Corollary \ref{HS_FJ in terms of HS_FI, Sigma, G_J Corollary} establishes our earlier claim that $M$ can describe the class $\mathrm{HS}_{\mathcal{F}_J}$ in terms of $\mathrm{HS}_{\mathcal{F}_I},\Sigma,$ and $\mathcal{G}_J$. The following corollary will be useful in the subsection 5.2.

\begin{corollary}\label{Sigma commutes with G_I corollary}
    If $\dot{x}\in\mathrm{HS}_{\mathcal{F}_I}$ and $\pi\in\mathcal{G}_I$, then $\Sigma(\pi\dot{x})=\pi\Sigma(\dot{x})$.
\end{corollary}
\begin{proof}
    Let $\dot{x}\in\mathrm{HS}_{\mathcal{F}_I}$, and let $\Sigma(\dot{x})=\dot{w}$. To prove the corollary, we will justify the line
    $$
    \Sigma(\pi\dot{x})=\Sigma(\pi\Gamma(\dot{w}))=\Sigma(\Gamma(\pi\dot{w}))=\pi\dot{w}=\pi\Sigma(\dot{x}).
    $$

    By theorem \ref{Sigma, Gamma Inverses on support I}, $\Gamma(\Sigma(\dot{x}))=\dot{x}$, so $\dot{x}$ can be replaced with $\Gamma(\dot{w})$ to show the first equality. Since we assume that $I\subseteq J$, we have that $\pi\in\mathcal{G}_I\preceq\mathcal{G}_J$; then, by lemma \ref{Gamma commutes with G_J}, we have that $\Gamma(\pi\dot{w})=\pi\Gamma(\dot{w})$ to justify the second equality. We then claim that $\mathrm{supp}(\pi\dot{w})\subseteq I$. We have that $\mathrm{supp}(\pi\dot{w})=\mathrm{supp}(\pi\Sigma(\dot{x}))=\pi\mathrm{supp}(\Sigma(\dot{x}))$; by lemma \ref{Sigma into HS_FJ^I}, we have that $\pi\mathrm{supp}(\Sigma(\dot{x}))\subseteq\pi\mathrm{supp}(\dot{x})$; because $\pi\in\mathcal{G}_I$ and $\mathrm{supp}(\dot{x})\subseteq I$, $\pi\mathrm{supp}(\dot{x})\subseteq I$ and the claim holds. Thus, because $\mathrm{supp}(\pi\dot{w})\subseteq I$ and by theorem \ref{Sigma, Gamma Inverses on support I}, we have that $\Sigma(\Gamma(\pi\dot{w}))=\pi\dot{w}$. Lastly, we had set $\dot{w}=\Sigma(\dot{x})$ by definition.
\end{proof}

To conclude the subsection, although our main theorems can be proven without using this next corollary, it is quick to prove and helps to fill out a mental image of the maps $\Sigma$ and $\Gamma$.

\begin{corollary}\label{Supports under Sigma, Gamma corollary}
    If $\dot{x}\in\mathrm{HS}_{\mathcal{F}_I}$, then $\mathrm{supp}(\dot{x})=\mathrm{supp}(\Sigma(\dot{x}))$; if $\dot{w}\in\mathrm{HS}_{\mathcal{F}_J}$ and $\mathrm{supp}(\dot{w})\subseteq I$, then $\mathrm{supp}(\dot{w})=\mathrm{supp}(\Gamma(\dot{w}))$.
\end{corollary}
\begin{proof}
    Let $\dot{x}\in\mathrm{HS}_{\mathcal{F}_I}$. By lemma \ref{Sigma into HS_FJ^I}, $\Sigma(\dot{x})\in\mathrm{HS}_{\mathcal{F}_J}$ and $\mathrm{supp}(\Sigma(\dot{x}))\subseteq\mathrm{supp}(\dot{x})$. By lemma \ref{Gamma into HS_FI}, $\mathrm{supp}(\Gamma(\Sigma(\dot{x})))\subseteq\mathrm{supp}(\Sigma(\dot{x}))$, and by theorem \ref{Sigma, Gamma Inverses on support I}, $\mathrm{supp}(\Gamma(\Sigma(\dot{x})))=\mathrm{supp}(\dot{x})$. Thus $\mathrm{supp}(\dot{x})\subseteq\mathrm{supp}(\Sigma(\dot{x}))$, and so $\mathrm{supp}(\dot{x})=\mathrm{supp}(\Sigma(\dot{x}))$.

    Let $\dot{w}\in\mathrm{HS}_{\mathcal{F}_J}$ with $\mathrm{supp}(\dot{w})\subseteq I$. By lemma \ref{Gamma into HS_FI}, $\Gamma(\dot{w})\in\mathrm{HS}_{\mathcal{F}_I}$ and $\mathrm{supp}(\Gamma(\dot{w}))\subseteq\mathrm{supp}(\dot{w})$. By lemma \ref{Sigma into HS_FJ^I}, $\mathrm{supp}(\Sigma(\Gamma(\dot{w})))\subseteq\mathrm{supp}(\Gamma(\dot{w}))$, and by theorem \ref{Sigma, Gamma Inverses on support I}, $\mathrm{supp}(\Sigma(\Gamma(\dot{w})))=\mathrm{supp}(\dot{w})$. Thus $\mathrm{supp}(\dot{w})\subseteq\mathrm{supp}(\Gamma(\dot{w}))$, and so $\mathrm{supp}(\dot{w})=\mathrm{supp}(\Gamma(\dot{w}))$.
\end{proof}

\subsection{Isomorphisms in $V\supseteq M$ and Elementary Embeddings in $M$}
We can now fix a bijection $b:I\to J$ in some outer model $V$ of $M$, in which $M$ is a class, and consider the natural isomorphisms in $V$ between $V^{\mathbb{P}(I,Q)}$ and $V^{\mathbb{P}(J,Q)}$, and between $\Vdash_{N(I,Q)}^V$ and $\Vdash_{N(J,Q)}^V$. We will show that, in $V$, these maps restrict to isomorphisms between the relativizations $(\mathrm{HS}_{\mathcal{F}_I})^M$ and $(\mathrm{HS}_{\mathcal{F}_J})^M$, and between the relativizations $\Vdash^M_{N(I,Q)}$ and $\Vdash^M_{N(J,Q)}$. We can first formally state the natural isomorphisms that $b$ extends to in $V$.

\begin{definition}\label{b and varphi_b definition}
    Suppose that $b:I\to J$ is a bijection in $V\vDash\mathrm{ZFC}$. By a slight misuse of notation, define the map
    $$
    b:\mathbb{P}(I,Q)\to\mathbb{P}(J,Q)
    $$
    by sending $\mathrm{supp}(b(p))=b(\mathrm{supp}(p))$ and for every $i\in\mathrm{supp}(b(p))$, setting
    $$
    b(p)(b(i))=p(i).
    $$
    From this, recursively define the map
    $$
    \varphi_b:V^{\mathbb{P}(I,Q)}\to V^{\mathbb{P}(J,Q)}
    $$
    by
    $$
    \varphi_b(\dot{x})=\{\langle b(p),\varphi_b(\dot{z})\rangle\ |\ \langle p,\dot{z}\rangle\in\dot{x}\}.
    $$
\end{definition}

\begin{observation}\label{observation}
Notice that for $p\in\mathbb{P}(I,Q)$, the definition of $b(p)$ can equivalently be expressed by setting $b(p)=\pi p$, for any $\pi\in\mathcal{G}_J$ so that $\pi(i)=b(i)$ for every $i\in\mathrm{supp}(p)$.
\end{observation}

It will be convenient to state the following lemma, which follows from this observation (or immediatly from definition \ref{b and varphi_b definition}).

\begin{lemma}\label{pi b lemma}
    Suppose that $b:I\to J$ is a bijection in $V\vDash\mathrm{ZFC}$, and that $\pi\in\mathcal{G}_J$. Then $b':I\to J$, defined by
    $$
    b'(i)=\pi b(i),
    $$
    is a bijection, and the consequent maps
    $$
    b':\mathbb{P}(I,Q)\to\mathbb{P}(J,Q),\ 
    \varphi_{b'}:V^{\mathbb{P}(I,Q)}\to V^{\mathbb{P}(J,Q)},
    $$
    from definition \ref{b and varphi_b definition}, are given by
    $$
    b'(p)=\pi b(p),\ 
    \varphi_{b'}(\dot{x})=\pi \varphi_b(\dot{x}).
    $$
\end{lemma}

The following lemma is also standard to verify.

\begin{lemma}\label{Checking that b, varphi_b are isos}
    Given a bijection $b:I\to J$, the map $b:\mathbb{P}(I,Q)\to \mathbb{P}(J,Q)$ in definition \ref{b and varphi_b definition} is a $\leq$-preserving isomorphism, and the map $\varphi_b:V^{\mathbb{P}(I,Q)}\to V^{\mathbb{P}(J,Q)}$ is an isomorphism in the sense that it is a bijection that preserves the relation
    $$
    \langle p,\dot{z}\rangle\in\dot{x}\iff\langle b(p),\varphi_b(\dot{z})\rangle\in\varphi_b(\dot{x}).
    $$
    Moreover, $\langle b,\varphi_b\rangle$ is an isomorphism between $\Vdash^V_{N(I,Q)}$ and $\Vdash^V_{N(J,Q)}$ in the sense that for every formula $\phi$, we have that for every $\dot{x}_0,...,\dot{x}_n\in(\mathrm{HS}_{\mathcal{F}_I})^V$
    $$
    p\Vdash^V_{N(I,Q)}\phi(\dot{x}_0,...,\dot{x}_n)\iff b(p)\Vdash^V_{N(J,Q)}\phi(\varphi_b(\dot{x}_0),...,\varphi_b(\dot{x}_n)),
    $$
    and for every $\dot{w}_0,...,\dot{w}_n\in(\mathrm{HS}_{\mathcal{F}_J})^V$,
    $$
    b^{-1}(p)\Vdash^V_{N(I,Q)}\phi(\varphi_b^{-1}(\dot{w}_0),...,\varphi_b^{-1}(\dot{w}_n))\iff p\Vdash^V_{N(J,Q)}\phi(\dot{w}_0,...,\dot{w}_n).
    $$
\end{lemma}

We can now consider how the map $\varphi_b$ behaves when restricted to $(\mathrm{HS}_{\mathcal{F}_I})^M$. Put informally, if $\dot{x}\in\mathrm{HS}_{\mathcal{F}_I},$ with $X=\mathrm{supp}(\dot{x})$, then $\varphi_b$ alters $\dot{x}$ in two ways:
    \begin{enumerate}[label=(\roman*)]
        \item it takes the part of $\dot{x}$ defined outside of its support $X$ and spreads it across $J\setminus X$;
        \item it copies the resulting name from (i), which still has the support $X$, and redefines it over the new support $b(X)=\{b(i)\ |\ i\in X\}$.
    \end{enumerate}
Point $\mathrm{(i)}$ can be expressed as taking $\dot{x}$ to $\Sigma(\dot{x})$; point $\mathrm{(ii)}$ can be expressed as taking $\Sigma(\dot{x})$ to $\pi\Sigma(\dot{x})$, for some $\pi$ so that $\pi(i)=b(i)$ for every $i\in X$. Theorem \ref{varphi_b in terms of Sigma, G_J} verifies that this intuition is correct, and establishes the central relationship between $\Sigma\in M$ and $\varphi_b\in V\setminus M$.

\begin{theorem}\label{varphi_b in terms of Sigma, G_J}
    Let $M\vDash\mathrm{ZFC}$, $I\subseteq J$ and $|I|<|J|$ in $M$; let $V\vDash\mathrm{ZFC}$ be an outer model of $M$, with $M$ a class in $V$, and let $b\in V$ be a bijection from $I$ to $J$. If $\dot{x}\in\mathrm{HS}_{\mathcal{F}_I}$, $\pi\in\mathcal{G}_J$, and $\pi(i)=b(i)$ for every $i\in\mathrm{supp}(\dot{x})$, then
    $$
    \varphi_b(\dot{x})=\pi\Sigma(\dot{x}).
    $$
\end{theorem}
\begin{proof}
    Let $\dot{x}\in\mathrm{HS}_{\mathcal{F}_I}$, $\mathrm{supp}(\dot{x})=X$, and $\pi\in\mathcal{G}_J$ for which $\pi(i)=b(i)$ for every $i\in X$. Suppose inductively that for every bijection $b:I\to J$, every $\dot{z}\in\mathrm{rng}(\dot{x})$, and every $\tau\in\mathcal{G}_J$, that if $\tau(i)=b(i)$ for every $i\in\mathrm{supp}(\dot{z})$, then $\varphi_b(\dot{z})=\tau\Sigma(\dot{z})$.
    
    Let $b'=\pi^{-1}b$, so that $b'(i)=i$ for every $i\in X$. We will show that $\varphi_{b'}(\dot{x})=\Sigma(\dot{x})$, in which case $\varphi_b(\dot{x})=\pi\Sigma(\dot{x})$ by lemma \ref{pi b lemma}. Recall that
    $$
    \varphi_{b'}(\dot{x})=\{\langle b'(p),\varphi_{b'}(\dot{z})\rangle\ |\ \langle p,\dot{z}\rangle\in\dot{x}\}
    $$
    and
    $$
    \Sigma(\dot{x})=\bigcup_{\langle p,\dot{z}\rangle\in \dot{x}}\mathrm{Orb}(\langle p,\dot{z}\rangle,\mathrm{fix}_{\mathcal{G}_J}(X)).
    $$
    
    To show that $\varphi_{b'}(\dot{x})\subseteq\Sigma(\dot{x})$, let $\langle b'(p),\varphi_{b'}(\dot{z})\rangle\in\varphi_{b'}(\dot{x})$ be arbitrary, with $\langle p,\dot{z}\rangle\in\dot{x}$. By the inductive hypothesis, note that
    $$
    \langle b'(p),\varphi_{b'}(\dot{z})\rangle=\langle b'(p),\tau\Sigma(\dot{z})\rangle
    $$
    for some $\tau\in\mathcal{G}_J$ so that $b'(i)=\tau(i)$ for every $i\in\mathrm{supp}(\dot{z})$. We can also assume that $\tau(i)=b'(i)$ for every $i\in\mathrm{supp}(p)$, so that $b'(p)=\tau p$ by observation \ref{observation}. Because $b'$ fixes $X$ pointwise, we can further assume that $\tau\in\mathrm{fix}_{\mathcal{G}_J}(X)$. In this case, because $\langle p,\dot{z}\rangle\in\dot{x}$ and $\tau\in\mathrm{fix}_{\mathcal{G}_J}(X)$, we have that
    $$
    \langle b'(p),\varphi_{b'}(\dot{z})\rangle=\langle \tau p,\tau\Sigma(\dot{z})\rangle\in \mathrm{Orb}(\langle p,\Sigma(\dot{z})\rangle,\mathrm{fix}_{\mathcal{G}_J}(X))\subseteq \Sigma(\dot{x}),
    $$
    so $\langle b'(p),\varphi_{b'}(\dot{z})\rangle\in\Sigma(\dot{x})$, and thus $\varphi_{b'}(\dot{x})\subseteq\Sigma(\dot{x})$.

    To show that $\Sigma(\dot{x})\subseteq\varphi_{b'}(\dot{x})$, let $\langle\tau p,\tau\Sigma(\dot{z})\rangle\in\Sigma(\dot{x})$, with $\langle p,\dot{z}\rangle\in\dot{x}$ and $\tau\in\mathrm{fix}_{\mathcal{G}_J}(X)$. We wish to express $\tau=\rho\sigma$ for some $\sigma\in\mathrm{fix}_{\mathcal{G}_I}(X)$, and some $\rho\in\mathcal{G}_J$ so that $\rho(i)=b'(i)$ for every $i\in\mathrm{supp}(\sigma p)\cup\mathrm{supp}(\sigma\dot{z})$. Let $S=\mathrm{supp}(p)\cup\mathrm{supp}(\dot{z})\subseteq I$, and fix $\sigma\in\mathcal{G}_J$ so that, for every $i\in S$, $\sigma(i)=b'^{-1}(\tau(i))$. Because $b'$ fixes $X$ pointwise and $\tau\in\mathrm{fix}_{\mathcal{G}_J}(X)$, it is possible to define such a $\sigma$ in $\mathrm{fix}_{\mathcal{G}_J}(X)$. Because $S\subseteq I$ and $\mathrm{rng}(b'^{-1})=I$, it is also possible to take $\sigma\in\mathrm{fix}_{\mathcal{G}_I}(X)$. Fix $\rho\in\mathcal{G}_J$ so that $\tau=\rho\sigma$, in which case for every $i\in S$, $\tau(i)=\rho(\sigma(i))=b'(\sigma(i)).$ Thus, $\rho(i)=b'(i)$ for every $i\in\mathrm{supp}(\sigma p)\cup\mathrm{supp}(\sigma\dot{z})$.

    We can prove that $\langle\tau  p,\tau\Sigma(\dot{z})\rangle\in\varphi_{b'}(\dot{x})$, and thus prove $\Sigma(\dot{x})\subseteq\varphi_{b'}(\dot{x})$, by justifying the line
    \begin{equation}\label{varphi_b' and tauSigma}
    \langle\tau  p,\tau\Sigma(\dot{z})\rangle=\langle\rho\sigma  p,\rho\sigma\Sigma(\dot{z})\rangle=\langle\rho\sigma  p,\rho\Sigma(\sigma\dot{z})\rangle=\langle b'(\sigma p),\varphi_{b'}(\sigma\dot{z})\rangle\in\varphi_{b'}(\dot{x}).
    \end{equation}
    By observation \ref{observation}, and because $\rho(i)=b'(i)$ for every $i\in\mathrm{supp}(\sigma p)$, $\tau p=\rho\sigma p=b(\sigma p)$; thus, equality holds on the first component across line (\ref{varphi_b' and tauSigma}). By corollary \ref{Sigma commutes with G_I corollary}, because $\sigma\in\mathcal{G}_I$, $\tau\Sigma(\dot{z})=\rho\sigma\Sigma(\dot{z})=\rho\Sigma(\sigma\dot{z})$. By the inductive hypothesis, and because $\rho(i)=b'(i)$ for every $i\in\mathrm{supp}(\sigma\dot{z})$, $\rho\Sigma(\sigma\dot{z})=\varphi_{b'}(\sigma\dot{z})$. Because $\sigma\in\mathrm{fix}_{\mathcal{G}_I}(X)$, $\langle \sigma p,\sigma\dot{z}\rangle\in\dot{x}$, and so $\langle b'(\sigma p),\varphi_{b'}(\sigma\dot{z})\rangle\in\varphi_{b'}(\dot{x})$.
    \end{proof}

\begin{remark}
    Although $\varphi_b\notin M$, by theorem \ref{varphi_b in terms of Sigma, G_J} it is still the case that $\varphi_b(\dot{x})\in M$ for every $\dot{x}\in(\mathrm{HS}_{\mathcal{F}_I})^M$. Theorem \ref{varphi_b in terms of Sigma, G_J} does not imply that the entire function $\varphi_b$ is in $M$, since this would require that $M$ has a function assigning a $\pi\in\mathcal{G}_J$ to every $\dot{x}\in(\mathrm{HS}_{\mathcal{F}_I})^M$, so that $\pi(i)=b(i)$ for every $i\in\mathrm{supp}(\dot{x})$. It is straightforward to see that the existence of such a function is equivalent to the existence of the bijection $b$.
\end{remark}

\begin{lemma}\label{B restricts to HS iso lemma}
Let $M\vDash\mathrm{ZFC}$, $I\subseteq J$ and $|I|<|J|$ in $M$; let $V\vDash\mathrm{ZFC}$ be an outer model of $M$, with $M$ a class in $V$, and let $b\in V$ be a bijection from $I$ to $J$. The map $\varphi_b$ restricts to an isomorphism from $(\mathrm{HS}_{\mathcal{F}_I})^M$ to $(\mathrm{HS}_{\mathcal{F}_J})^M$.
\end{lemma}
\begin{proof}
    We will prove that $\varphi_b$ restricts to a bijection between $(\mathrm{HS}_{\mathcal{F}_I})^M$ and $(\mathrm{HS}_{\mathcal{F}_J})^M$. If this can be done, then because the map $\varphi_b:V^{\mathbb{P}(I,Q)}\to V^{\mathbb{P}(J,Q)}$ is an isomorphism, in the sense of lemma \ref{Checking that b, varphi_b are isos}, the fact that $\varphi_b((\mathrm{HS}_{\mathcal{F}_I})^M)=(\mathrm{HS}_{\mathcal{F}_J})^M$ shows that the restriction of $\varphi$ to $(\mathrm{HS}_{\mathcal{F}_I})^M$ will also be an isomorphism in the sense of lemma \ref{Checking that b, varphi_b are isos}.

    Theorem \ref{varphi_b in terms of Sigma, G_J} shows, in particular, that $\varphi_b((\mathrm{HS}_{\mathcal{F}_I})^M)\subseteq(\mathrm{HS}_{\mathcal{F}_J})^M$; for $\dot{x}\in(\mathrm{HS}_{\mathcal{F}_I})^M$, theorem \ref{varphi_b in terms of Sigma, G_J} gives us that $\varphi_b(\dot{x})=\pi\Sigma(\dot{x})$, with $\pi\in\mathcal{G}_J\in M$ and $\Sigma(\dot{x})\in(\mathrm{HS}_{\mathcal{G}_J})^M$ by lemma \ref{Sigma into HS_FJ^I}.
    
    To prove that $(\mathrm{HS}_{\mathcal{F}_J})^M\subseteq\varphi_b((\mathrm{HS}_{\mathcal{F}_I})^M)$, we will apply corollary \ref{HS_FJ in terms of HS_FI, Sigma, G_J Corollary}. Let $\dot{w}\in(\mathrm{HS}_{\mathcal{F}_J})^M$ and take $\pi\in\mathcal{G}_J$ so that for every $i\in\mathrm{supp}(\dot{w})$, $\pi(i)=b^{-1}(i)$. In this case, $\mathrm{supp}(\pi\dot{w})\subseteq I$, and by corollary \ref{HS_FJ in terms of HS_FI, Sigma, G_J Corollary},
    $$
    \dot{w}=\pi^{-1}\Sigma(\Gamma(\pi\dot{w})).
    $$
    Let $\dot{x}=\Gamma(\pi\dot{w})$. Because $\mathrm{supp}(\pi\dot{w})\subseteq I$, by lemma \ref{Gamma into HS_FI}, $\dot{x}\in\mathrm{HS}_{\mathcal{F}_I}$. We can then calculate that, by corollary \ref{Supports under Sigma, Gamma corollary}
    $$
    \mathrm{supp}(\dot{x})=\mathrm{supp}(\Gamma(\pi\dot{w}))=\mathrm{supp}(\pi\dot{w})=\pi\mathrm{supp}(\dot{w}).
    $$
    Because we took $\pi$ so that for every $i\in\mathrm{supp}(\pi\dot{w})=\mathrm{supp}(\dot{x})$, $\pi^{-1}(i)=b(i)$, we can apply theorem \ref{varphi_b in terms of Sigma, G_J} to conclude that
    $$
    \dot{w}=\pi^{-1}\Sigma(\Gamma(\pi\dot{w}))=\pi^{-1}\Sigma(\dot{x})=\varphi_b(\dot{x}).
    $$
    Since $\dot{w}\in\mathrm{HS}_{\mathcal{F}_J}$ was arbitrary, $(\mathrm{HS}_{\mathcal{F}_J})^M\subseteq\varphi_b((\mathrm{HS}_{\mathcal{F}_I})^M)$, which completes the proof.
\end{proof}

\begin{remark}
    Note the important distinction that $\varphi_b$ does \textit{not} restrict to an $\in$-preserving isomorphism between the relativized classes $M^{\mathbb{P}(I,Q)}$ and $M^{\mathbb{P}(J,Q)}$. If $\dot{A}_I$ and $\dot{A}_J$ are the respective names for the sets of mutually $Q$-generic filters added by $\mathbb{P}(I,Q)$ and $\mathbb{P}(J,Q)$, as in definition \ref{Dedekind Finite A}, then $\varphi_b(\dot{A}_I)=\dot{A}_J$. However, there is a $\mathbb{P}(I,Q)$-name $\dot{g}$ for a bijection between $\dot{A}_I$ and $\check{|I|}$; $\varphi_b(\dot{g})$ would also be a name for a bijection between $\dot{A}_J$ and $\check{|I|}$. Since $|I|<|J|$ in $M$, $\varphi_b(\dot{g})\in V^{\mathbb{P}(J,Q)}\setminus M^{\mathbb{P}(J,Q)}$.
\end{remark}

\begin{theorem}\label{KS new proof by outer model}
    Let $M\vDash\mathrm{ZFC}$, $I\subseteq J$ and $|I|<|J|$ in $M$; let $V\vDash\mathrm{ZFC}$ be an outer model of $M$, with $M$ a class in $V$, and let $b\in V$ be a bijection from $I$ to $J$. Then the map $\langle b,\varphi_b\rangle$ from lemma \ref{Checking that b, varphi_b are isos} restricts to an isomorphism in $V$ between the relativized forcing relations $\Vdash^M_{N(I,Q)}$ and $\Vdash^M_{N(J,Q)}$.
\end{theorem}
\begin{proof}
We have that $\langle b,\varphi_b\rangle$ is an isomorphism in $V$ between $\Vdash^V_{N(I,Q)}$ and $\Vdash^V_{N(J,Q)}$, in the sense of lemma \ref{Checking that b, varphi_b are isos}. Because $\varphi_b$ restricts to an isomorphism between $(\mathrm{HS}_{\mathcal{F}_I})^M$ and $(\mathrm{HS}_{\mathcal{F}_J})^M$, it is straightforward to check that $\langle b,\varphi_b\upharpoonright_{(\mathrm{HS}_{\mathcal{F}_I})^M}\rangle$ is an isomorphism between the relativized forcing relations $\Vdash^M_{N(I,Q)}$ and $\Vdash^M_{N(J,Q)}$, in the sense that for every $p\in\mathbb{P}(I,Q)$ and $\dot{x}_0,...,\dot{x}_n\in(\mathrm{HS}_{\mathcal{F}_I})^M$,
    $$
    p\Vdash^M_{N(I,Q)}\phi(\dot{x}_0,...,\dot{x}_0)\iff b(p)\Vdash^M_{N(J,Q)}\phi(\varphi_b(\dot{x}_0),...,\varphi_b(\dot{x}_n),
    $$
    and that for every $p\in\mathbb{P}(J,Q)$ and $\dot{w}_0,...,\dot{w}_n\in(\mathrm{HS}_{\mathcal{F}_J})^M$,
    $$
    p\Vdash^M_{N(J,Q)}\phi(\dot{w}_0,...,\dot{w}_0)\iff b^{-1}(p)\Vdash^M_{N(I,Q)}\phi(\varphi_b^{-1}(\dot{w}_0),...,\varphi_b^{-1}(\dot{w}_n).
    $$
\end{proof}

\begin{corollary}[Halpern, L\'evy, \cite{halpern1971boolean}, Stefanovi\'c, \cite{stefanovic2023alternatives}]\label{N(I,Q) BPI corollary.}
Theorem 4.1 holds (for every $I$).
\end{corollary}
\begin{proof}
    Let $J=\theta$ be large enough for theorem \ref{N(theta,Q) Models BPI} to apply. Working in the $\mathrm{Coll}(I,\theta)$-extension $M[H]$ (in which $M$ is a class), with $b\in M[H]$ a bijection from $I$ to $\theta$, the isomorphism $\langle b,\varphi_b\rangle$ from theorem \ref{KS new proof by outer model} can be used to transfer witnesses $q,y'$ to the filter extension property, with $K=\mathrm{fix}_{\mathcal{G}_\theta}(F\cup X)$, for the symmetric system $\langle\mathbb{P}(\theta,Q),\mathcal{G},\mathcal{F}\rangle$, to witnesses $b^{-1}(q),b^{-1}(y')$ to the filter extension property, with $K=\mathrm{fix}_{\mathcal{G}_I}(b^{-1}(F)\cup b^{-1}(X))$, for the symmetric system $\langle\mathbb{P}(I,Q),\mathcal{G},\mathcal{F}\rangle$.
\end{proof}

\begin{remark}
    Given the proof of corollary \ref{N(I,Q) BPI corollary.}, which relies on theorems \ref{N(theta,Q) Models BPI} and \ref{KS new proof by outer model}, note the broader potential of Ramsey-theoretic methods in models without choice. If a model is more or less ``agnostic" about the size of a given set, one can treat the set as arbitrarily large to apply a Ramsey-theoretic argument.
\end{remark}

\begin{corollary}[Karagila, Schlicht, \cite{karagila2020have}]
    Let $M\vDash\mathrm{ZFC}$ and let $I,J\in M$ be infinite. For every $\mathbb{P}(I,Q)$-generic $G$, there is a $\mathbb{P}(J,Q)$-generic $F$ so that $N^G(I,Q)=N^F(J,Q)$.
\end{corollary}
\begin{proof}
    Let $G$ be an arbitrary $\mathbb{P}(I,Q)$-generic filter over $M$. Let $H$ be $\mathrm{Coll}(I,J)$-generic over $M$, and let $V$ be any outer model for which both $G,H\in V$, and in which $M$ is a class. Let $b\in V$ be an bijection from $I$ to $J$. Extend $b$ to the isomorphism between $\mathbb{P}(I,Q)$ and $\mathbb{P}(J,Q)$ from definition \ref{b and varphi_b definition}. From this isomorphism, it is straightforward to show that $b(G)$ is $\mathbb{P}(J,Q)$-generic over $M$.

    We then claim that for every $\dot{x}\in\mathrm{HS}_{\mathcal{F}_I}$, $\dot{x}^G=\varphi_b(\dot{x})^F$. If this can be shown, the fact that $\varphi_b\upharpoonright_{(\mathrm{HS}_{\mathcal{F}_I})^M}$ surjects onto $(\mathrm{HS}_{\mathcal{F}_J})^M$ shows the desired result that
    $$
    N^G(I,Q)=N^F(J,Q).
    $$

    Assume inductively that for every $\dot{z}\in\mathrm{rng}(\dot{x})$, $\dot{z}^G=\varphi_b(\dot{z})^F.$ Using theorem \ref{KS new proof by outer model}, we have the following sequence of equivalences in $V$:
    \begin{align*}
        \dot{z}^G\in\dot{x}^G&\iff(\exists p\in G)p\Vdash^M_{N(I,Q)}\dot{z}\in\dot{x}\\&\iff(\exists b(p)\in b(G)=F)b(p)\Vdash^M_{N(J,Q)}\varphi_b(\dot{z})\in\varphi_b(\dot{x})\\&\iff\varphi_b(\dot{z})^F\in\varphi_b(\dot{x})^F.
    \end{align*}
    By the inductive hypothesis, $\dot{z}^G=\varphi_b(\dot{z})^F$.
    By the definition of $\varphi_b(\dot{x})$, every element of $\varphi_b(\dot{x})^F$ is of the form $\varphi_b(\dot{z})^F$, for some $\dot{z}\in\mathrm{rng}(\dot{x})$. We have thus shown that $\dot{x}^G$ and $\varphi_b(\dot{x})^F$ have the same elements, and so they are equal.
\end{proof}

\begin{theorem}\label{Sigma Elementary Embedding Theorem}
    For every $p\in\mathbb{P}(I,Q)$, every $\dot{x}_0,...,\dot{x}_n\in\mathrm{HS}_{\mathcal{F}_I}$, and every formula $\phi$,
    $$
    p\Vdash_{N(I,Q)}\phi(\dot{x}_0,...,\dot{x}_n)\iff p\Vdash_{N(J,Q)}\phi(\Sigma(\dot{x}_0),...,\Sigma(\dot{x}_n))
    $$
\end{theorem}
\begin{proof}
Let $G$ be $\mathrm{Coll}(I,J)$-generic, and let $V=M[G]$. Fix a $p\in\mathbb{P}(I,Q)$, $\dot{x}_0,...,\dot{x}_n\in\mathrm{HS}_{\mathcal{F}_I}$, and a formula $\phi$. Let $X=\mathrm{supp}(p)\cup\bigcup_{i\leq n}\mathrm{supp}(\dot{x})$. There is a bijection $b_X\in V$ that fixes $X$ pointwise, in which case $b_X(p)=p$ and by theorem \ref{varphi_b in terms of Sigma, G_J}, for every $i\leq n$, $\varphi_{b_X}(\dot{x}_i)=\Sigma(\dot{x})$. Applying theorem \ref{KS new proof by outer model} results in the conclusion of the theorem for $p,\dot{x}_0,...,\dot{x}_n,$ and $\phi$, which were arbitrary.
\end{proof}

\begin{corollary}\label{Sigma* elementary embedding}
    Let $M\vDash\mathrm{ZFC}$, let $I\subseteq J$ be infinite sets in $M$, and let $F$ be $\mathbb{P}(J,Q)$-generic over $M$. Then $G=\{p\in F\ |\ \mathrm{supp}(p)\subseteq I\}$ is $\mathbb{P}(I,Q)$-generic over $M$, $N^G(I,Q),N^F(J,Q)\subseteq M[F]$, and the map
    $$
    \Sigma^*:N^G(I,Q)\to N^F(J,Q)
    $$
    in $M[F]$, defined by
    $$
    \Sigma^*(\dot{x}^G)=\Sigma(\dot{x})^F
    $$
    is well-defined and an elementary embedding.
\end{corollary}
\begin{proof}
    The claim that $G$ is $\mathbb{P}(I,Q)$-generic over $M$ is straightforward to verify. To prove that $\Sigma^*$ is well-defined, let $\dot{x}_0,\dot{x}_1\in(\mathrm{HS}_{\mathcal{F}_I})^M$ so that $\dot{x}_0^G=\dot{x}_1^G$. Then there is a $p\in G\subseteq\mathbb{P}(I,Q)$ so that
    $$
    p\Vdash_{N(I,Q)}^M\dot{x}_0=\dot{x}_1.
    $$
    By theorem \ref{Sigma Elementary Embedding Theorem}, we get that
    $$
    p\Vdash_{N(J,Q)}^M\Sigma(\dot{x}_0)=\Sigma(\dot{x}_1).
    $$
    Since $G\subseteq F$, $p\in F$, so in $M[F]$, $\Sigma^*(\dot{x}_0^G)=\Sigma(\dot{x}_0)^F=\Sigma(\dot{x}_1)^F=\Sigma^*(\dot{x}_1^G)$, making $\Sigma^*$ well-defined.

    Nearly the same argument works to show that $\Sigma^*$ is an elementary embedding. Let $\dot{x}_0^G,...,\dot{x}_n^G\in N^G(I,Q)$, and let $\phi$ be an arbitrary formula. The corollary then follows by justifying that
    \begin{align*}
        N^G(I,Q)\vDash\phi(\dot{x}_0^G,...,\dot{x}_n^G)&\iff(\exists p\in G)p\Vdash_{N(I,Q)}^M\phi(\dot{x}_0,...,\dot{x}_n)\\&\iff(\exists q\in F) q\Vdash_{N(J,Q)}^M\phi(\Sigma(\dot{x}_0),...,\Sigma(\dot{x}_n))\\&\iff N^F(J,Q)\vDash\phi(\Sigma^*(\dot{x}_0)^F,...,\Sigma^*(\dot{x}_n)^F).
    \end{align*}

    The second equivalence is the only one that needs further justification. To prove the direction $(\Rightarrow)$, note that $p\in G\subseteq F$, so by theorem \ref{Sigma Elementary Embedding Theorem} we can take $p=q\in F$ so that the latter statement holds. To prove the direction $(\Leftarrow)$, note that for every $i\leq n$, $\mathrm{supp}(\Sigma(\dot{x}_i))\subseteq I$ by lemma \ref{Sigma into HS_FJ^I}. By lemma \ref{Support Restriction Lemma}, taking $p=q\upharpoonright_I\in F$, we have that
    $$
    p\Vdash^M_{N(J,Q)}\phi(\Sigma(\dot{x}_0),...,\Sigma(\dot{x}_n)),
    $$
    and that $p\in G$. By theorem \ref{Sigma Elementary Embedding Theorem}, we then have that
    $$
    p\Vdash_{N(I,Q)}^M\phi(\dot{x}_0,...,\dot{x}_n).
    $$
    
\end{proof}

\subsection{An Alternate Proof That $N(I,Q)\vDash\mathrm{BPI}$}

We can conclude this section by noting that a case of our main argument answers the $8^\mathrm{th}$ question listed on Karagila's website, which outlines an alternate proof that $N(I,Q)\vDash\mathrm{BPI}$ based on theorems \ref{Blass SVC BPI Theorem} and \ref{SVC in sym ext}, both of which appear in \cite{blass1979injectivity}. Blass notes that theorem \ref{Blass SVC BPI Theorem}, and its proof, were first shown to him ``in a slightly different context" by Pincus; this seems to correspond to 4.1.8$^*$ in \cite{pincus1974independence}.

\begin{theorem}[Halpern, L\'evy, \cite{halpern1971boolean}, Stefanovi\'c, \cite{stefanovic2023alternatives}]\label{Karagila 8}
    For every $I,Q$, $N(I,Q)\vDash\mathrm{BPI}$.
\end{theorem}
\begin{proof}
    By theorem \ref{SVC in sym ext}, $N(I,Q)\vDash\mathrm{SVC}$; by theorem \ref{Blass SVC BPI Theorem} there is a particular set $S\in N(I,Q)$, and a particular filter $F'=\{\{p\in S^{<\omega}\ |\ s\in\mathrm{rng}(p)\}\ |\ s\in S\}\in N(I,Q)$ on $S^{<\omega}$, so that any ultrafilter $U'\in N(I,Q)$ extending $F'$ can be used to prove $\mathrm{BPI}$ uniformly in $N(I,Q)$. Karagila's $8^\mathrm{th}$ question asks whether one can describe such an ultrafilter in the Cohen model. Our methods, in particular, describe such an ultrafilter in the generalized Cohen model $N(I,Q)$, from which $\mathrm{BPI}$ follows.
\end{proof}

We can also observe that a more precise version of theorem \ref{Halpern-Levy FEP Theorem} is available by connecting the results from \cite{halpern1971boolean}/\cite{stefanovic2023alternatives} and \cite{blass1979injectivity}. Recall that, in general, $\mathrm{supp}(\dot{x})$ may be a proper subgroup of $\mathrm{sym}_\mathcal{G}(\dot{x})$.

\begin{corollary}\label{Halpern-Levy FEP General}
Theorem \ref{Halpern-Levy FEP Theorem} holds when $K=\mathrm{fix}_\mathcal{G}(F)\cap\mathrm{fix}_\mathcal{G}(X)$, with $F=\mathrm{supp}(\dot{F})$ and $X=\mathrm{supp}_\mathcal{G}(X)$, is replaced with $K=\mathrm{sym}_\mathcal{G}(\dot{F})\cap\mathrm{sym}_\mathcal{G}(\dot{x})$.
\end{corollary}
\begin{proof}
    Let $\dot{S}$ and $\dot{F}'$ be names for the sets $S$ and $F'$ in the proof of theorem \ref{Karagila 8}. From our proof sketch of theorem \ref{SVC in sym ext}, it is possible to take $\dot{S}$ so that $\mathrm{sym}_\mathcal{G}(\dot{S})=\mathcal{G}$. In this case, by the definition of $F'$, it is also possible to take $\dot{F}'$ so that $\mathrm{sym}_\mathcal{G}(\dot{F}')=\mathcal{G}$. By theorem \ref{Halpern-Levy FEP Theorem} and our proof of theorem \ref{fep implies bpi in sym ext}, there is a name $\dot{U}'$ for an ultrafilter on $S^{<\omega}$ extending $F'$, with $\dot{U}'=\mathcal{G}$.

    By our proof of theorem \ref{BPI equivalent to FEP for sym ext}, given $\dot{F},\dot{x},\dot{y}\in\mathrm{HS}_\mathcal{F}$, as in theorem \ref{Halpern-Levy FEP Theorem}, there is a name for an ultrafilter $\dot{U}$ on $\dot{x}$ that extends $\dot{F}$, with
    $$
    \mathrm{sym}_{\mathcal{G}}(\dot{F})\cap\mathrm{sym}_\mathcal{G}(\dot{x})\cap\mathrm{sym}_\mathcal{G}(\dot{U}')\preceq\mathrm{sym}_\mathcal{G}(\dot{U}).
    $$
Because $\dot{U}$ is a name for an ultrafilter on $\dot{x}$, for every $q\leq p$, there is a $y'\in\{\dot{y},\dot{x}\setminus\dot{y}\}$ so that $\langle q,y'\rangle\in \dot{U}$. The corollary follows from the fact that $\mathrm{sym}_\mathcal{G}(\dot{U}')=\mathcal{G}$, in which case
$$
\dot{F}\cup\mathrm{Orb}(\langle q,y'\rangle,\mathrm{sym}_\mathcal{G}(\dot{F})\cap\mathrm{sym}_\mathcal{G}(\dot{x}))\subseteq\dot{U},
$$
and so
$$
1\Vdash``\mathrm{Orb}(\langle q,y'\rangle,\mathrm{sym}_\mathcal{G}(\dot{F})\cap\mathrm{sym}_\mathcal{G}(\dot{x}))\subseteq2^{\dot{x}}\mathrm{\ has\ the\ FIP}."
$$
\end{proof}

\section{The (Virtual) Ramsey Property}
To prove the general case of theorem \ref{Halpern-Levy FEP Theorem}, we handled the symmetric system $\langle\mathbb{P}(I,Q),\mathcal{G},\mathcal{F}\rangle$ differently when $I$ was large compared to when $I$ was small. In this section, we can formally justify this difference in treatment. We will do so by introducing a dynamical property for symmetric systems that implies $\mathrm{BPI}$ and generalizes our main argument from section 4, but provably fails for the symmetric system $\langle\mathbb{P}(\omega,2^{<\omega}),\mathcal{G},\mathcal{F}\rangle$ that presents the Cohen model.

This property can be motivated by remark \ref{search and shift remark}, where we call our main proof from section 4 a ``search and shift" argument. We wish to captures the ability to prove instances of the filter extension property, for names $\dot{F},\dot{x}\in\mathrm{HS}_{\mathcal{F}}$, by a search and shift argument using \textit{any} $\mathbb{P}$-name $\dot{U}$ for an ultrafilter on $\dot{x}$ that extends $\dot{F}$.

In subsection 6.1, we show that the Ramsey property from \cite{blass1986prime} captures this ability for permutation models. Given a filter $F$ on a set $x$, with both in $M(A,\mathcal{G},\mathcal{F})$, we directly prove that the Ramsey property guarantees the success of a search and shift argument through an arbitrary ultrafilter $U\in M\vDash\mathrm{ZFCA}$ on $x$ that extends $F$. In this manner, the Ramsey property directly implies the filter extension property, and thus $\mathrm{BPI}$. This reproves one direction of the equivalence from \cite{blass1986prime} between $\mathrm{BPI}$ in a permutation model $M(A,\mathcal{G},\mathcal{F})$ and the Ramsey property on $\mathcal{F}$. By the other direction of this equivalence, our proof shows that $\mathrm{BPI}$ in a permutation model can always be explained in terms of ultrafilters in the outer model $M\vDash\mathrm{ZFC}$ coordinating the ability to extend symmetric filters in $M(A,\mathcal{G},\mathcal{F})$.

In subsection 6.2, we state the virtual Ramsey property for symmetric systems. This is a dynamical property that generalizes the search and shift argument used to prove theorem \ref{N(theta,Q) Models BPI}, and generalizes the Ramsey property from \cite{blass1986prime}. We confirm that the virtual Ramsey property allows for a search and shift argument, using arbitrary names for ultrafilters, to prove instances of the filter extension property. Thus, we show that $\mathrm{BPI}$ can be proven in (all) permutation models and in (a class of) symmetric extensions by the same direct, dynamical argument: the (virtual) Ramsey property directly implies the filter extension property, which implies $\mathrm{BPI}$.

In subsection 6.3, we confirm that $\langle\mathbb{P}(I,Q),\mathcal{G},\mathcal{F}\rangle$ has the virtual Ramsey property when $I$ is large. This formalizes the notion that our main proof in section 4 is essentially dynamical, and verifies that the Ramsey property is not vacuous. We then justify the need for a new method of proof in section 5, by showing that the virtual Ramsey property fails for the symmetric system $\langle\mathbb{P}(\omega,2^{<\omega}),\mathcal{G},\mathcal{F})$, despite the fact that this symmetric system has the filter extension property and $\mathrm{BPI}$ holds in the Cohen model. We conclude by posing a question about the relationship between the virtual Ramsey property and $\mathrm{BPI}$ in symmetric extensions.

\subsection{The Ramsey Property, the Filter Extension Property, and $\mathrm{BPI}$ in Permutation Models}
Blass introduced the Ramsey property in \cite{blass1986prime}, but it will be more convenient to use the equivalent formulation of the property given in \cite{blass2011partitions}. If the reader is unfamiliar with the Ramsey property, they may find it helpful to skip ahead to our proof of theorem \ref{RP,FEP, BPI in perm models}. In this proof, definition \ref{RP definition} can be recognized as the exact property needed to successfully search through \textit{any} ultrafilter $U$ on a set $x\in M(A,\mathcal{G},\mathcal{F})$ extending a filter $F\in M(A,\mathcal{G},\mathcal{F})$, then shift the results of the search by some $\sigma\in K\preceq\mathrm{sym}_\mathcal{G}(F)\cap\mathrm{sym}_\mathcal{G}(x)$ to prove instances of the filter extension property.

\begin{definition}[Blass, \cite{blass1986prime}, \cite{blass2011partitions}]\label{RP definition}\noindent Let $\mathcal{G}$ be a group.
    \begin{enumerate}[label=(\roman*)]
        \item A group action $\Gamma$ of $\mathcal{G}$ on a set $A$ has the Ramsey property if for every coloring $c:A\to 2$ and for every $\Delta\in[A]^{<\omega}$, there is a $\sigma\in G$ so that the restriction $c\upharpoonright_{\sigma\Delta}$ is monochromatic.
        \item A subgroup $H\preceq \mathcal{G}$ is called a Ramsey subgroup of $\mathcal{G}$ if the action of $\mathcal{G}$ on the cosets $\mathcal{G}/H$ (by left multiplication) has the Ramsey property. 
        \item A normal filter $\mathcal{F}$ of subgroups of a group $\mathcal{G}$ is said to have the Ramsey property if $\mathcal{F}$ has a base of subgroups $K\in \mathcal{F}$ so that every subgroup $H$ of $K$ in $\mathcal{F}$ is a Ramsey subgroup of $K$.
    \end{enumerate}
\end{definition}

The filter extension property and its relationship to the Ramsey property are not identified explicitly in \cite{blass1986prime}. By our observations in section 3, their relationship is implicit in Blass' proofs that use Halpern's contradiction framework, so we attribute the entire theorem \ref{RP,FEP, BPI in perm models} to Blass.

\begin{theorem}[Blass, \cite{blass1986prime}]\label{RP,FEP, BPI in perm models}
    If $M\vDash\mathrm{ZFCA}$, the following are equivalent:
    \begin{enumerate}[label=(\roman*)]
        \item $M(A,\mathcal{G},\mathcal{F})\vDash\mathrm{BPI}$,
        \item $M(A,\mathcal{G},\mathcal{F})$ has the filter extension property,
        \item $\mathcal{F}$ has the Ramsey property.
    \end{enumerate}
\end{theorem}
\begin{proof}
    The equivalence $(i)\iff(iii)$ is proven explicitly in \cite{blass1986prime}. The implication $(ii)\Rightarrow(i)$ is theorem \ref{FEP proves BPI in permutation models}. In Blass' contradiction proof that the Ramsey property implies $\mathrm{BPI}$, he considers a filter $U\in M(A,\mathcal{G},\mathcal{F})$ on $x$ that extends $F$ and is maximal with respect to inclusion and the property that $K\preceq\mathrm{sym}_\mathcal{G}(U)$, where $K\in\mathcal{F}$ witnesses the Ramsey property. He then argues by contradiction that $U$ is an ultrafilter on $x$ in $M(A,\mathcal{G},\mathcal{F})$.
    
    One can replace $U$ in Blass' proof with an arbitrary filter $F$, for which $K\preceq \mathrm{sym}_\mathcal{G}(F)\cap\mathrm{sym}_\mathcal{G}(x)$ is a witness of the Ramsey property. Then Blass' proof can be used to show, by contradiction, that $K$ also witnesses the filter extension property. In this case, $\mathcal{F}$ has a base of witnesses for the filter extension property. Thus, $(iii)\Rightarrow(ii)$ holds.
\end{proof}

We can now give a direct proof of $(i)\Rightarrow(ii)$ from the theorem, using a search and shift argument through an arbitrary ultrafilter $U$ on $x\in M(A,\mathcal{G},\mathcal{F})$ that extends $F\in M(A,\mathcal{G},\mathcal{F})$.

\begin{proof}[Proof of $(iii)\Rightarrow(ii)$]
    Let $K\in\mathcal{F}$ witness the Ramsey property. We will show that $K$ also witnesses the filter extension property; then, the base of subgroups in $\mathcal{F}$ witnessing the Ramsey property also witnesses the filter extension property. Let $F,x,y\in M(A,\mathcal{G},\mathcal{F})$, $K\preceq\mathrm{sym}_\mathcal{G}(F)\cap\mathrm{sym}_\mathcal{G}(x)$, and
    $$
    M(A,\mathcal{G},\mathcal{F})\vDash``F\subseteq2^x\mathrm{\ has\ the\ FIP}"\wedge y\subseteq x.
    $$
    Let $H=K\cap \mathrm{sym}_\mathcal{G}(y)=K\cap \mathrm{sym}_\mathcal{G}(x\setminus y)$. Because $K$ is a witness of the Ramsey property and $H$ is a subgroup of $K$ in $\mathcal{F}$, $H$ is a Ramsey subgroup of $K$.

    Let $U\in M$ be an ultrafilter on $x$ extending $F$; $U$ itself is not necessarily hereditarily $\mathcal{F}$-symmetric, but exists in $M\vDash\mathrm{ZFCA}$. Define the coloring $c:K/H\to\{\mathrm{green,\ red}\}$ by setting
    $$
    c(\pi H)=
    \begin{cases}
        &\mathrm{green\ if\ }\pi y\in U\\
        &\mathrm{red\ if\ }\pi (x\setminus y)\in U.
    \end{cases}
    $$
    Note that $c(\pi H)$ is well-defined, since every $h\in H$ fixes both $y$ and $x\setminus y$. For every $\Delta=\{\pi_0H,...,\pi_nH\}\in[K/H]^{<\omega}$, the fact that $H$ is a Ramsey subgroup of $K$ guarantees a $\sigma\in K$ so that $c\upharpoonright\sigma \Delta$ is monochromatic. If $c$ colors the elements of $\sigma \Delta=\{\sigma\pi_0H,...,\sigma\pi_nH\}$ green, then
    \begin{equation}
    (F\cup\{\sigma\pi_iy|i\leq n\})\subseteq U.
    \end{equation}
    Because $\sigma,\sigma^{-1}\in K\subseteq\mathrm{fix}(F)\cap\mathrm{fix}(x)$, we can apply $\sigma^{-1}$ to conclude that
    $$
    (F\cup\{\pi_iy|i\leq n\})\subseteq \sigma^{-1}U.
    $$
    By symmetry (lemma \ref{automorphisms of permutation model lemma}), $\sigma^{-1}U$ is still an ultrafilter on $\sigma^{-1}x=x$ in $M$, so we can conclude that
    $$
    M(A,\mathcal{G},\mathcal{F})\vDash``F\cup\{\pi_iy|i\leq n\}\subseteq2^x\mathrm{\ has\ the\ FIP.}"
    $$
    If $c$ colors the elements of $\sigma \Delta=\{\sigma\pi_0H,...,\sigma\pi_nH\}$ red, the same argument allows us to conclude
    $$
    M(A,\mathcal{G},\mathcal{F})\vDash``F\cup\{\pi_i(x\setminus y)|i\leq n\}\subseteq2^x\mathrm{\ has\ the\ FIP.}"
    $$
    
    Now define a coloring $C:[K/H]^{<\omega}\to\{\mathrm{green,\ red}\}$ by
    $$
    C(\Delta)=
    \begin{cases}
        &\mathrm{green\ if\ }(\exists\sigma\in K)c\upharpoonright \sigma \Delta\mathrm{\ is\ green.}\\
        &\mathrm{red\ otherwise.}
    \end{cases}
    $$
    Ordering $[K/H]^{<\omega}$ by reverse-inclusion, $[K/H]^{<\omega}$ is downward-directed by unions. This ensures that $C$ colors a dense set of its nodes to be either mutually green or mutually red. In the former case, we claim that $F\cup\mathrm{Orb}(y,K)$ has the FIP; in the latter case, we claim that $F\cup\mathrm{Orb}(x\setminus y,K)$ has the FIP.

    Without loss of generality, we prove the former claim. For $\pi_0,...,\pi_n\in K$ arbitrary, we aim to show that $F\cup\{\pi_iy|i\leq n\}$ has the FIP. Let $\Delta=\{\pi_0H,...,\pi_nH\}\in [K/H]^{<\omega}$. By the assumption that $C$ colors a dense set of nodes in $[K/H]^{<\omega}$ green, we can find some $\Delta'\supseteq \Delta$ so that $C(\Delta')=\mathrm{green}$. In other words, there is a $\sigma\in K$ so that $c\upharpoonright \sigma \Delta'$ is green, and so
    $$
    M(A,\mathcal{G},\mathcal{F})\vDash``F\cup\{\pi y|\pi H\in \Delta'\}\mathrm{\ has\ the\ FIP}."
    $$
    Since $\Delta\subseteq \Delta'$, this gives the desired result that
    $$
    M(A,\mathcal{G},\mathcal{F})\vDash``F\cup\{\pi_i y|i\leq n\}\mathrm{\ has\ the\ FIP}."
    $$
    Since $\pi_0,...,\pi_n\in K$ were arbitrary, we can conclude that
    $F\cup\mathrm{Orb}(y,K)$ has the FIP. Because $F,x,y\in M(A,\mathcal{G},\mathcal{F})$ were taken to arbitrarily meet the hypothesis of the filter extension property, with $K\preceq\mathrm{sym}_\mathcal{G}(F)\cap\mathrm{sym}_\mathcal{G}(x)$, we can conclude that $K$ witnesses the filter extension property.
\end{proof}

\begin{remark}
In the same sense as in remark \ref{search and shift remark}, the Ramsey property is a dynamical condition that allows ultrafilters in the outer model $M\vDash\mathrm{ZFCA}$ to coordinate symmetric filter extension in the inner model $M(A,\mathcal{G},\mathcal{F})$. By the implication $(i)\Rightarrow(iii)$ from theorem , originally from \cite{blass1986prime}, whenever $\mathrm{BPI}$ holds in $M(A,\mathcal{G},\mathcal{F})$ it can be attributed to this mechanism of inheritance from the outer model $M$.
    
\end{remark}

\subsection{The Virtual Ramsey Property, the Filter Extension Property, and $\mathrm{BPI}$ in Symmetric Extensions}
To generalize the Ramsey property to a symmetric system $\langle\mathbb{P},\mathcal{G},\mathcal{F}\rangle$, recall that we would like to capture the ability to verify the filter extension property by using a search and shift argument with \textit{any} $\mathbb{P}$-name $\dot{U}$ for an ultrafilter on $\dot{x}\in\mathrm{HS}_\mathcal{F}$ that extends $\dot{F}\in\mathrm{HS}_\mathcal{F}$. To motivate this definition, it will be helpful to consider the main obstacle in directly generalizing our proof for permutation models.

The main ingredient in our proof of theorem \ref{RP,FEP, BPI in perm models}$(iii)\Rightarrow(ii)$ was that, if $U$ is an ultrafilter on $x$ extending $F$ (with $x,F\in M(A,\mathcal{G},\mathcal{F})$), then for every $\Delta\in [\mathrm{Orb}(y',K)]^{<\omega}$, we can identify some $\sigma\in K\preceq\mathrm{sym}_\mathcal{G}(F)\cap\mathrm{sym}_\mathcal{G}(x)$ for which
$$
F\cup\sigma \Delta\subseteq U.
$$
For the symmetric system $\langle\mathbb{P}(\theta,Q),\mathcal{G},\mathcal{F}\rangle$ from theorem \ref{N(theta,Q) Models BPI}, this is generally not possible. For for an arbitrary $\mathbb{P}(\theta,Q)$-name $\dot{U}$ for an ultrafilter on $\dot{x}\in\mathrm{HS}_{\mathcal{F}}$ extending $\dot{F}\in\mathrm{HS}_\mathcal{F}$, it may be the case that for every $q\leq p$, $y'\in\{\dot{y},\dot{x}\setminus\dot{y}\}$, $K\in\mathcal{F}$, and $\sigma\in K$, there is a $\Delta\in[\mathrm{Orb}(\langle q,y'\rangle,K)]^{<\omega}$ so that
\begin{equation}\label{failure line}
\dot{F}\cup\sigma \Delta\not\subseteq\dot{U}.
\end{equation}

To address this obstacle, recall two distinctions between our proof of theorem \ref{N(theta,Q) Models BPI} and our proof of theorem \ref{RP,FEP, BPI in perm models}$(iii)\Rightarrow(ii)$. First, to prove that
$$
1\Vdash``\dot{F}\cup\mathrm{Orb}(\langle q,y'\rangle,K)\subseteq2^{\dot{x}}\mathrm{\ has\ the\ FIP},"
$$
we only considered subsets $\Delta=\{\langle\pi p,\pi y'\rangle\ |\ \pi\in\delta\}\in[\mathrm{Orb}(\langle q,y'\rangle,K)]^{<\omega}$ that corresponded to sets $\delta\in[K]^{<\omega}$ for which the conditions in $\{\pi q\ |\ \pi\in\delta\}$ are mutually compatible. This makes it so the names in $\{\pi y'\ |\ \pi\in\delta\}$ can mutually appear in an interpretation of $\Delta\subseteq\mathrm{Orb}(\langle q,y'\rangle,K)$. Second, we did not aim to show that there was a $\sigma\in K$ for which
$$
\dot{F}\cup\sigma \Delta\subseteq\dot{U}.
$$
Instead, in line (\ref{Refined Forcing Search Conclusion, Main Proof}), we showed that there was a condition $r_\delta=\bigcup_{\vec{\alpha}\in S}\pi_{\vec{\alpha}}q_{\vec{\alpha}}$ and a $\sigma\in K$, for the right choice of $S\subseteq\prod_{i\in Y\setminus(F\cup X)}H^m_i$, so that
$$
r_\delta\Vdash\dot{F}\cup\{\langle1,\sigma\pi y'\rangle\ |\ \pi\in\delta\}\subseteq\dot{U}.
$$
We then showed that $r_\delta$ could be restricted to $\sigma\bigcup_{\pi\in \delta}\pi q$, to deduce from the line above that
$$
\sigma\bigcup_{\pi\in \delta}\pi q\Vdash``\dot{F}\cup\{\langle1,\sigma\pi y'\rangle\ |\ \pi\in\delta\}\subseteq2^{\dot{x}}\mathrm{\ has\ the\ FIP}."
$$
From here, we shift back by $\sigma^{-1}$ to get the desired conclusion. In other words, we don't mind the failure in line (\ref{failure line}) -- we only need to find information from $\dot{U}$ that can be used to \textit{deduce} that an instance of the FIP is forced.

By formally stating these two features of our proof of theorem \ref{N(theta,Q) Models BPI}, we can capture the desired ability to prove instances of the filter extension property by a search and shift argument that uses arbitrary names for ultrafilters. We call this the virtual Ramsey property for symmetric systems, which we state in definition \ref{Virtual Ramsey property definition}. Definition \ref{Reducible definition} formalizes a notion of ``deduction," as we use it in the sketch above, and lemma \ref{reducibility lemma in N(I,Q)} verifies the connection between definition \ref{Reducible definition} and the sketch.

\begin{definition}\label{Reducible definition}
    Let $\langle\mathbb{P},\mathcal{G},\mathcal{F}\rangle$ be a symmetric system, $H\in\mathcal{F}$, and $p,q\in\mathbb{P}$. We say that $p$ is $H$-reducible to $q$ if whenever $H\preceq \bigcap_{i\leq n}\mathrm{sym}_\mathcal{G}(\dot{x}_i)$, for $\mathbb{P}$-names $\dot{x}_0,...,\dot{x}_n$, it is the case that
    $$
    p\Vdash\varphi(\dot{x}_0,...,\dot{x}_n)\Rightarrow q\Vdash\varphi(\dot{x}_0,...,\dot{x}_n).
    $$

\end{definition}

\begin{lemma}\label{reducibility lemma in N(I,Q)}
    Let $\langle\mathbb{P}(I,Q),\mathcal{G},\mathcal{F}\rangle$ be a symmetric system for a generalized Cohen model. Let $\dot{F},\dot{x},\dot{y}$ meet the hypothesis for the filter extension property, and have the supports $F,X,Y\in[I]^{<\omega}$. Let $H=\mathrm{fix}_{\mathcal{G}}(F\cup X\cup Y)$ and $K=\mathrm{fix}_\mathcal{G}(F\cup X)$, so that $H,K\in\mathcal{F}$ and $H\preceq K$. For every $\delta\in[K]^{<\omega}$ and $\sigma\in K$, $r_\delta$ is $\bigcap_{\pi\in \delta}(\sigma\pi)H(\sigma\pi)^{-1}$-reducible to $\bigcup_{\pi\in\delta}\sigma\pi q$ if $r_\delta\upharpoonright_{E}=\bigcup_{\pi\in\delta}\sigma\pi q$ for 
    $$
    E=\bigcup_{\pi\in \delta}\sigma\pi (F\cup X\cup Y)\in[I]^{<\omega}.
    $$
    In this case, if $\delta=\{\pi_0,...,\pi_n\}$ and $y'\in\{\dot{y},\dot{x}\setminus\dot{y}\}$, then for every formula $\phi$,
    $$
    r_\delta\Vdash\phi(\dot{F},\dot{x},\sigma\pi_0y',...,\sigma\pi_ny')\Rightarrow \bigcup_{\pi\in\delta}\sigma\pi q\Vdash\phi(\dot{F},\dot{x},\sigma\pi_0y',...,\sigma\pi_ny').
    $$
    
\end{lemma}
\begin{proof}
    Note first that
    \begin{equation}
    (\sigma\pi)\mathrm{fix}_\mathcal{G}(F\cup X\cup Y)(\sigma\pi)^{-1}=\mathrm{fix}_\mathcal{G}(\sigma\pi(F\cup X\cup Y)),
    \end{equation}
    and consequently that
    \begin{equation}
    \bigcap_{\pi\in \delta}(\sigma\pi)\mathrm{fix}_\mathcal{G}(F\cup X\cup Y)(\sigma\pi)^{-1}=\mathrm{fix}_\mathcal{G}(\bigcup_{\pi\in \delta}\sigma\pi (F\cup X\cup Y))=\mathrm{fix}_\mathcal{G}(E).
    \end{equation}
    If $\mathrm{fix}_\mathcal{G}(E)\preceq\bigcap_{i\leq n}\mathrm{sym}_\mathcal{G}(\dot{x}_i)$, then by applying lemma \ref{Support Restriction Lemma}, we can conclude that $r_\delta$ is $\mathrm{fix}_{\mathcal{G}}(E)$-reducible to $r_\delta\upharpoonright_{E}=\bigcup_{\pi\in\delta}\sigma\pi q$.

    Because $\sigma\in K,\delta\subseteq K$, and $K=\mathrm{fix}_\mathcal{G}(F\cup X)$,
    $$
    E=\bigcup_{\pi\in\delta}\sigma\pi(F\cup X\cup Y)=F\cup X\cup \bigcup_{\pi\in \delta}\sigma\pi Y.
    $$
    Thus, because each $\mathrm{supp}(\sigma\pi y')=\sigma\pi Y$, $\mathrm{fix}_\mathcal{G}(E)\preceq\mathrm{sym}_\mathcal{G}(\dot{F})\cap\mathrm{sym}_\mathcal{G}(\dot{x})\cap\bigcap_{\pi\in\delta}\mathrm{sym}_\mathcal{G}(\sigma\pi y')$, and so the remainder of the lemma follows from definition \ref{Reducible definition}.
\end{proof}

If definition \ref{Virtual Ramsey property definition} appears overly technical, the reader may find it helpful to skip ahead to read the proof of theorem \ref{RP, FEP, BPI for symmetric extensions}$(iii)\Rightarrow(ii)$. This proof generalizes our proof of theorem \ref{N(theta,Q) Models BPI}, and our proof of theorem \ref{RP,FEP, BPI in perm models}$(iii)\Rightarrow(ii)$ -- definition \ref{Virtual Ramsey property definition} is stated with only this application in mind.

\begin{definition}\label{Virtual Ramsey property definition}
Let $\mathbb{P}$ be a poset and $H,K\preceq\mathrm{Aut}(\mathbb{P})$.
    \begin{enumerate}[label=(\roman*)]
        \item A subgroup $H\preceq K$ is called a $\mathbb{P}$-virtual Ramsey subgroup of $K$ if for every $\mathbb{P}$-name $\dot{c}$ for a function from $K/H$ to $\{\mathrm{green,\ red}\}$, and for every $p\in \mathbb{P}$, there is a $q\leq p$ and a color $c'\in\{\mathrm{green,\ red}\}$ so that the following holds: if $\delta\in[K]^{<\omega}$ and the conditions $\{\pi q|\pi\in \delta\}$ are mutually compatible, there is a $\sigma\in K$ and a condition $r_\delta$ that is $\bigcap_{\pi\in \delta}(\sigma \pi)H(\sigma\pi)^{-1}$-reducible to $\bigcup_{\pi\in \delta}\sigma\pi q$, so that for every $\pi\in \delta$,
        $$
        r_\delta\Vdash\dot{c}(\check{\sigma\pi H})=\check{c}'.
        $$
        \item The symmetric system $\langle\mathbb{P},\mathcal{G},\mathcal{F}\rangle$ has the virtual Ramsey property if it has a base of subgroups $K$ for which every subgroup $H$ of $K$ in $\mathcal{F}$ is a $\mathbb{P}$-virtual Ramsey subgroup of $K$.
    \end{enumerate}
    
\end{definition}

\begin{remark}
    In definition \ref{Virtual Ramsey property definition}, we refer to the union of conditions. In the poset $\mathbb{P}(I,Q)$, this is well-defined. For an arbitrary poset, one can work in the regular open algebra of the poset, and take this union to mean the greatest lower bound of the given conditions.
\end{remark}

\begin{theorem}\label{RP, FEP, BPI relationship for symmetric systems}
    If $M\vDash\mathrm{ZFC}$, then for the conditions labeled below, $(i)\iff(ii)\Leftarrow(iii)$ and $(ii)\not\Rightarrow(iii)$.
    \begin{enumerate}[label=(\roman*)]
        \item $\mathrm{BPI}$ holds in the symmetric extension of $M$ based on $\langle\mathbb{P},\mathcal{G},\mathcal{F}\rangle$,
        \item $\langle\mathbb{P},\mathcal{G},\mathcal{F}\rangle$ has the filter extension property,
        \item $\langle\mathbb{P},\mathcal{G},\mathcal{F}\rangle$ has the virtual Ramsey property,
    \end{enumerate}
\end{theorem}

\begin{proof}
Theorem \ref{BPI equivalent to FEP for sym ext} is that $(i)\iff(ii)$, so we prove the relationship between $(ii)$ and $(iii)$.

\underline{$(iii)\Rightarrow(ii)$} 
Suppose that $K\in\mathcal{F}$ is a witness of the $\mathbb{P}$-virtual Ramsey property, that the $\mathbb{P}$-names $\dot{F},\dot{x},\dot{y}\in\mathrm{HS}_\mathcal{F}$, meet the hypothesis of the filter extension property, with $K\preceq\mathrm{sym}_\mathcal{G}(\dot{F})\cap\mathrm{sym}_\mathcal{G}(\dot{x})$, and let $p\in\mathbb{P}$. We aim to show that $K$ also witnesses the filter extension property by producing $q\leq p$ and $y'\in\{\dot{y},\dot{x}\setminus\dot{y}\}$, for which
\begin{equation}\label{FEP conclusion, section 6}
1\Vdash``\dot{F}\cup\mathrm{Orb}(\langle q,y'\rangle,K)\subseteq2^{\dot{x}}\mathrm{\ has\ the\ FIP}."
\end{equation}
If this can be done, then the base of subgroups of $\mathcal{F}$ that witnesses the virtual Ramsey property also witnesses the filter extension property.

Let $H=K\cap\mathrm{sym}_\mathcal{G}(\dot{y})=K\cap\mathrm{sym}_\mathcal{G}(\dot{x}\setminus\dot{y})\preceq K\preceq\mathrm{sym}_\mathcal{G}(\dot{F})\cap\mathrm{sym}_\mathcal{G}(\dot{x})$; note that $H\in\mathcal{F}$, so $H$ is a $\mathbb{P}$-virtual Ramsey subgroup of $K$. Fix a $\mathbb{P}$-name $\dot{U}$ for an ultrafilter on $\dot{x}$ that extends $\dot{F}$, and let $\dot{c}$ be a $\mathbb{P}$-name for a $2$-coloring of the cosets $K/H$ so that
$$
    1\Vdash(\dot{c}(\check{\pi H})=\check{\mathrm{green}}\iff \pi\dot{y}\in\dot{U})\wedge (\dot{c}(\check{\pi H})=\check{\mathrm{red}}\iff \pi(\dot{x}\setminus\dot{y})\in\dot{U}).
$$
Because $H=K\cap\mathrm{sym}_\mathcal{G}(\dot{y})=K\cap\mathrm{sym}_\mathcal{G}(\dot{x}\setminus\dot{y})$, $\dot{c}$ is well-defined.\footnote{Note that it may be the case that $\pi\dot{y}$ and $\tau\dot{y}$ are forced to be equal, even through $\pi H\neq\tau H$, and that this does not affect the proof.} Because $H$ is a $\mathbb{P}$-virtual Ramsey subgroup of $K$, let $q\leq p$ and $c'\in\{\mathrm{green,\ red}\}$ be as in definition \ref{Virtual Ramsey property definition}(i). We aim to show that $q,y'=\dot{y}$ satisfy line (\ref{FEP conclusion, section 6}) if $c'=\mathrm{green}$ and that $q,y'=\dot{x}\setminus\dot{y}$ satisfy line (\ref{FEP conclusion, section 6}) if $c'=\mathrm{red}$. Without loss of generality, we assume $c'=\mathrm{green}$ and prove the former case.

Take $\delta=\{\pi_0,...,\pi_n\}\in[K]^{<\omega}$ for which the conditions in $\{\pi q|\pi\in \delta\}$ are mutually compatible, so that the names $\pi_0\dot{y},...,\pi_n\dot{y}$ can mutually appear in an interpretation of $\mathrm{Orb}(\langle q,\dot{y}\rangle,K)$. Using definition \ref{Virtual Ramsey property definition}(i), fix $\sigma\in K$ and $r_\delta\in\mathbb{P}$ so that for every $i\leq n$,
    $$
    r_\delta\Vdash``\dot{c}(\check{\sigma\pi_i H})=\mathrm{green}."
    $$
By our definition for $\dot{c}$, this means that
$$
r_\delta\Vdash\bigcap_{\pi\in \delta}\sigma\pi \dot{y}\in\dot{U},
$$
and thus that
    \begin{equation}\label{line in U, section 6}
    r_\delta\Vdash``\dot{F}\cup\{\langle 1,\sigma\pi\dot{y}\rangle\ |\ \pi\in\delta\}\subseteq2^{\dot{x}}\mathrm{\ has\ the\ FIP}."
    \end{equation}

By definition \ref{Virtual Ramsey property definition}(i), we can take $r_\delta$ to be $\bigcap_{\pi\in \delta}(\sigma\pi)H(\sigma\pi)^{-1}$-reducible to $\bigcup_{\pi\in \delta}\sigma\pi q$. Because we took $H=K\cap\mathrm{sym}_\mathcal{G}(\dot{y})$, and because $\sigma\in K$ and $\delta\subseteq K$, for $\pi\in \delta$ we have that
    $$
    (\sigma\pi)H(\sigma\pi)^{-1}=K\cap\mathrm{sym}_\mathcal{G}(\sigma\pi \dot{y}),
    $$
and thus
$$
\bigcap_{\pi\in \delta}(\sigma\pi)H(\sigma\pi)^{-1}\preceq\mathrm{sym}_\mathcal{G}(\dot{F})\cap\mathrm{sym}_\mathcal{G}(\dot{x})\cap\bigcap_{\pi\in \delta}\mathrm{sym}_\mathcal{G}(\sigma\pi \dot{y}).
$$
By definition \ref{Reducible definition} and line (\ref{line in U, section 6}), we have that
$$
\bigcup_{\pi\in \delta}\sigma\pi q\Vdash``\dot{F}\cup\{\langle1,\sigma\pi\dot{y}\rangle\ |\ \pi\in\delta\}\subseteq2^{\dot{x}}\mathrm{\ has\ the\ FIP}."
$$
By applying the symmetry lemma with $\sigma^{-1}\in K\preceq\mathrm{sym}_\mathcal{G}(\dot{F})\cap\mathrm{sym}_\mathcal{G}(\dot{x})$, we get that
$$
\bigcup_{\pi\in \delta}\pi q\Vdash``\dot{F}\cup\{\langle1,\pi\dot{y}\rangle\ |\ \pi\in\delta\}\subseteq2^{\dot{x}}\mathrm{\ has\ the\ FIP}."
$$
Since $\delta\in [K]^{<\omega}$ was arbitrary, we can conclude line (\ref{FEP conclusion, section 6}) with $q\leq p$ and $y'=\dot{y}$, and thus that $K$ witnesses the filter extension property for $\langle\mathbb{P},\mathcal{G},\mathcal{F}\rangle$.

 \underline{$(ii)\not\Rightarrow(iii)$} 
We defer the full proof of this non-implication to the next subsection, where we show that the virtual Ramsey property fails for $\langle\mathbb{P}(\omega,2^{<\omega}),\mathcal{G},\mathcal{F}\rangle$. Conditional on this theorem, the non-implication is proven -- $\langle\mathbb{P}(\omega,2^{<\omega}),\mathcal{G},\mathcal{F}\rangle$ has the filter extension property by corollary \ref{N(I,Q) models BPI from FEP}.
\end{proof}

In section 5, we used to fact that distinct symmetric systems could be used to present the same symmetric extension. We can thus rephrase theorem \ref{RP, FEP, BPI relationship for symmetric systems} in terms of symmetric extensions.

\begin{theorem}\label{RP, FEP, BPI for symmetric extensions}
    If $M\vDash\mathrm{ZFC}$ and $N$ is a symmetric extension of $M$, then for the conditions labeled below, $(i)\iff(ii)\Leftarrow(iii)$.
    \begin{enumerate}[label=(\roman*)]
        \item $N\vDash\mathrm{BPI}$,
        \item every symmetric system presenting $N$ has the filter extension property,
        \item There is a symmetric system presenting $N$ that has the virtual Ramsey property.
    \end{enumerate}
\end{theorem}
\begin{proof}
  Let $\langle \mathbb{P},\mathcal{G},\mathcal{F}\rangle$ be a symmetric system presenting $N$ that has the virtual Ramsey property. By theorem \ref{RP, FEP, BPI relationship for symmetric systems}, $\langle\mathbb{P},\mathcal{G},\mathcal{F}\rangle$ has the filter extension property, and hence $N\vDash\mathrm{BPI}$; thus, $(iii)\Rightarrow(ii)\Rightarrow(i)$.
  
  Let $\mathbb{P}',\mathcal{G}',\mathcal{F}'\rangle$ be an arbitrary symmetric system that presents $N$, and assume that $N\vDash\mathrm{BPI}$. By theorem \ref{BPI equivalent to FEP for sym ext} (or by theorem \ref{RP, FEP, BPI relationship for symmetric systems}), $\langle\mathbb{P}',\mathcal{G}',\mathcal{F}\rangle$ has the filter extension property; thus, $(i)\Rightarrow(ii)$.
\end{proof}

\subsection{Revisiting $\mathrm{BPI}$ in the Generalized Cohen Model}
We can now verify that when $I$ is large, $\langle\mathbb{P}(I,Q),\mathcal{G},\mathcal{F}\rangle$ has the virtual Ramsey property, and that the virtual Ramsey property fails for $\langle\mathbb{P}(\omega,2^{<\omega}),\mathcal{G},\mathcal{F}\rangle$.

\begin{theorem}\label{fix(F cup X cup Y) a P(theta,Q)-virtual ramsey subgroup of fix(F cup X)}
    Let $\theta$ be large enough for theorem \ref{N(theta,Q) Models BPI}, and let $\langle\mathbb{P}(\theta,Q),\mathcal{G},\mathcal{F}\rangle$ be a symmetric system for a generalized Cohen model. If $X\subseteq Y\in[\theta]^{<\omega}$, then $\mathrm{fix}_\mathcal{G}(Y)$ is a $\mathbb{P}(\theta,Q)$-virtual Ramsey subgroup of $\mathrm{fix}_\mathcal{G}(X)$.
\end{theorem}
\begin{proof}
    Instead of repeating the full details of this proof, we will refer to the details of our proof of theorem \ref{N(theta,Q) Models BPI}. Fix a name for a coloring $\dot{c}$ as in the theorem statement, and $p\in\mathbb{P}(\theta,Q)$; for the same reason as in the proof of theorem \ref{N(theta,Q) Models BPI}, we can assume without loss of generality that $\mathrm{supp}(p)=X\cup Y$. We will treat the coloring name $\dot{c}$ in the same way that we treated the ultrafilter name $\dot{U}$ in our proof of theorem \ref{N(theta,Q) Models BPI}. Let $\dot{x},\dot{y}$ be names for which $\mathrm{sym}_\mathcal{G}(\dot{x})=\mathrm{fix}_\mathcal{G}(X)=K$ and $\mathrm{sym}_\mathcal{G}(\dot{y})=\mathrm{fix}_\mathcal{G}(Y)=H$.\footnote{Compared with the proof of theorem \ref{N(theta,Q) Models BPI}, this is akin to taking $\dot{F}$ to be the empty filter.} Then $\dot{c}$ can be replaced, without loss of generality, with the name $\dot{C}$ for the coloring of $\mathrm{Orb}(\dot{y},K)$ defined by
    $$
    1\Vdash (\dot{c}(\check{\pi H})=\mathrm{green}\iff\dot{C}(\pi\dot{y})=\mathrm{green})\wedge(\dot{c}(\check{\pi H})=\mathrm{red}\iff\dot{C}(\pi\dot{y})=\mathrm{red}).
    $$
    For $i\in Y\setminus X$, let $\theta_i\subseteq\theta\setminus X$ be disjoint sets of cardinality $\theta$; for every $\vec{\alpha}\in\prod_{i\in Y\setminus X}\theta_i$, either
    \begin{enumerate}
        \item $\pi_{\vec{\alpha}} p\Vdash\dot{C}(\pi\dot{y})=\mathrm{red}$, or
        \item there is some $q_{\vec{\alpha}}\leq p$ so that $\pi_{\vec{\alpha}}q_{\vec{\alpha}}\Vdash\dot{C}(\pi\dot{y})=\mathrm{green}$.
    \end{enumerate}
By taking $\theta$ as in theorem \ref{N(theta,Q) Models BPI}, for $i\in Y\setminus X$, there are disjoint infinite sets $\mathcal{H}_i$ so that either
\begin{enumerate}
    \item[(a)] the first alternative above happens for all $\vec{\alpha}\in \prod_{i\in Y\setminus X}\mathcal{H}_i$, or
    \item[(b)] there is a type $t$ and a condition $q$ with $\mathrm{supp}(q)= X\cup Y$ that distinguishes its support, so that for all $\vec{\alpha}\in \prod_{i\in Y\setminus X}\mathcal{H}_i$, $\pi_{\vec{\alpha}}q_{\vec{\alpha}}$ satisfies the second alternative above, is of type $t$, and restricts to $q_{\vec{\alpha}}|_{X\cup Y}=q$.
\end{enumerate}
We can assume without loss of generality that the second case occurs, and we argue that $q\leq p$ and $c'=\mathrm{green}$ witness, for the coloring $\dot{c}$, that $H$ is a $\mathbb{P}(\theta,Q)$-virtual Ramsey subgroup of $K$.

We can complete the proof by using more details of the proof of theorem \ref{N(theta,Q) Models BPI}; for every claim made below, the corresponding claim is proven in theorem \ref{N(theta,Q) Models BPI}. Let $\delta\in[K]^{<\omega}$, so that the conditions in the set $\{\pi q\ |\ \pi\in\delta\}$ are mutually compatible. For every $m<\omega$, there are disjoint sets $H^m_i\subseteq\mathcal{H}_i$, for $i\in Y\setminus X$, so that the conditions $\pi_{\vec{\alpha}}q_{\vec{\alpha}},$ for $\vec{\alpha}\in\prod_{i\in Y\setminus X}H^m_i$ are mutually compatible. We can produce, for some $S\subseteq\prod_{i\in Y\setminus X}H^m_i$, a $r_\delta=\bigcup_{\vec{\alpha}\in S}\pi_{\vec{\alpha}}q_{\vec{\alpha}}$ and a $\sigma\in K$ so that for every $\pi\in\delta$,
$$
r_\delta\Vdash\dot{C}(\sigma\pi\dot{y})=\mathrm{green}.
$$
We can then show that $r_\delta$ restricts to $\bigcup_{\pi\in\delta}\sigma\pi q$ over the union of supports $\bigcup_{\pi\in\delta}\sigma\pi(X\cup Y)$. By lemma \ref{reducibility lemma in N(I,Q)}, with $F=\emptyset,$ this shows that $r_\delta$ is $\bigcap_{\pi\in \delta}(\sigma\pi)H(\sigma\pi)^{-1}$-reducible to $\bigcup_{\pi\in\delta}\sigma\pi q$. This completes the proof for $\dot{C}$, and thus for the original coloring $\dot{c}$.
\end{proof}

We can now check that theorem \ref{fix(F cup X cup Y) a P(theta,Q)-virtual ramsey subgroup of fix(F cup X)} implies that every subgroup of $K=\mathrm{fix}_{\mathcal{G}}(E)$ in $\mathcal{F}$ is a $\mathbb{P}(\theta,Q)$-virtual Ramsey subgroup of $K$.

\begin{lemma}\label{vrp subgroup squeeze lemma}
    Let $\mathbb{P}$ is a poset and let $\mathcal{D},H,$ and $K$ are groups for which $\mathcal{D}\preceq H\preceq K\preceq\mathrm{Aut}(\mathbb{P})$. If $\mathcal{D}$ is a $\mathbb{P}$-virtual Ramsey subgroup of $K$, then $H$ is also a $\mathbb{P}$-virtual Ramsey subgroup of $K$.
\end{lemma}
\begin{proof}
    First note that every $\mathbb{P}$-name $\dot{c}$ for a function from $K/H$ to $\{\mathrm{green,\ red}\}$ is also (up to a minor change in notation) a $\mathbb{P}$-name for a function from $K/\mathcal{D}$ to $\{\mathrm{green,\ red}\}$. Next, note that because $\mathcal{D}\preceq H$, we have that for every $\delta\in[K]^{<\omega}$ and $\sigma\in K$,
    $$
    \bigcap_{\pi\in \delta}(\sigma\pi)\mathcal{D}(\sigma\pi)^{-1}\preceq \bigcap_{\pi\in \delta}(\sigma\pi)H(\sigma\pi)^{-1}.
    $$
    Then, by definition \ref{Reducible definition}, if $r$ is $\bigcap_{\pi\in \delta}(\sigma\pi)\mathcal{D}(\sigma\pi)^{-1}$-reducible to $\bigcup_{\pi\in \delta}\sigma\pi q$, it is also $\bigcap_{\pi\in \delta}(\sigma\pi)H(\sigma\pi)^{-1}$-reducible to $\bigcup_{\pi\in \delta}\sigma\pi q$. From these details, the lemma holds.
\end{proof}

\begin{corollary}If $\theta$ is large enough for theorem \ref{N(theta,Q) Models BPI}, then $\langle\mathbb{P}(\theta,Q),\mathcal{G},\mathcal{F}\rangle$ has the virtual Ramsey property.
\end{corollary}
\begin{proof}
We claim that the base $\{\mathrm{fix}_\mathcal{G}(X)\ |\ X\in[\theta]^{<\omega}\}$ of $\mathcal{F}$ witnesses the virtual Ramsey property. Let $K=\mathrm{fix}_\mathcal{G}(X)$, and let $H\in\mathcal{F}$ be a subgroup of $K$. We can find a $Y\supseteq X$ so that $\mathrm{fix}_\mathcal{G}(Y)\preceq H$. By theorem \ref{fix(F cup X cup Y) a P(theta,Q)-virtual ramsey subgroup of fix(F cup X)}, $\mathrm{fix}_\mathcal{G}(Y)$ is a $\mathbb{P}(\theta,Q)$-virtual Ramsey subgroup of $K$; by lemma \ref{vrp subgroup squeeze lemma}, so is $H$.
\end{proof}

Theorem \ref{FEP without RP} formally completes our proof of theorem \ref{RP, FEP, BPI for symmetric extensions}, by establishing that the filter extension property does not imply the virtual Ramsey property for a symmetric system.

\begin{theorem}\label{FEP without RP}
    $\langle\mathbb{P}(\omega,2^{<\omega}),\mathcal{G},\mathcal{F}\rangle$ does not have the virtual Ramsey property.
\end{theorem}

\begin{proof}   
We will show that there are names $\dot{F},\dot{x},\dot{y}\in\mathrm{HS}_\mathcal{F}$ that meet the hypotheses for the filter extension property, and a $\mathbb{P}(\omega,2^{<\omega})$-name $\dot{U}$ for an ultrafilter on $\dot{x}$ that extends $\dot{F}$, but cannot be used in a search and shift argument to produce witnesses $q\in\mathbb{P}(\omega,2^{<\omega})$ and $y'\in\{\dot{y},\dot{x}\setminus\dot{y}\}$ for the filter extension property.

To do so, take $\dot{F}=\emptyset,\dot{x}=\check{\omega}$ and $\dot{y}=\dot{a}_n$ for some $n\in\omega$. Recall from definition \ref{Dedekind Finite A} that $\dot{a}_n$ is the standard name for the Cohen real in $N(\omega,2^{<\omega})$ indexed by $n\in\omega$; note that when $n\not\in E\in[\omega]^{<\omega}$, $\mathrm{Orb}(\dot{a}_n,\mathrm{fix}_\mathcal{G}(E))=\{\dot{a}_m\ |\ m\in\omega\setminus E\}$. Let $\dot{F}^*$ be a name for a subset of $2^\omega$, so that for every $m\in\omega$,
    \begin{equation}
    (p\Vdash\dot{a}_m\in\dot{F}\iff p\Vdash\check{m}\in\dot{a}_m)\wedge(p\Vdash\check{\omega}\setminus\dot{a}_m\in\dot{F}\iff p\Vdash\check{m}\not\in\dot{a}_m).
    \end{equation}
Note that $\dot{F}^*\not\in\mathrm{HS}_\mathcal{F}$. By a standard density argument,
    $$
    1\Vdash``\dot{F}^*\subseteq2^{\check{\omega}}\mathrm{\ has\ the\ FIP}."
    $$
    Let $\dot{U}$ be any name for an ultrafilter on $\omega$ that extends $\dot{F}^*$. Because $\dot{F}^*$ decides between each $\dot{a}_m$ and its complement on a maximal antichain of conditions, note that for all $m\in\omega$
    $$
    1\Vdash(\dot{a}_m\in\dot{F}^*\iff\dot{a}_m\in\dot{U})\wedge(\check{\omega}\setminus\dot{a}_m\in\dot{F}^*\iff\check{\omega}\setminus\dot{a}_m\in\dot{U}).
    $$

    We claim that, given an arbitrary $K\in\mathcal{F}$, we can choose some $n\in\omega$ to specify $\dot{a}_n$, so that $\dot{U}$ cannot be used in a search and shift argument to produce witnesses $q\leq 1$ and $y'\in\{\dot{a}_n,\check{\omega}\setminus\dot{a}_n\}$ for
    $$
    1\Vdash``\mathrm{Orb}(\langle q,y'\rangle,K)\mathrm{\ has\ the\ FIP}."
    $$
    Formally, what we mean is that given an arbitrary $K\in\mathcal{F}$, we can choose some $n\in\omega$, and consider the name for a $2$-coloring $\dot{c}$ of the cosets $K/H$, with $H=K\cap\mathrm{fix}_\mathcal{G}(\{n\})$, defined by
    $$
    1\Vdash(\dot{c}(\check{\pi H})=\mathrm{green}\iff\pi\dot{a}_n\in\dot{U})\wedge(\dot{c}(\check{\pi H})=\mathrm{red}\iff\pi(\check{\omega}\setminus\dot{a}_n)\in\dot{U}).
    $$
    We then claim that no pair $q\leq 1$ and $c'\in\{\mathrm{green, red}\}$ witnesses that $H$ is a $\mathbb{P}(\omega,2^{<\omega})$-virtual Ramsey subgroup of $K$, for the coloring $\dot{c}$.
    
    Let $K\in\mathcal{F}$ be arbitrary, suppose for the sake of contradiction that $K$ witnesses the virtual Ramsey property for $\mathcal{F}$, and suppose that $\mathrm{fix}_\mathcal{G}(E)\preceq K$ for some $E\in[\omega]^{<\omega}$. Fix $n\in\omega\setminus E$, and let $q\in\mathbb{P}(\omega,2^{<\omega})$ and $c'\in\{\mathrm{green, red}\}$ arbitrary. We will show that $q$ and $c'$ do not witness that $H=K\cap\mathrm{fix}_\mathcal{G}(\{n\})$ is a $\mathbb{P}(\omega,2^{<\omega})$-virtual Ramsey subgroup of $K$, for the coloring $\dot{c}$ defined above.
    
    By the definition of $\dot{U}$, $q$ can only force whether $\dot{a}_m\in\dot{U}$, or equivalently whether $\check{m}\in\dot{a}_m$, if $\mathrm{ht}(q(m))\geq m$. For any $\delta\in[K]^{<\omega}$ and $\sigma\in K$, let
    \begin{equation}
    N=\mathrm{max}\{\mathrm{ht}(q(k))|k\in\mathrm{supp}(q)\}=\mathrm{max}\{\mathrm{ht}((\bigcup_{\pi\in f}\sigma\pi q)(k))|k\in\mathrm{supp}\bigcup_{\pi\in f}\sigma\pi q\}.
    \end{equation}
    Thus, there is an $N\in\omega$ so that for every $m>N$, $\sigma\in K$, and $\delta\in[K]^{<\omega}$, $\bigcup_{\pi\in \delta}\sigma\pi q$ does not force whether $\dot{a}_m\in\dot{U}$, or equivalently whether $\check{m}\in\dot{a}_m$.
    
    Let $m_0,...,m_{N+1}\not\in (E\cup\mathrm{supp}(q))$ be distinct, let $\pi_{m_k}=(n\ m_k)$, and let $\delta=\{\pi_{m_k}\ |\ k\leq N+1\}$. Because $n\not\in E$ and each $m_k\not\in E$, we have that $\delta\subseteq \mathrm{fix}_\mathcal{G}(E)\preceq K$; because each $m_k\not\in\mathrm{supp}(q)$, the conditions in $\{\pi q|\pi\in\delta\}$ are mutually compatible. Let $\sigma\in K$ be arbitrary. By the pidgeonhole principle, there is $k\leq N+1$ for which $\sigma (m_k)=\sigma\pi_{m_k}(n)>N$. From this, we can conclude that $\bigcup_{\pi\in \delta}\sigma\pi q$ cannot force whether $\check{\sigma (m_k)}\in\dot{a}_{\sigma (m_k)}$, or equivalently whether $\dot{a}_{\sigma(m_k)}\in\dot{U}$.

    If $H$ was a $\mathbb{P}(\omega,2^{<\omega})$-virtual Ramsey subgroup of $K$, then for the coloring $\dot{c}$, the witnesses $q$ and $c'$, and the set $\delta\in [K]^{<\omega}$, there would be a $\sigma\in K$ and a condition $r_\delta$ so that for every $\pi\in \delta,$
    $$
    r_\delta\Vdash\dot{c}(\check{\sigma\pi H})=c',
    $$
    so that $r_\delta$ is $\bigcap_{\pi\in \delta}(\sigma\pi)H(\sigma\pi)^{-1}$-reducible to $\bigcup_{\pi\in\delta}\sigma\pi q$. By our definition of $\dot{c}$ from $\dot{U}$, this would imply that $r_\delta$ forces whether or not $\sigma\pi_{m_k}\dot{a}_n=\dot{a}_{\sigma(m_k)}\in\dot{U}$, or equivalently whether $\check{\sigma(m_k)}\in\dot{a}_{\sigma(m_k)}$, for every $k\leq N+1$. By the above considerations, there is some $k\leq N+1$ so that $\bigcup_{\pi\in\delta}\sigma\pi q$ does not force whether $\check{\sigma(m_k)}\in\dot{a}_{\sigma(m_k)}$. Because $\bigcap_{\pi\in\delta}(\sigma\pi)H(\sigma\pi)^{-1}\preceq\mathrm{sym}_\mathcal{G}(\check{\sigma(m_k)})\cap\mathrm{sym}_\mathcal{G}(\dot{a}_{\sigma(m_k)})$, the fact that
    $$
    r_\delta\Vdash\check{\sigma(m_k)}\in\dot{a}_{\sigma(m_k)}\not\Rightarrow \bigcup_{\pi\in \delta}\sigma\pi q\Vdash\check{\sigma(m_k)}\in\dot{a}_{\sigma(m_k)}
    $$
    shows that $r_\delta$ is not $\bigcap_{\pi\in\delta}(\sigma\pi)H(\sigma\pi)^{-1}$-reducible to $\bigcup_{\pi\in\delta}\sigma\pi q$, which gives the desired contradiction to $H$ being a $\mathbb{P}(\omega,2^{<\omega})$-virtual Ramsey subgroup of $K$. Because $K\in\mathcal{F}$ was arbitrary, $\langle\mathbb{P}(\omega,2^{<\omega},\mathcal{G},\mathcal{F}\rangle$ does not have the virtual Ramsey property.
\end{proof}

This raises the question of the missing implication arrow in theorem \ref{RP, FEP, BPI for symmetric extensions}.

\begin{question}\label{RP and BPI Question}
    If $\mathrm{BPI}$ holds in a symmetric extension $N$, is there a symmetric system that presents $N$ and has the virtual Ramsey property?
\end{question}

Recall that in this paper, we take symmetric systems to be set-sized as opposed to proper class-sized. This is used critically to prove, in theorem \ref{fep equiv to bpi in permutation models}, that the filter extension property for a symmetric system is necessary for $\mathrm{BPI}$ to hold in the corresponding symmetric extension $N$. Our proof of theorem \ref{fep equiv to bpi in permutation models} used the fact that $N\vDash\mathrm{SVC}$, and theorem \ref{SVC iff AC in extension} shows that choice can be restored over $N$ by forcing with a set-size poset $\mathbb{P}$. We expect that there are simple class-sized symmetric extensions that model $\mathrm{BPI}$, but their corresponding class-sized symmetric systems have neither the filter extension property, nor the virtual Ramsey property.

We also expect that none of the known (set-sized) symmetric extensions that model $\mathrm{BPI}$ are counterexamples to question \ref{RP and BPI Question}. Every such model that the author is aware of seems to critically rely on mimicking aspects of the dynamical features of the Cohen model. We conjecture that question \ref{RP and BPI Question} can be answered negatively, that doing so will establish fundamentally new reasons for $\mathrm{BPI}$ to hold in symmetric extensions, and that this will result in new independence theorems over $\mathrm{ZF}+\mathrm{BPI}$.

\section{Future Work and Questions}
\subsection{Future Work}
In part 2 of this series, we plan to apply the ideas from this paper to provide simple new proofs of the independence of $\mathrm{AC}$ from $\mathrm{BPI}+\mathrm{DC}_{\kappa}$ (originally proven in \cite{pincus1977addingdep}) and of the transfer theorem for $\mathrm{BPI}$ from permutation models to symmetric extensions (originally proven in \cite{pincus1974independence}). We will also prove the following new results: for every symmetric extension used in the Jech-Sochor transfer theorem, every wellorderable set has a (nonprincipal) ultrafilter, and  $\mathrm{BPI}+\mathrm{DC}$ does not imply that successor cardinals are regular. Part 2 can be read independently of this paper.

We are also working to apply the tools developed in this paper to prove $\mathrm{BPI}$ in the conjectured model of $\mathrm{BPI}+\mathrm{DC}$ from \cite{pincus1977adding}. From this, we are working to positively resolve Pincus' conjectured version of the Halpern-L\"auchli theorem from the same paper. We intend this to be a case study to later answer the question from \cite{karagila2019iterated} regarding how $\mathrm{BPI}$ is preserved in symmetric iterations. It is noteworthy that recently, in \cite{holy2025orderinga} and \cite{holy2025orderingb}, adjusted versions of Pincus' model from \cite{pincus1977adding} were defined using symmetric iterations, to study the consistency of the ordering principle and $\mathrm{DC}/\mathrm{DC}_\kappa$.

\subsection{Questions}
For the reader's convenience, we restate questions \ref{Iso btw sym systems question}, \ref{Iso btw sym systems in M question}, and \ref{RP and BPI Question} that were posed throughout the paper.

\begin{question-non}{\bf \ref{Iso btw sym systems question}.}
    What is a good notion of isomorphism between symmetric systems?
\end{question-non}

\begin{question-non}{\bf \ref{Iso btw sym systems in M question}.}
    To what extent can question \ref{Iso btw sym systems question} be answered using maps in the given model $M\vDash\mathrm{ZFC}$, as opposed to maps in some outer model $V$ of $M$?
\end{question-non}

We give partial answers to these questions in the case of the generalized Cohen model $N(I,Q)$, and in theorems \ref{varphi_b in terms of Sigma, G_J}, \ref{KS new proof by outer model}, \ref{Sigma Elementary Embedding Theorem}, and \ref{Sigma* elementary embedding}.

\begin{question-non}{\bf \ref{RP and BPI Question}.}
    If $\mathrm{BPI}$ holds in a symmetric extension $N$, is there a symmetric system that presents $N$ and has the virtual Ramsey property?
\end{question-non}

In section 6, we conjectured that question \ref{RP and BPI Question} can be answered negatively. Note the connection between questions \ref{Iso btw sym systems question} and \ref{RP and BPI Question}, as a counterexample to question \ref{RP and BPI Question} would involve examining every symmetric system that presents the given symmetric extension.

\appendix

\section{The Halpern-L\"auchli Theorem}

We state the variant of the Halpern-L\"auchli theorem that is relevant to this paper and introduce the basic definitions required to do so. This can be found in many sources, including chapter 6 of \cite{todorchevich1995some} and \cite{stefanovic2023alternatives}.

\begin{definition}
    Given a partial order $T$, we say that $T$ is a \textit{tree} if for every $t\in T$, the set
    $$
    \{s\in T\ |\ s\leq_T t\}
    $$
    is wellordered by $\leq_T$. Each $t\in T$ is called a node of the tree. We will write $\leq$ to mean $\leq_T$.
\end{definition}

\begin{definition}
    Define the height $\mathrm{ht}(t)$ of a node $t\in T$ as the unique order type of the set $\{s\in T\ |\ s\leq t\}$ and the height $\mathrm{ht}(T)$ of the tree $T$ as $\mathrm{sup}\{\mathrm{ht}(t)\ |\ t\in T\}$.
\end{definition}

\begin{definition}
    Given $l\in\mathrm{ON}$, define the $l$th level of $T$ to be $T(l)=\{t\in T\ |\ \mathrm{ht}(t)=l\}$.
\end{definition}

\begin{definition}
    Given $t\in T$, let $T[t]=\{s\in T\ |\ s\leq t\vee t\leq s\}$ be a tree using the order from $T$.
\end{definition}

\begin{definition}
    Given a tree $T$ and a node $t\in T$, define the set of successors of $t$ in $T$ as
    $$
    s(t)=\{u\in T\ |\ t\leq u\wedge \mathrm{ht}(u)=\mathrm{ht}(t)+1\}.
    $$
\end{definition}

\begin{definition}
We say that a tree is \textit{finitely branching} if for every $t\in T$, $|s(t)|<\omega$.
\end{definition}

\begin{definition}
    Given trees $T_0,...,T_{d-1}$ of height $\omega$, the level product $\otimes_{i\leq d}T_i$ is the set 
    $$
    \bigcup_{n\in\omega}(T_0(n)\times...\times T_{d-1}(n))
    $$
    ordered by $(a_0,...,a_{d-1})\leq (b_0,...,b_{d-1})$ if and only if $b_i\leq_{T_i} a_i$ for each $i<d$.
\end{definition}

\begin{definition}
    Let $T$ be a tree, $l\leq m\leq n$ be levels of $T$, and $t\in T(l)$. A set $S\subseteq T(n)$ is called $(m,n)$-dense above $t$ if for every $u\in T(m)$ with $t\leq u$, there is some $v\in S$ with $v\leq u$.
\end{definition}

\begin{theorem}[Halpern-L\"auchli]
    Given $d,e<\omega$, a level product $T=\otimes_{i<d}T_i$ of finitely branching trees of height $\omega$, and a coloring $c:T\to e$, there is some node $\vec{t}=\langle t_0,...,t_{d-1}\rangle\in T$ so that for every $m<\omega$ there are infinitely many $n<\omega$ with sets $S_i\subseteq T_i(n)$ (for $i<d$) that are $(m,n)$-dense above $t_i$ in the tree $T_i$, and so that the restriction $c\upharpoonright _{S_0\times...\times S_{d-1}}$ is constant.
\end{theorem}

\section{A Sketch of Harrington's Proof of the Halpern-L\"auchli Theorem}
Harrington's proof of the Halpern-L\"auchli theorem was originally unpublished. We base this sketch on the complete exposition in chapter 6 of \cite{todorchevich1995some}, but we introduce notation that will be relevant to the proofs in this paper. The proof in \cite{todorchevich1995some} treats the case of the Halpern-L\"auchli theorem where the level product is of the form $T=\otimes_{i<d}T'$, where each component tree is the same tree $T'$, and the coloring $c:\otimes_{i<d}T'\to\{\mathrm{green,\ red}\}$ only has two colors. They note that the proof only changes in notation to consider an arbitrary level product $\otimes_{i<d}T_i$, and that a simple induction proof establishes the case where $e<\omega$ many colors are used from the case where $|e|=2$. This case of the Halpern-L\"auchli theorem, for colorings $c:\otimes_{i<d}T'\to\{\mathrm{green,\ red}\}$, is also the case that is most directly connected to the proofs in this paper. As such, our exposition will make the same simplifying assumptions as in \cite{todorchevich1995some}.

Fix $T=\otimes_{i<d}T'$ and a coloring $c:T\to\{\mathrm{green,\ red}\}$, as in the statement of theorem $\mathrm{A}.8$. We can reverse the order on the tree $T'$, then consider the poset $\mathbb{P}(I,T')$.\footnote{See definition \ref{P(I,Q)} for a full definition of the poset $\mathbb{P}(I,T)$.} The first central idea in Harrington's proof is to view the Halpern-L\"auchli theorem as a density condition on the forcing $\mathbb{P}(I,T')$, and to identify a property of the $\mathbb{P}(I,T')$-extension that implies this density condition.

\begin{definition}
    We say that $p\in\mathbb{P}(I,T')$ is a level condition if there is an $n\in\omega$, so that for every $k\in\mathrm{supp}(p),$ $\mathrm{ht}(p(k))=n$. For such an $n\in \omega$ and $p\in\mathbb{P}(I,T')$, we denote $\mathrm{ht}(p)=n$.
\end{definition}

\begin{definition}
    For $D\in[I]^d$, define the subset
    $$
    T_D=\{p\in\mathbb{P}(I,T)\ |\ \mathrm{supp}(p)=D\wedge p\mathrm{\ is\ a\ level\ condition.}\}
    $$
\end{definition}

$T_D$ is isomorphic to $T=\otimes_{i<d}T'$, so we can define a coloring $c_D$ of $T_D$ that is isomorphic to $c$.

\begin{definition}\label{T_D isomorphism}
    Fix an order preserving isomorphism $f:T_D\to T$ and define
    $$
    c_D:T_D\to\{\mathrm{green,\ red}\}
    $$
    by $c_D(q)=c(f(q))$.
\end{definition}

\begin{lemma}
    Let $f:T_D\to T$ be the isomorphism fixed in definition $\mathrm{B}.3$. Then $q\in T_D$ and $c'\in\{\mathrm{green,\ red}\}$ witness the Halpern-L\"auchli theorem for $T_D$ and $c_D$ if and only if $f(q)\in T$ and $c'$ witness the Halpern-L\"auchli theorem for $T$ and $c$.
\end{lemma}

We can now consider what property of the extension $M[G]$ corresponds to the density condition given by the Halpern-L\"aucli theorem for $T_D$ and $c_D$.

\begin{definition}\label{Harrington real definition}
    We call the $\mathbb{P}(I,T')$-name $\dot{y}$ a (name for a) \textit{Harrington real} (for $T=\otimes_{i<d}T'$ and $c$, witnessed by the support $D$) if there is some $D\in[I]^d$ and a $c'\in\{\mathrm{green,\ red}\}$ so that (for some choice of $f:T_D\to T$ in definition \ref{T_D isomorphism}),
    $$
    \dot{y}=\{\langle p,\check{n}\rangle\ |\ p\in T_D\wedge c_D(p)=c'\wedge \mathrm{ht}(p)=n\}.
    $$
\end{definition}

\begin{lemma}
    If the $\mathbb{P}(I,T')$ name $\dot{y}$ is a Harrington real, then
$1\Vdash\dot{y}\subseteq\check{\omega}.$
\end{lemma}

Consider the symmetric system $\langle\mathbb{P}(I,T'),\mathcal{G},\mathcal{F}\rangle$, as in definition \ref{Generalized Cohen System}, corresponding to the generalized Cohen model $N(I,T')$. With this framing, if $\dot{y}$ is a Harrington real witnessed by the support $D\in[I]^d$, then $\mathrm{supp}(\dot{y})=D$ as in definition \ref{Support Definition}. We can also consider the action of $\mathcal{G}$ on the Harrington reals. This will be convenient language to state theorem \ref{Harrington theorem}, and it will begin to suggest the connections that we develop in this paper.

\begin{lemma}
    If $\dot{y}$ is a Harrington real for $\otimes_{i<d}T_i$ and $c$, witnessed by the support $D\in [I]^d$, and if $\pi\in\mathcal{G}$, then $\pi\dot{y}$ is also a Harrington real for $\otimes_{i<d}T_i$ and $c$, witnessed by the support $\pi D\in[I]^d$.
\end{lemma}

In theorem \ref{Harrington theorem}, we take $\pi_{\vec{\alpha}}$ as in definition \ref{pi_vec{alpha}}.

\begin{theorem}\label{Harrington theorem}
    Let the $\mathbb{P}(I,T')$-name $\dot{y}$ be a Harrington real for $\otimes_{i<d}T'$ and $c$, witnessed by the support $D\in[I]^d$. Let $g(\mathrm{green})=\dot{y}$ and $g(\mathrm{red})=\check{\omega}\setminus\dot{y}$. For the conditions labeled below, $(i)\Rightarrow(ii)$.
    \begin{enumerate}[label=(\roman*)]
        \item $q\in T_D$ and $c'\in\{\mathrm{green,\ red}\}$ witness the following property: for every $m\in\omega$, there are disjoint sets $H^m_i\subseteq I$ of cardinality $m$, for $i\in D$, so that the conditions $\pi_{\vec{\alpha}}q\in\mathbb{P}(I,T')$, for $\vec{\alpha}\in\prod_{i\in D}H^m_i$, are mutually compatible and
        $$
        \bigcup_{\vec{\alpha}\in\prod_{i<d}H^m_i}\pi_{\vec{\alpha}}q\Vdash``\bigcap_{\vec{\alpha}\in\prod_{i<d}H^m_i}\pi_{\vec{\alpha}}g(c')\mathrm{\ is\ infinite},"
        $$
        \item $q\in T_D$ and $c'\in\{\mathrm{green,\ red}\}$ witness the Halpern-L\"auchli theorem for $T_D$ and $c_D$.
    \end{enumerate}
\end{theorem}
\begin{proof}
This is proven in \cite{todorchevich1995some}, modulo the fact that we have applied lemma \ref{Support Restriction Lemma} to normalize the conditions used in their exposition to be of the form $\pi_{\vec{\alpha}}q$. The details for a similar normalization are given explicitly in section 4 of this paper.
\end{proof}

The second central idea of Harrington's proof is to search through a $\mathbb{P}(I,T)$-name $\dot{U}$ for an ultrafilter on $\omega$ to identify, for every $m<\omega$, disjoint sets $H^m_i\subseteq I_i$, for $i<d$, a condition $q\in T_D$, and mutually compatible conditions $\pi_{\vec{\alpha}}q_{\vec{\alpha}}\leq \pi_{\vec{\alpha}}q$, for $\vec{\alpha}\in\prod_{i\in D}H^m_i$, so that (without loss of generality)
$$
\langle \pi_{\vec{\alpha}}q_{\vec{\alpha}},\pi_{\vec{\alpha}}C(c')\rangle\in\dot{U}.
$$
Because $1$ forces that $\dot{U}$ is an ultrafilter on $\omega$, this gives $(i)$ from theorem \ref{Harrington theorem} with $\bigcup_{\vec{\alpha}\in\prod_{i<d}H^m_i}\pi_{\vec{\alpha}}q_{\vec{\alpha}}$ in place of $\bigcup_{\vec{\alpha}\in\prod_{i<d}H^m_i}\pi _{\vec{\alpha}}q$. Using lemma \ref{Support Restriction Lemma} to restrict from $\bigcup_{\vec{\alpha}\in\prod_{i<d}H^m_i}\pi _{\vec{\alpha}}q_{\vec{\alpha}}$ to $\bigcup_{\vec{\alpha}\in\prod_{i<d}H^m_i}\pi_{\vec{\alpha}}q$ is the normalization that we mentioned above.

The third central idea in the proof is that, by taking the index set $I$ in $\mathbb{P}(I,T')$ to be large enough, a Ramsey-theoretic argument can guarantee that this search is successful.

As a historical note, one can compare Harrington's proof to a certain style of proof in recursion theory, where one constructs an object with reference to a certain forcing poset. In this type of argument, the forcing relation is used to inform which conditions to incorporate into the object, to guarantee certain desirable features. This can be thought of as conducting a search through a forcing poset, guided by the forcing relation.

\section{A Proof of Lemma 4.6 From Lemma 6.8 of \cite{todorchevich1995some}}
The poset $\mathcal{C}_\theta$ and the $\mathcal{C}_\theta$-name $\dot{r}_\alpha$, defined in \cite{todorchevich1995some}, can be compared to $\mathbb{P}(\theta,2^{<\omega})$, as in definition \ref{P(I,Q)}, and the $\alpha^{\mathrm{th}}$ Cohen real added by $\mathbb{P}(\theta,2^{<\omega})$, which we denoted $\dot{a}_\alpha$ in definition \ref{Dedekind Finite A}.

\begin{definition}[Todor\v cevi\'c, Farah, \cite{todorchevich1995some}]
    Let $\mathcal{C}_\theta$ denote the poset of finite partial functions $p:\theta\to2$, ordered by reverse-inclusion. For $\alpha<\theta,$ define the $\mathcal{C}_\theta$-name $\dot{r}_\alpha$ so that
    $$
    (p\Vdash\dot{r}_\alpha(m)=\check{n})\iff(\omega\alpha+m\in\mathrm{dom}(p)\wedge p(\omega\alpha+m)=n).
    $$
\end{definition}

In \cite{todorchevich1995some}, the term ``isomorphism type" or just ``type" is used. We replace this with ``TF-type" and so as not to conflict with our use of ``type" in definition 4.5.

\begin{definition}[Todor\v cevi\'c, Farah, \cite{todorchevich1995some}]
    Two conditions $p,q\in\mathcal{C}_\theta$ are said to be of the same \textit{TF-type} if there is an order-preserving map $\pi$ from $\mathrm{dom}(p)$ to $\mathrm{dom}(q)$, so that $p(\xi)=q(e(\xi))$. We use the term \textit{TF-type} to refer to an equivalence class in $\mathcal{C}_\theta$ under this relation.
\end{definition}

Note that the notion of TF-type is distinct from the notion of type in definition 4.5, given the difference in the underlying posets  $\mathcal{C}_\theta$ and $\mathbb{P}(\theta,Q)$, but that these capture the same idea.

\begin{lemma}[Lemma 6.8, \cite{todorchevich1995some}]\label{Lemma 6.8}
    For every TF-type $t$ of a condition in $\mathcal{C}_\theta$ and integers $m$ and $d$, there is an integer $M=M(t,m,d)$ such that for every set $\{p_{\vec{\alpha}}\ |\ \vec{\alpha}\in M^d\}$ of elements of $\mathcal{C}_\theta$ of TF-type $t$, there exist $H_i\subseteq M$ (for $i<d$) of size $m$, so that the conditions in $\{p_{\vec{\alpha}}\ |\ \vec{\alpha}\in \prod_{i<d}H_i\}$ are pairwise compatible.
\end{lemma}

\begin{corollary}
    Lemma 4.6 holds.
\end{corollary}
\begin{proof}
    Let $t$ be a type in $\mathbb{P}(\theta,Q)$ and notice that for every $p,r\in t$ (which is well-defined, since $t\subseteq\mathbb{P}(\theta,Q)$ is an equivalence class), $\mathrm{rng}(p)=\mathrm{rng}(r)\subseteq Q$. Let $\mathrm{rng}(t)\subseteq Q$ be the finite range of any condition in $t$. For every $q\in\mathrm{rng}(t)\subseteq Q$, fix some $t_q\in2^{\omega}$, so that for every $q,q'\in\mathrm{rng}(t)$,
    $$
    t_q\ ||\ t_{q'}\iff q\ ||\ q'.
    $$
    Define the map
    $$
    f:t\to \mathcal{C}_\theta
    $$
    by setting $f(p)$ to be the finite partial function
    $$
    f(p):\bigcup_{\alpha\in\mathrm{supp}(p)}[\omega \alpha,\omega(\alpha+1))\to 2,
    $$
    so that when $p(\alpha)=q$, $f(p)\upharpoonright_{[\omega\alpha,\omega\alpha+\mathrm{ht}(t_q)]}$ is isomorphic to $t_q\in2^{<\omega}$. This ensures that conditions $p,r\in t$ are compatible iff and only if $f(p),f(r)\in 2^{<\omega}$ are compatible in $2^{<\omega}$.  

    The map $f$ allows us to translate between lemma \ref{Lemma 6.8} for FT-types in $\mathcal{C}_\theta$ and lemma \ref{Compatible Type Ramsey Lemma} for types in $\mathbb{P}(\theta,Q)$, which concludes the proof.
    
\end{proof}

\printbibliography

\end{document}